\PassOptionsToPackage{linktocpage}{hyperref}
\documentclass[review]{elsarticle}
\usepackage[fleqn]{amsmath}
\usepackage{amssymb,latexsym,amsmath}

\usepackage{overpic}
\usepackage{color}
 \usepackage[hyperindex,breaklinks]{hyperref}
\usepackage{graphicx}
\usepackage{amsfonts}
\usepackage{subfigure}
\usepackage{epstopdf}
\usepackage{subfigure}
\usepackage{epsfig}
\graphicspath{{../Plots/}}

\newtheorem{thm}{Theorem}[section]
\newtheorem{defn}{Definition}[section]

\newtheorem{lem}{Lemma}[section]

\newtheorem{exmp}{Example}[section]
\numberwithin{equation}{section}

\begin{document}
\begin{frontmatter}
\title{A direct discontinuous
Galerkin method for fractional convection-diffusion and
Schr\"{o}dinger type equations}
\author[]{Tarek Aboelenen}
\ead{tarek.aboelenen@aun.edu.eg}




\address{Department of Mathematics, Assiut University, Assiut 71516, Egypt}
\begin{abstract}
    A direct discontinuous Galerkin   (DDG) finite element method  is   developed for  solving   fractional convection-diffusion and Schr\"{o}dinger type equations with a fractional Laplacian operator of order $\alpha$ $(1<\alpha<2)$. The  fractional operator of order $\alpha$ is expressed as a composite of second order derivative and a fractional integral of order $2-\alpha$. These problems     have been expressed as a system of parabolic   equation and low order integral equation. This allows us to apply the DDG method which is based on the direct weak formulation for solutions of fractional convection-diffusion and Schr\"{o}dinger type equations   in each computational cell, letting cells communicate via the numerical flux $(\partial_{x}u)^{*}$ only.   Moreover, we prove stability and optimal order of convergence $O(h^{N+1})$ for the general fractional convection-diffusion  and Schr\"{o}dinger problems where $h$, $N$ are the space step size and polynomial degree. The DDG method has the advantage of easier formulation and implementation as well as the high order accuracy. Finally, numerical experiments confirm the theoretical results of the method.




\end{abstract}
\end{frontmatter}
\section{Introduction}

Fractional calculus is a useful tool in various areas of physics and engineering \cite{mandelbrot1983fractal,tarasov2011fractional,herrmann2014fractional,uchaikin2013fractional}. Several examples of applications can be found in wide areas, such as fractals \cite{tarasov2011fractional}, kinetic theories of systems with chaotic dynamics \cite{moffatt2013topological,zaslavsky2002chaos,saichev1997fractional}, pseudochaotic dynamics \cite{zaslavsky2001weak}, anomalous transport \cite{metzler2004restaurant}, electrochemistry \cite{oldham2010fractional} and image processing \cite{cuesta2012image}, etc.\\
In this paper, we consider the  fractional convection-diffusion equation
\begin{equation}\label{25n}
\begin{split}
&\frac{\partial u}{\partial t}+\varepsilon(-\Delta)^{\frac{\alpha}{2}}u+\frac{\partial}{\partial x} f(u)=0, \quad x\in \mathbb{R},\,\,t\in(0,T],\\
&u(x,0) = u_{0}(x),\quad x\in \mathbb{R}.
\end{split}
\end{equation}
The generalized nonlinear fractional Schr\"{o}dinger  equation
\begin{equation}\label{sch1vnz}
\begin{split}
&i\frac{\partial u}{\partial t}- \varepsilon_{1}(-\Delta)^{\frac{\alpha}{2}}u+ \varepsilon_{2}f(|u|^{2})u=0,\quad x\in \mathbb{R},\,\,t\in(0,T],\\
&u(x,0) = u_{0}(x),\quad x\in \mathbb{R},
\end{split}
\end{equation}
and the strongly coupled nonlinear fractional Schr\"{o}dinger equations
\begin{equation}\label{sch2}
\begin{split}
&i\frac{\partial u}{\partial t}- \varepsilon_{1}(-\Delta)^{\frac{\alpha}{2}}u+ \varpi_{1}u+\varpi_{2} v+ \varepsilon_{2}f(|u|^{2},|v|^{2})u=0,\quad x\in \mathbb{R},\,\,t\in(0,T],\\
&i\frac{\partial v}{\partial t}- \varepsilon_{3}(-\Delta)^{\frac{\alpha}{2}}v+\varpi_{2}u+ \varpi_{1}v+\varepsilon_{4} g(|u|^{2},|v|^{2})v=0,\quad x\in \mathbb{R},\,\,t\in(0,T],\\
&u(x,0) = u_{0}(x),\quad x\in \mathbb{R},\\
&v(x,0) =  v_{0}(x),\quad x\in \mathbb{R},
\end{split}
\end{equation}
and homogeneous boundary conditions. $f(u)$ and $g(u)$ are arbitrary (smooth) nonlinear real functions and  $\varepsilon,\,\varepsilon_{i}$, $i=1,2,3,4$  are a real constants, $\varpi_{1}$ is normalized birefringence constant and $\varpi_{2}$ is the linear coupling parameter which accounts for the effects that arise from twisting and elliptic deformation of the fiber \cite{cai2010multisymplectic}. The  fractional Laplacian $-(-\Delta)^{\frac{\alpha}{2}}$, which can be defined using Fourier analysis as  \cite{el2006finite,muslih2010riesz, yang2010numerical}

$$-(-\Delta)^{\frac{\alpha}{2}}u(x,t)=\mathcal{F}^{-1}(|\xi|^{\alpha}\hat{u}(\xi,t))$$
where $\mathcal{F}$ is the Fourier transform.\\
 The equation \eqref{25n} is involved in many different physical problems, such as such as geomorphology \cite{fowler2002evolution,alfaro2012general}, overdriven detonations in gases \cite{alfaro2012general,clavin2002instabilities},
signal processing \cite{azerad2012simultaneous}, and anomalous diffusion in semiconductor growth \cite{shlesinger1995levy}. Numerical studies of  fractional convection diffusion equations  have
attracted a lot of interest in recent years.
 Several authors have proposed a variety of high-order finite difference
schemes for solving time-fractional convection-diffusion equations, for examples \cite{cui2014high,lin2007finite,wang2014high,zhai2014unconditionally}, and solving space-fractional convection-diffusion equations \cite{liu2007stability,chen2014second}.
 Furthermore, numerical methods for fractional diffusion problems and financial models with fractional Laplacian operators or Riesz fractional derivatives have been  studied in a number of papers, such
as \cite{deng2008finite,ding2012new,yang2010numerical,matache2005fast, cont2005finite}.\\
The equations \eqref{sch1vnz} and \eqref{sch2} arise in many physical fields, especially in
in fluid mechanics, nonlinear optics, solid state physics and plasma waves and for two interacting nonlinear packets in a dispersive and conservative system, see, e.g.,\cite{benney1967propagation, bullough1979solitons,yang1997classification,ran2016conservative} and reference therein. In recent years, developing various numerical algorithms for solving these problems has received much attention.  Wang and Huang \cite{wang2015energy} studied  an energy conservative Crank-Nicolson difference scheme for nonlinear Riesz space-fractional Schr\"{o}dinger equation. Yang \cite{yang2016class} proposed a class of linearized energy-conserved finite difference schemes for nonlinear space-fractional Schr\"{o}dinger equation. Galerkin finite element method for nonlinear fractional Schr\"{o}dinger equations were considered \cite{li2016galerkin}. Ran and Zhang \cite{ran2016conservative} proposed a conservative difference scheme for solving the strongly coupled nonlinear fractional Schr\"{o}dinger equations. A numerical study based on an implicit fully discrete LDG for the time-fractional coupled Schr\"{o}dinger systems is presented \cite{wei2012numerical}.\\
The discontinuous Galerkin (DG) method is a class of finite element methods using
discontinuous, piecewise polynomials as the solution and the test spaces in
the spatial direction. There have been various DG methods suggested in the literature to solve diffusion problem, including the method originally proposed by Bassi and Rebay \cite{bassi1997high} for compressible Navier-Stokes
equations, its generalization called the local discontinuous Galerkin (LDG) methods introduced
in \cite{cockburn1989tvb} by Cockburn and Shu and further studied in \cite{cockburn2002approximation, cockburn2005locally}.  For application of the method to fractional problems, Mustapha and McLean \cite{Mustapha2011,mustapha2012uniform,mustapha2013superconvergence,deng2013local}  have developed and analyzed discontinuous Galerkin methods for time fractional diffusion and wave equations.  Xu and Hesthaven \cite{doi:10.1137/130918174}  proposed a LDG method for  fractional convection-diffusion
equations. They proved stability and optimal order of convergence $N+1$
for the fractional diffusion problem when polynomials of degree $N$, and an order of convergence of $N+\frac{1}{2}$ is established for the general fractional convection-diffusion problem with general monotone flux for the nonlinear term. Aboelenen and El-Hawary \cite{cann} proposed a high-order nodal discontinuous Galerkin  method for  a linearized fractional  Cahn-Hilliard equation. They proved stability and optimal order of convergence $N+1$
for the linearized fractional Cahn-Hilliard  problem. A nodal discontinuous Galerkin  method  was   developed to  solve the nonlinear
 Riesz space fractional Schr\"{o}dinger equation and  the strongly coupled nonlinear Riesz space fractional Schr\"{o}dinger equations  \cite{Aboelenen2018428}. They proved, for both problems, $L^{2}$ stability  and optimal order of convergence $O(h^{N+1})$. Huang et al.\cite{huang2015numerical}   solved and analyzed  the time fractional diffusion equations by using a fully discrete DDG method.\\
The key of the local discontinuous Galerkin method for the fractional convection-diffusion and
Schr\"{o}dinger type equations \cite{doi:10.1137/130918174,Aboelenen2018428} is to rewrite the fractional operator as a composite of first order derivatives and a fractional integral and convert these problems into a system of low order equations by introducing an auxiliary variable. By solving the system, one obtains the solutions of the fractional convection-diffusion and
Schr\"{o}dinger type equations. The shortcoming is computational cost larger whether you use explicit and implicit step to solve the fully discrete system by the LDG method. Other DG method is called the DDG method introduced in  \cite{liu2009direct,liu2010direct, liu2015optimal} which is based on weak formulation for the solution of the parabolic equation in each computational cell and cells communicate via the
numerical flux $(\partial_{x}u)^{*}$ only. Here, we rewrite the fractional Laplacian operator as a
composite of second order derivative and a fractional integral and convert
 the fractional convection-diffusion and
Schr\"{o}dinger type equations into a system of parabolic equation and low order integral equation. This allows us to apply the DDG method which is based on the direct weak formulation of the equations \eqref{25n}-\eqref{sch2} and the construct of the suitable numerical flux on the cell edges. This method is called DDG method for not introducing any auxiliary variables  in contrast to the LDG method  \cite{doi:10.1137/130918174,Aboelenen2018428}.\\
The rest of this paper is organized as follows. Section \ref{sc1} introduce some basic definitions and recall a few central results. In section \ref{sc2}, we derive the DDG formulation for the  fractional convection-diffusion  problem. We present a stability and convergence analysis for the fractional convection-diffusion  equations in section \ref{sc3}.
We derive the DDG  formulation for the  nonlinear fractional Schr\"{o}dinger equation in section \ref{sc4}. Then  we prove a theoretical result of $L^{2}$ stability for the nonlinear case as well as an error estimate for the linear case in section \ref{sc5}. In section \ref{sc6} we present a DDG method for the strongly coupled nonlinear  fractional Schr\"{o}dinger equations and give a theoretical result of $L^{2}$ stability  for the nonlinear case and an error estimate for the linear case. Section \ref{s5n} presents some numerical examples to  illustrate the efficiency of the scheme. A few concluding remarks are offered in section \ref{s6n}.
\section{Preliminary definitions}\label{sc1}
In this section, we make some preparation including the definitions of fractional derivatives \cite{miller1993introduction} and associated functional setting for
the subsequent numerical schemes and theoretical analysis.\\
 Apart from the definitions of  the fractional
Laplacian based on the Fourier and the integral form, it can also be defined using
ideas of fractional calculus \cite{el2006finite,muslih2010riesz, yang2010numerical}, as
\begin{equation}\label{n5}
\frac{\partial^{\alpha}}{\partial |x|^{\alpha}}u(x,t)=-(-\Delta)^{\frac{\alpha}{2}}u(x,t)=-\frac{
{}^{\,\,\,\,\,\,C}_{-\infty}\mathcal{D}_{x}^{\alpha}u(x,t)+{}^{C}_{\,x}\mathcal{D}_{\infty}^{\alpha}u(x,t)}
{2\cos\big(\frac{\pi\alpha}{2}\big)},
\end{equation}
where ${}^{\,\,\,\,\,\,C}_{-\infty}\mathcal{D}_{x}^{\alpha}$ and ${}^{C}_{\,x}\mathcal{D}_{\infty}^{\alpha}$ refer to the left and right Caputo fractional derivatives, respectively, of $\alpha$th order. This definition is also known as a  Riesz derivative. To prepare
we introduce a few definitions and recall some properties of fractional integrals and
derivatives.\\
 The left-sided and right-sided Riemann-Liouville integrals of order $\alpha$, when $0 < \alpha < 1$, are defined, respectively, as
 \begin{equation}\label{1}
\big({}^{\,\,RL}_{-\infty}\mathcal{I}_{x}^{\alpha}f\big)(x)=\frac{1}{\Gamma(\alpha)}\int_{-\infty}^{x}
\frac{f(s)ds}{(x-s)^{1-\alpha}}, \quad x > -\infty,
\end{equation}
and
\begin{equation}\label{1111}
\big({}^{RL}_{\,\,\,\,x}\mathcal{I}_{\infty}^{\alpha}f\big)(x)=
\frac{1}{\Gamma(\alpha)}\int_{x}^{\infty}\frac{f(s)ds}{(s-x)^{1-\alpha}}, \quad x < \infty,
\end{equation}
where $\Gamma$ represents the Euler Gamma function. The corresponding inverse operators, i.e., the left-sided and
right-sided  fractional derivatives of order $\alpha$, are then defined based on \eqref{1} and \eqref{1111}, as
\begin{equation}\label{2}
\big({}^{\,\,RL}_{-\infty}\mathcal{D}_{x}^{\alpha}f\big)(x)=\frac{d}{dx}\big({}^{\,\,RL}_{-\infty}
\mathcal{I}_{x}^{1-\alpha}f\big)(x)=\frac{1}{\Gamma(1-\alpha)}
\frac{d}{dx}\int_{-\infty}^{x}\frac{f(s)ds}{(x-s)^{\alpha}}, \quad x > -\infty,
\end{equation}
and
\begin{equation}\label{3}
\big({}^{RL}_{\,\,\,\,x}\mathcal{D}_{\infty}^{\alpha}f\big)(x)
=\frac{-d}{dx}\big(^{RL}_{\,\,\,\,x}\mathcal{I}_{\infty}
^{1-\alpha}f\big)(x)=\frac{1}{\Gamma(1-\alpha)}
\bigg(\frac{-d}{dx}\bigg)\int_{x}^{\infty}\frac{f(s)ds}{(s-x)^{\alpha}}, \quad x < \infty.
\end{equation}
This allows for the definition of the left and right Riemann-Liouville fractional derivatives of order $\alpha$ $ (n-1 <\alpha<n),\,\,n\in \mathbb{N}$ as
\begin{equation}\label{2}
\big({}^{\,\,RL}_{-\infty}\mathcal{D}_{x}^{\alpha}f\big)(x)=\bigg(\frac{d}{dx}\bigg)^{n}\big({}^{\,\,RL}_{-\infty}
\mathcal{I}_{x}^{n-\alpha}f\big)(x)=\frac{1}{\Gamma(n-\alpha)}
\bigg(\frac{d}{dx}\bigg)^{n}\int_{-\infty}^{x}\frac{f(s)ds}{(x-s)^{-n+1+\alpha}}, \quad x > -\infty,
\end{equation}
and
\begin{equation}\label{3}
\big({}^{RL}_{\,\,\,\,x}\mathcal{D}_{\infty}^{\alpha}f\big)(x)=\bigg(\frac{-d}{dx}\bigg)^{n}\big(^{RL}_{\,\,\,\,x}
\mathcal{I}_{\infty}
^{n-\alpha}f\big)(x)=\frac{1}{\Gamma(n-\alpha)}
\bigg(\frac{-d}{dx}\bigg)^{n}\int_{x}^{\infty}\frac{f(s)ds}{(s-x)^{-n+1+\alpha}}, \quad x < \infty.
\end{equation}
Furthermore, the corresponding left-sided and right-sided  Caputo derivatives of order $\alpha$ $ (n-1 <\alpha<n)$ are obtained as
\begin{equation}\label{4n}
\big({}^{\,\,\,\,\,\,C}_{-\infty}\mathcal{D}_{x}^{\alpha}f\big)(x)=\bigg({}^{RL}_{-\infty}\mathcal{I}_{x}^{n-\alpha}
\frac{d^{n}f}{dx^{n}}\bigg)(x)
=\frac{1}
{\Gamma(n-\alpha)}\int_{-\infty}^{x}\frac{f^{(n)}(s)ds}{(x-s)^{-n+1+\alpha}}, \quad x > -\infty,
\end{equation}
and
\begin{equation}\label{5b}
\big({}^{C}_{\,x}\mathcal{D}_{\infty}^{\alpha}f\big)(x)=(-1)^{n}\bigg({}^{RL}_{\,\,\,\,x}
\mathcal{I}_{\infty}^{n-\alpha}\frac{d^{n}f}{dx^{n}}
\bigg)(x)=
\frac{1}{\Gamma(n-\alpha)}\int_{x}^{\infty}\frac{(-1)^{n}f^{(n)}(s)ds}{(s-x)^{-n+1+\alpha}}, \quad x < \infty.
\end{equation}
If $\alpha<0$, the fractional Laplacian becomes the fractional integral operator. In this case, for any $0 <\mu<1$, we define
\begin{equation}\label{nc5}
\Delta_{-\mu/2}u(x)=-\frac{
{}^{\,\,\,\,\,\,C}_{-\infty}\mathcal{D}_{x}^{-\mu}u(x)
+{}^{C}_{\,x}\mathcal{D}_{\infty}^{-\mu}u(x)}{2\cos\big(\frac{\pi(2-\mu)}{2}\big)}=\frac{
{}^{\,\,\,\,\,\,C}_{-\infty}\mathcal{D}_{x}^{-\mu}u(x)
+{}^{C}_{\,x}\mathcal{D}_{\infty}^{-\mu}u(x)}{2\cos\big(\frac{\pi\mu}{2}\big)}=\frac{
{}^{RL}_{-\infty}\mathcal{I}_{x}^{-\mu}u(x)
+{}^{RL}_{\,\,\,\,x}
\mathcal{I}_{\infty}^{-\mu}u(x)}{2\cos\big(\frac{\pi\mu}{2}\big)}.
\end{equation}
When $1 <\alpha<2$, using \eqref{4n}, \eqref{5b} and \eqref{nc5}, we can rewrite the fractional Laplacian in
the following form:
\begin{equation}\label{n5}
-(-\Delta)^{\frac{\alpha}{2}}u(x)=\Delta_{\frac{(\alpha-2)}{2}}\bigg(\frac{d^{2}u(x)}{dx^{2}}\bigg).
\end{equation}
To carry out the analysis, we introduce the appropriate fractional spaces.
 \begin{defn} (The right and left fractional spaces \cite{Ervin_variationalformulation}). We define the seminorm
\begin{equation}\label{7}
|u|_{J_{R}^{\alpha}(\mathbb{R})}=\big\|{}^{RL}_{\,\,\,\,x}
\mathcal{D}_{x_{R}}^{\alpha}u\big\|_{L^{2}(\mathbb{R})},
\end{equation}
\begin{equation}\label{7}
|u|_{J_{L}^{\alpha}(\mathbb{R})}=\big\|{}^{RL}_{\,x_{L}}\mathcal{D}_{x}^{\alpha}u\big\|_{L^{2}(\mathbb{R})}.
\end{equation}
and the norm
\begin{equation}\label{7}
\|u\|_{J_{R}^{\alpha}(\mathbb{R})}=(|u|_{J_{R}^{\alpha}(\mathbb{R})}^{2}+\|u\|_{L^{2}(\mathbb{R})}^{2})^{\frac{1}{2}},
\end{equation}
\begin{equation}\label{7}
\|u\|_{J_{L}^{\alpha}(\mathbb{R})}=(|u|_{J_{L}^{\alpha}(\mathbb{R})}^{2}+\|u\|_{L^{2}(\mathbb{R})}^{2})^{\frac{1}{2}},
\end{equation}
and let the two spaces $J_{R}^{\alpha}(\mathbb{R})$ and $J_{L}^{\alpha}(\mathbb{R})$ denote the closure of $C_{0}^{\infty}(\mathbb{R})$ with respect to $\|.\|_{J_{R}^{\alpha}(\mathbb{R})}$ and $\|.\|_{J_{L}^{\alpha}(\mathbb{R})}$ respectively.
\end{defn}
\begin{defn} (symmetric fractional space \cite{Ervin_variationalformulation}). We define the seminorm
\begin{equation}\label{7}
\|u\|_{J_{S}^{\alpha}(\mathbb{R})}=\big|\big({}^{RL}_{\,x_{L}}\mathcal{D}_{x}^{\alpha}u,{}^{RL}_{\,\,\,\,x}
\mathcal{D}_{x_{R}}^{\alpha}u\big)_{L^{2}(\mathbb{R})}\big|^{\frac{1}{2}},
\end{equation}
and the norm
\begin{equation}\label{7}
\|u\|_{J_{S}^{\alpha}(\mathbb{R})}=\big(|u|_{J_{S}^{\mu}(\mathbb{R})}^{2}+\|u\|_{L^{2}(\mathbb{R})}^{2}\big)^{\frac{1}{2}}.
\end{equation}
and let $J_{S}^{\alpha}(\mathbb{R})$ denote the closure of $C_{0}^{\infty}(\mathbb{R})$ with respect to $\|.\|_{J_{S}^{\alpha}(\mathbb{R})}$.
\end{defn}
\begin{lem}\label{lk}(see \cite{Ervin_variationalformulation}). For any $0 <s<1$, the fractional integral satisfies the following
property:
\begin{equation}\label{7}
({}^{\,\,RL}_{-\infty}\mathcal{I}_{x}^{s}u,{}^{RL}_{\,\,\,\,x}\mathcal{I}_{\infty}^{s}u)_{\mathbb{R}}=
\cos(s\pi)|u|_{J_{L}^{-s}(\mathbb{R})}^{2}=
\cos(s\pi)|u|_{J_{R}^{-s}(\mathbb{R})}^{2}.
\end{equation}
\end{lem}
\begin{lem}\label{lkb} For any $0 <\mu<1$, the fractional integral satisfies the following
property:
\begin{equation}\label{7}
(\Delta_{-\mu}u,u)_{\mathbb{R}}=
|u|_{J_{L}^{-\mu}(\mathbb{R})}^{2}=
|u|_{J_{R}^{-\mu}(\mathbb{R})}^{2}.
\end{equation}
\end{lem}
Generally, we consider the problem in a bounded domain instead of $\mathbb{R}$. Hence,
we restrict the definition to the domain $\Omega = [a, b]$.
\begin{defn} Define the spaces $J_{R,0}^{\alpha}(\Omega),J_{L,0}^{\alpha}(\Omega),J_{S,0}^{\alpha}(\Omega)$ as the closures of
$C_{0}^{\infty}(\Omega)$ under their respective norms.
\end{defn}
\begin{lem}\label{lg}
(fractional Poincar$\acute{e}$-Friedrichs, \cite{Ervin_variationalformulation}). For $u \in J_{L,0}^{\alpha}(\Omega)$ and $\alpha \in \mathbb{R}$, we have
\begin{equation}\label{7}
\|u\|_{L^{2}(\Omega)}\leq C |u|_{J_{L,0}^{\alpha}(\Omega)},
\end{equation}
and for $u \in J_{R,0}^{\alpha}(\Omega)$, we have
\begin{equation}\label{7}
\|u\|_{L^{2}(\Omega)}\leq C |u|_{J_{R,0}^{\alpha}(\Omega)}.
\end{equation}
\end{lem}
\begin{lem}\label{lga234} (See \cite{Kilbas:2006:TAF:1137742})
For any $0 <\mu<1$, the fractional integration operator ${}^{RL}_{\,\,\,\,a}\mathcal{I}_{x}^{\mu}$ is bounded in $L^{2}(\Omega)$:
\begin{equation}\label{7}
\|{}^{RL}_{\,\,\,\,a}\mathcal{I}_{x}^{\mu}u\|_{L^{2}(\Omega)}\leq K \|u\|_{L^{2}(\Omega)}.
\end{equation}
The fractional integration operator ${}^{RL}_{\,\,\,\,x}\mathcal{I}_{b}^{\mu}$ is bounded in $L^{2}(\Omega)$:
\begin{equation}\label{7}
\|{}^{RL}_{\,\,\,\,x}\mathcal{I}_{b}^{\mu}u\|_{L^{2}(\Omega)}\leq K \|u\|_{L^{2}(\Omega)}.
\end{equation}
\end{lem}
\begin{lem}\label{lga2} (See \cite{Aboelenen2018428})
The fractional integration operator $\Delta_{-\mu}$ is bounded in $L^{2}(\Omega)$:
\begin{equation}\label{7}
\|\Delta_{-\mu}u\|_{L^{2}(\Omega)}\leq K \|u\|_{L^{2}(\Omega)}.
\end{equation}
\end{lem}
\section{ The DDG scheme for  the fractional convection-diffusion equation}\label{sc2}
In this section, we  construct DDG method for the  fractional convection-diffusion problem \eqref{25n}. So, we rewrite
the Riesz  fractional derivative of order $\alpha$ $(1 < \alpha< 2)$  as a composite of second order derivative and low order fractional integral
  operator. Since the integral operator naturally connects the discontinuous function, we need not add a
penalty term or introduce a numerical fluxes for the integral equation.\\
 We introduce  the auxiliary variables $p,q$ and set
\begin{equation}\label{1a}
\begin{split}
&p=\Delta_{(\alpha-2)/2}q, \quad q=\frac{\partial^{2}}{\partial x^{2}}u,
\end{split}
\end{equation}
then, the nonlinear fractional convection-diffusion  problem can be rewritten as
\begin{equation}\label{1bchcf}
\begin{split}
&\frac{\partial u}{\partial t}-\varepsilon p+ \frac{\partial}{\partial x} f(u)=0,\\
&p=\Delta_{(\alpha-2)/2}q,\quad q=\frac{\partial^{2}}{\partial x^{2}}u.\\
\end{split}
\end{equation}
The weak solution formulation for this problem is to find a functions $u,p,q \in C(0, T; H^{1}(\Omega))$ such that for all $v,\psi,\phi\in  H_{0}^{1}(\Omega)$
 \begin{equation}\label{1bch2znmm}
\begin{split}
&\big( \frac{\partial u}{\partial t},v\big)-\varepsilon\big(p,v\big)-\big( f(u),\partial_{x} v\big)=0,\\
&\big(p,\psi\big)=\big(\Delta_{(\alpha-2)/2}q,\psi\big),\\
&\big(q,\phi\big)=-\big(\partial_{x}u,\partial_{x}\phi\big).\\
\end{split}
\end{equation}
 To discretize this weak formulation, we set up a partition of the domain $\Omega$ into $K$ non-overlapping elements such that $\Omega=\bigcup_{k=1}^{K}D^{k}$ with  mesh  $D^{k}=[x_{k-\frac{1}{2}},x_{k+\frac{1}{2}}]$,  $\Delta x_{k}=x_{k+\frac{1}{2}}-x_{k-\frac{1}{2}}$ and $k = 1,...,K$. Let $u_{h}, p_{h}, q_{h}\in V_{k}^{N}$ be the approximation of $u,p,q$ respectively, where  the approximation space  is defined as
\begin{equation}\label{82aa1}
V_{k}^{N}=\{v:v_{k}\in\mathbb{P}^{N}(D^{k}), \, \forall D^{k}\in \Omega\},
\end{equation}
where $\mathbb{P}^{N}(D^{k})$ denotes the set of polynomials of degree up to $N$  defined  on  the  element  $D^{k}$.\\
We define the local inner product and $L^{2}(D^{k})$ norm
\begin{equation}\label{82aa1}
(u,v)_{D^{k}}=\int_{D^{k}}uvdx,\quad \|u\|^{2}_{D^{k}}=(u,u)_{D^{k}},
\end{equation}
\begin{equation}\label{82aa1}
(u,v)=\sum_{k=1}^{K}(u,v)_{D^{k}},\quad \|u\|^{2}_{L^{2}(\Omega)}=\sum_{k=1}^{K}(u,u)_{D^{k}}.
\end{equation}
 We define  DDG scheme as follows: find $u_{h},p_{h}, q_{h}\in V_{k}^{N}$, such that for all test functions $v,\psi,\phi\in V_{k}^{N}$,
\begin{equation}\label{1bch1nb}
\begin{split}
&\big(\frac{\partial u_{h}}{\partial t},v\big)_{D^{k}}-\varepsilon\big(p_{h},v\big)_{D^{k}}+ \big(\frac{\partial}{\partial x}f(u_{h}),v\big)_{D^{k}}=0,\\
&\big(p_{h},\psi\big)_{D^{k}}=\big(\Delta_{(\alpha-2)/2}q_{h},\psi\big)_{D^{k}},\\
&\big(q_{h},\phi\big)_{D^{k}}=\big(\frac{\partial^{2}}{\partial x^{2}}u_{h},\phi\big)_{D^{k}}.\\
\end{split}
\end{equation}
To complete the DDG scheme, we introduce some notation and the numerical flux.\\
Define
\begin{equation}\label{82aa1}
\{u\}=\frac{u^{+}+u^{-}}{2},\quad [u]=u^{+}-u^{-}.
\end{equation}

 The numerical flux involves the average $\{\partial_{x}u\}$ and the jumps of even order derivatives of $u$ \cite{liu2010direct}:
\begin{equation}\label{flf}
\begin{split}
(\partial_{x}u)^{*}=\frac{\beta_{0}}{h}[u]+\{\partial_{x}u\}+\beta_{1}h[\partial_{x}^{2}u].\\
\end{split}
\end{equation}
where  $\beta_{0}$ and $\beta_{1}$ are chosen to ensure the stability and accuracy of the scheme.\\
The idea of the DDG method  is to enforce the weak
formulation \eqref{1bch2znmm} in such a way that both $u_{h}$ and $\phi$ are approximated in $V_{k}^{N}$. The discontinuous
nature of numerical solutions and test functions crossing interfaces necessarily requires some
interface corrections, leading to the following:
\begin{equation}\label{1bch2znm}
\begin{split}
&\big( (u_{h})_{t},v\big)_{D^{k}}-\varepsilon\big(p_{h},v\big)_{D^{k}}-\big( f(u_{h}),\partial_{x} v\big)_{D^{k}}+\big(\hat{n}.f(u_{h})^{*},v\big)_{\partial D^{k}}=0,\\
&\big(p_{h},\psi\big)_{D^{k}}=\big(\Delta_{(\alpha-2)/2}q_{h},\psi\big)_{D^{k}},\\
&\big(q_{h},\phi\big)_{D^{k}}=-\big(\partial_{x}u_{h},\partial_{x}\phi\big)_{D^{k}}
-\bigg((\partial_{x}u_{h})^{*}[\phi]+\{\partial_{x}\phi\}[u_{h}]\bigg)_{k+\frac{1}{2}}.\\
\end{split}
\end{equation}

\section{ Stability and error estimates}\label{sc3}
 In the following we discuss stability and accuracy of the proposed scheme, for the fractional diffusion  problem and the nonlinear fractional convection-diffusion  problem.
\subsection{ Stability analysis }  In order to carry out the analysis of the DDG scheme, we have the following results.
\begin{defn}\label{tt4hzx} (Admissibility \cite{liu2010direct}). We call a numerical flux $(\partial_{x}u_{h})^{*}$ of the form \eqref{flf} admissible if there exists a $\gamma \in (0, 1)$ and $0 < \mu \leq 1$ such that
\begin{equation}\label{1bch2hj}
\begin{split}
\gamma\big(\partial_{x}u_{h},\partial_{x}u_{h}\big)+
\sum_{k=1}^{K}((\partial_{x}u_{h})^{*}[u_{h}]
+\{\partial_{x}u_{h}\}[u_{h}]\big)_{k+\frac{1}{2}}\geq \mu\sum_{k=1}^{K}\frac{[u_{h}]^{2}_{k+\frac{1}{2}}}{h}.\\
\end{split}
\end{equation}
holds for  any piecewise  polynomials  of degree $N$,i.e., $u\in V_{k}^{N}$
\end{defn}
\begin{thm}\label{tt4hc} ($L^{2}$ stability).
The semidiscrete scheme \eqref{1bch2znm} is   $L^{2}$ stable, and for all $T>0$ its solution satisfies $$
\|u_{h}(x,T)\|_{L^{2}(\Omega)}\leq c\|u_{0}(x)\|_{L^{2}(\Omega)}. $$
\end{thm}
\textbf{Proof.}
Set $(v,\psi,\phi)=(u_{h},p_{h}-q_{h},u_{h})$ in \eqref{1bch2znm}, we get
\begin{equation}\label{1bch2z}
\begin{split}
&\big( (u_{h})_{t},u_{h}\big)_{D^{k}}-\big( f(u_{h}),\partial_{x} u_{h}\big)_{D^{k}}+\big(\hat{n}.f(u_{h})^{*},u_{h}\big)_{\partial D^{k}}+\big(p_{h},p_{h}\big)_{D^{k}}+\big(\Delta_{(\alpha-2)/2}q_{h},q_{h}\big)_{D^{k}}\\
&=
\big(p_{h},q_{h}\big)_{D^{k}}+\varepsilon\big(p_{h},u_{h}\big)_{D^{k}}-\big(q_{h},u_{h}\big)_{D^{k}}
-\big(\partial_{x}u_{h},\partial_{x}u_{h}\big)_{D^{k}}+\big(\Delta_{(\alpha-2)/2}q_{h},p_{h}\big)_{D^{k}}
-\bigg((\partial_{x}u_{h})^{*}[u_{h}]+\{\partial_{x}u_{h}\}[u_{h}]\bigg)_{k+\frac{1}{2}}.\\
\end{split}
\end{equation}
Define $\theta(u)=\int^{u}f(s)ds$; then
\begin{equation}\label{1bch2z}
\begin{split}
&\big( (u_{h})_{t},u_{h}\big)_{D^{k}}+\big(\hat{n}.(f(u_{h})^{*},u_{h}\big)_{\partial D^{k}}-(\theta(u_{h}))^{-}_{k+\frac{1}{2}}
+(\theta(u_{h}))^{+}_{k-\frac{1}{2}}+\big(p_{h},p_{h}\big)_{D^{k}}+\big(\Delta_{(\alpha-2)/2}q_{h},q_{h}\big)_{D^{k}}\\
&=
\big(p_{h},q_{h}\big)_{D^{k}}+\varepsilon\big(p_{h},u_{h}\big)_{D^{k}}-\big(q_{h},u_{h}\big)_{D^{k}}
-\big(\partial_{x}u_{h},\partial_{x}u_{h}\big)_{D^{k}}+\big(\Delta_{(\alpha-2)/2}q_{h},p_{h}\big)_{D^{k}}
-\bigg((\partial_{x}u_{h})^{*}[u_{h}]+\{\partial_{x}u_{h}\}[u_{h}]\bigg)_{k+\frac{1}{2}}.\\
\end{split}
\end{equation}
Employing Young's inequality and  Lemma \ref{lga2}, we obtain
\begin{equation}\label{1bch2z}
\begin{split}
&\big( (u_{h})_{t},u_{h}\big)_{D^{k}}+\big(\hat{n}.(f(u_{h})^{*},u_{h}\big)_{\partial D^{k}}-(\theta(u_{h}))^{-}_{k+\frac{1}{2}}
+(\theta(u_{h}))^{+}_{k-\frac{1}{2}}+\big(p_{h},p_{h}\big)_{D^{k}}+\big(\Delta_{(\alpha-2)/2}q_{h},q_{h}\big)_{D^{k}}\\
&\leq c_{3}\|u_{h}\|^{2}_{L^{2}(D^{k})}+c_{1}\|p_{h}\|^{2}_{L^{2}(D^{k})}+c_{2}\|q_{h}\|^{2}_{L^{2}(D^{k})}
-\big(\partial_{x}u_{h},\partial_{x}u_{h}\big)_{D^{k}}-\bigg((\partial_{x}u_{h})^{*}[u_{h}]+\{\partial_{x}u_{h}\}[u_{h}]\bigg)_{k+\frac{1}{2}}.\\
\end{split}
\end{equation}
Recalling Lemma \ref{lg}, provided $c_{1}$ is sufficiently small such that $c_{1}\leq1$, we obtain that
\begin{equation}\label{1bch2}
\begin{split}
&\big( (u_{h})_{t},u_{h}\big)_{D^{k}}+\big(\hat{n}.(f(u_{h})^{*},u_{h})_{\partial D^{k}}-(\theta(u_{h}))^{-}_{k+\frac{1}{2}}
+(\theta(u_{h}))^{+}_{k-\frac{1}{2}}\\
&\leq c_{3}\|u_{h}\|^{2}_{L^{2}(D^{k})}
-\big(\partial_{x}u_{h},\partial_{x}u_{h}\big)_{D^{k}}-\bigg((\partial_{x}u_{h})^{*}[u_{h}]+\{\partial_{x}u_{h}\}[u_{h}]\bigg)_{k+\frac{1}{2}},\\
\end{split}
\end{equation}
such that $\big(\hat{n}.f(u_{h})^{*},u_{h}\big)_{\partial D^{k}}-(\theta(u_{h}))^{-}_{k+\frac{1}{2}}
+(\theta(u_{h}))^{+}_{k-\frac{1}{2}}\geq 0$, summering over all elements, we immediately recover
\begin{equation}\label{1bch2}
\begin{split}
&\big( (u_{h})_{t},u_{h}\big)\leq c_{3}\|u_{h}\|^{2}_{L^{2}(\Omega)}
-\big(\partial_{x}u_{h},\partial_{x}u_{h}\big)-
\sum_{k=1}^{K}\bigg((\partial_{x}u_{h})^{*}[u_{h}]+\{\partial_{x}u_{h}\}[u_{h}]\bigg)_{k+\frac{1}{2}}.\\
\end{split}
\end{equation}
From  the admissible  condition   \ref{tt4hzx}  of the numerical  flux  defined  in  \eqref{flf}, we obtain that
\begin{equation}\label{1bch2}
\begin{split}
&\big( (u_{h})_{t},u_{h}\big)+(1-\gamma)\|\partial_{x}u_{h}\|^{2}_{L^{2}(\Omega)}
+\mu\sum_{k=1}^{K}\frac{[u_{h}]^{2}_{k+\frac{1}{2}}}{h}\leq c_{3}\|u_{h}\|^{2}_{L^{2}(\Omega)}.\\
\end{split}
\end{equation}

Employing Gronwall's lemma, we obtain $\|u_{h}(x,T)\|_{L^{2}(\Omega)}\leq c\|u_{0}(x)\|_{L^{2}(\Omega)}$. $\quad\Box$

\subsection{Error estimates}
we list some inverse properties and special projections  $\mathcal{P}^{+}$ of the finite element space $V_{k}^{N}$ that will be used in our error analysis.
\begin{equation}\label{prh}
\begin{split}
&(\mathcal{P}^{+}u-u,v)_{D^{k}}=0,\quad\forall  v\in\mathbb{P}_{N}^{k-2}(D^{k}),\quad k=1,...,K,\\
&(\partial_{x}\mathcal{P}^{-}u)^{*}=\big(\beta_{0}h^{-1}[\mathcal{P}^{-}u]+\{\partial_{x}\mathcal{P}^{-}u\}
+\beta_{1}h[\partial_{x}^{2}\mathcal{P}^{-}u]\big)_{x_{k+\frac{1}{2}}}=(\partial_{x}u(x_{k+\frac{1}{2}}))^{*}.
\end{split}
\end{equation}
\begin{lem}\label{lga2zz} (See \cite{liu2015optimal})
For $(\beta_{0}, \beta_{1})$ such that $\beta_{0} > \Gamma(\beta_{1})$, the projection $\mathcal{P}^{+}$ defined in \eqref{prh} exists, and
\begin{equation}\label{sth}
\begin{split}
&\|\mathcal{P}^{+}u(.)-u(.)\|_{L^{2}(\Omega)}\leq  Ch^{N+1}.\\
\end{split}
\end{equation}
\end{lem}
 Let $\pi^{e}$ denote the projection error, then  the following inequality  \cite{liu2015optimal}:
\begin{equation}\label{sth11}
\begin{split}
&\|\pi^{e}\|_{L^{2}(\Omega)}+h\|\pi^{e}\|_{\infty}+h^{\frac{1}{2}}\| \pi\|_{\Gamma_{h}}\leq  Ch^{N+1}.\\
\end{split}
\end{equation}
For any function $u_{h}\in V_{k}^{N}$, the following inverse inequalities hold \cite{Ciarlet:2002:FEM:581834}:
\begin{equation}\label{sth22}
\begin{split}
&\|\partial_{x} u_{h}\|_{L^{2}(\Omega)}\leq Ch^{-1}\|u_{h}\|_{L^{2}(\Omega_{h})},\\
&\| u_{h}\|_{\Gamma_{h}}\leq Ch^{-1/2}\|u_{h}\|_{L^{2}(\Omega)},\\
&\| u_{h}\|_{\infty}\leq Ch^{-1/2}\|u_{h}\|_{L^{2}(\Omega)},\\
\end{split}
\end{equation}
where $\Gamma_{h}$ denotes the set of interface points of all  the elements, $D^{k}$, $k = 1, 2, ... , K$ and here and below $C$ is a positive constant (which may have a different value in each
occurrence) depending solely on u and its derivatives but not of  $h$.
\begin{thm}\label{th1} (Diffusion without convection $f(u)=0$).
Let $u$ be the exact solutions of the problem \eqref{25n}, and let $u_{h}$ be the numerical solutions of the semi-discrete DDG scheme \eqref{1bch2znm}.  Then for small enough $h$, we have the following error estimates:
\begin{equation}\label{bnkl}
\begin{split}
&\|u-u_{h}\|_{L^{2}(\Omega)}\leq Ch^{N+1},\\
\end{split}
\end{equation}
\end{thm}
where the constant $C$ is dependent upon $T$ and some norms of the solutions. \\
\textbf{Proof}.
We consider the fractional diffusion equation
\begin{equation}\label{25nl}
\begin{split}
&\frac{\partial u}{\partial t}+\varepsilon(-\Delta)^{\frac{\alpha}{2}}u=0.
\end{split}
\end{equation}
It is easy to verify that the exact solution of the above  \eqref{25nl} satisfies
\begin{equation}\label{1bch2kj}
\begin{split}
\big( u_{t},v\big)_{D^{k}}-\varepsilon\big(p,v\big)_{D^{k}}+
\big(p,\psi\big)_{D^{k}}-\big(\Delta_{(\alpha-2)/2}q,\psi\big)_{D^{k}}+\big(q,\phi\big)_{D^{k}}
+\big(\partial_{x}u,\partial_{x}\phi\big)_{D^{k}}+\bigg((\partial_{x}u_{h})^{*}[\phi]+\{\partial_{x}\phi\}[u_{h}]\bigg)_{k+\frac{1}{2}}=0.\\
\end{split}
\end{equation}
Subtracting  \eqref{1bch2kj} from  the fractional diffusion equation \eqref{1bch2znm}, we have the following error equation
\begin{equation}\label{1bch2kj2c}
\begin{split}
\big( &(u-u_{h})_{t},v\big)_{D^{k}}-\varepsilon\big(p-p_{h},v\big)_{D^{k}}+
\big(p-p_{h},\psi\big)_{D^{k}}-\big(\Delta_{(\alpha-2)/2}(q-q_{h}),\psi\big)_{D^{k}}\\
&+\big(q-q_{h},\phi\big)_{D^{k}}+\big(\partial_{x}(u-u_{h}),\partial_{x}\phi\big)_{D^{k}}
-((\partial_{x}(u-u_{h}))^{*}[\phi]+(\partial_{x}\phi)^{*}[u-u_{h}]\big)_{k+\frac{1}{2}}=0.\\
\end{split}
\end{equation}
Denoting
\begin{equation}\label{91h}
\begin{split}
&\pi=\mathcal{P}^{+}u-u_{h},\quad \pi^{e}=\mathcal{P}^{+}u-u,\quad \epsilon=\mathcal{P}^{+}p-p_{h},\quad \epsilon^{e}=\mathcal{P}^{+}p-p,\\
 &\varphi^{e}=\mathcal{P}^{+}q-q,\quad \varphi=\mathcal{P}^{+}q-q_{h}.
\end{split}
\end{equation}
From the Galerkin orthogonality \eqref{1bch2kj2c}, we get
\begin{equation}\label{1bch2kj2h}
\begin{split}
\big( &(\pi-\pi^{e})_{t},v\big)_{D^{k}}-\varepsilon\big(\epsilon-\epsilon^{e},v\big)_{D^{k}}+
\big(\epsilon-\epsilon^{e},\psi\big)_{D^{k}}-\big(\Delta_{(\alpha-2)/2}(\varphi-\varphi^{e}),\psi\big)_{D^{k}}\\
&+\big(\varphi-\varphi^{e},\phi\big)_{D^{k}}+\big(\partial_{x}(\pi-\pi^{e}),\partial_{x}\phi\big)_{D^{k}}
+\bigg((\partial_{x}(\pi-\pi^{e}))^{*}[\phi]+\{\partial_{x}\phi\}[\pi-\pi^{e}]\bigg)_{k+\frac{1}{2}}=0.\\
\end{split}
\end{equation}
We take the test functions
\begin{equation}\label{91h}
\begin{split}
&v=\pi,\quad\psi=\epsilon-\varphi,\quad \phi=\pi,\\
\end{split}
\end{equation}
 we obtain
\begin{equation}\label{1bch2kj2kp}
\begin{split}
\big( &(\pi-\pi^{e})_{t},\pi\big)_{D^{k}}-\varepsilon\big(\epsilon-\epsilon^{e},\pi\big)_{D^{k}}+
\big(\epsilon-\epsilon^{e},\epsilon-\varphi\big)_{D^{k}}
-\big(\Delta_{(\alpha-2)/2}(\varphi-\varphi^{e}),\epsilon-\varphi\big)_{D^{k}}\\
&+\big(\varphi-\varphi^{e},\pi\big)_{D^{k}}+\big(\partial_{x}(\pi-\pi^{e}),\partial_{x}\pi\big)_{D^{k}}
+\bigg((\partial_{x}(\pi-\pi^{e}))^{*}[\pi]+\{\partial_{x}\pi\}[\pi-\pi^{e}]\bigg)_{k+\frac{1}{2}}=0.\\
\end{split}
\end{equation}
Summing over $k$ and from  the admissible  condition  \eqref{1bch2hj}  of the numerical  flux  defined  in  \eqref{flf}, we obtain that
\begin{equation}\label{1bch2kj2kp}
\begin{split}
&\big( \pi_{t},\pi\big)+
\big(\epsilon,\epsilon\big)+\big(\Delta_{(\alpha-2)/2}\varphi,\varphi\big)+
(1-\gamma)\|\partial_{x}\pi\|^{2}_{L^{2}(\Omega)}+\mu\sum_{k=1}^{K}\frac{[\pi]^{2}_{k+\frac{1}{2}}}{h}\\
&\leq \big( \pi^{e}_{t},\pi\big)-\big(\Delta_{(\alpha-2)/2}\varphi^{e},\epsilon-\varphi\big) +\big(\Delta_{(\alpha-2)/2}\varphi,\epsilon)+\big(\epsilon^{e},\epsilon-\varphi\big)\\
&\quad+\big(\epsilon,\varphi\big)+\varepsilon\big(\epsilon-\epsilon^{e},\pi\big)
-\big(\varphi-\varphi^{e},\pi\big)
+\big(\partial_{x}\pi^{e},\partial_{x}\pi\big)
+\sum_{k=1}^{K}\bigg((\partial_{x}(\pi^{e}))^{*}[\pi]+\{\partial_{x}\pi\}[\pi^{e}]\bigg)_{k+\frac{1}{2}}.\\
\end{split}
\end{equation}
Using the definitions of the projections $\mathcal{P}^{+}$ in \eqref{prh}, we get
\begin{equation}\label{1bch2kj2kp}
\begin{split}
&\big( \pi_{t},\pi\big)+
\big(\epsilon,\epsilon\big)+\big(\Delta_{(\alpha-2)/2}\varphi,\varphi\big)+
(1-\gamma)\|\partial_{x}\pi\|^{2}_{L^{2}(\Omega)}+\mu\sum_{k=1}^{K}\frac{[\pi]^{2}_{k+\frac{1}{2}}}{h}\\
&\leq \big( \pi^{e}_{t},\pi\big)-\big(\Delta_{(\alpha-2)/2}\varphi^{e},\epsilon-\varphi\big) +\big(\Delta_{(\alpha-2)/2}\varphi,\epsilon)+\big(\epsilon^{e},\epsilon-\varphi\big)\\
&\quad+\big(\epsilon,\varphi\big)+\varepsilon\big(\epsilon-\epsilon^{e},\pi\big)
-\big(\varphi-\varphi^{e},\pi\big).\\
\end{split}
\end{equation}
From the approximation results \eqref{sth} and Young's inequality, we obtain
\begin{equation}\label{1bch2kj2kp}
\begin{split}
&\big( \pi_{t},\pi\big)+
\big(\epsilon,\epsilon\big)+\big(\Delta_{(\alpha-2)/2}\varphi,\varphi\big)+
(1-\gamma)\|\partial_{x}\pi\|^{2}_{L^{2}(\Omega)}+\mu\sum_{k=1}^{K}\frac{[\pi]^{2}_{k+\frac{1}{2}}}{h}\\
&\leq c_{3}\|\pi\|^{2}_{L^{2}(\Omega)}+c_{2}\|\varphi\|^{2}_{L^{2}(\Omega)}
+c_{1}\|\epsilon\|^{2}_{L^{2}(\Omega)}
+Ch^{2N+2}.\\
\end{split}
\end{equation}
Recalling Lemma \ref{lg}, we get
\begin{equation}\label{1bch2kj2kp}
\begin{split}
&\big( \pi_{t},\pi\big)+
\big(\epsilon,\epsilon\big)+
(1-\gamma)\|\partial_{x}\pi\|^{2}_{L^{2}(\Omega)}+\mu\sum_{k=1}^{K}\frac{[\pi]^{2}_{k+\frac{1}{2}}}{h}\leq c_{3}\|\pi\|^{2}_{L^{2}(\Omega)}+c_{1}\|\epsilon\|^{2}_{L^{2}(\Omega)}
+Ch^{2N+2},\\
\end{split}
\end{equation}
provided $c_{1}$ is sufficiently small such that $c_{1}\leq1$, we obtain
\begin{equation}\label{1bch2kj2kp}
\begin{split}
&\big( \pi_{t},\pi\big)+
(1-\gamma)\|\partial_{x}\pi\|^{2}_{L^{2}(\Omega)}+\mu\sum_{k=1}^{K}\frac{[\pi]^{2}_{k+\frac{1}{2}}}{h}\leq c_{3}\|\pi\|^{2}_{L^{2}(\Omega)}
+Ch^{2N+2}.\\
\end{split}
\end{equation}
Employing Gronwall's lemma  and standard approximation theory, we can get \eqref{bnkl}. $\quad\Box$\\
For the more general fractional convection-diffusion problem, we introduce a few
results and then give the error estimate.
\begin{lem}\label{as1}  (see \cite{zhang2004error}). For any piecewise smooth function $\pi \in L^{2}(\Omega)$, on each
cell boundary point we define
\begin{equation}\label{91}
\begin{split}
\kappa(f^{*};\pi)\equiv \kappa(f^{*};\pi^{-},\pi^{+})=\left\{
  \begin{array}{ll}
    [w]^{-1}(f(\pi)-f^{*}(\pi)), & \hbox{if $[\pi]\neq0$;} \\
    \frac{1}{2}|f^{'}(\overline{\pi})|, & \hbox{if $[\pi]=0$,}
  \end{array}
\right.
\end{split}
\end{equation}
\end{lem}
where $f^{*}(\pi) \equiv f^{*}(\pi^{−},\pi^{+})$ is a monotone numerical flux consistent with the given flux
$f$. Then $\kappa(f^{*},\pi)$ is nonnegative and bounded for any ($\pi^{−},\pi^{+}) \in \mathbb{R}$. \\
To estimate the nonlinear part,  we can write it into the following
form
\begin{equation}\label{91bm44j}
\begin{split}
\sum_{k=1}^{K}\mathcal{H}_{k}(f;u,u_{h};\pi)=\sum_{k=1}^{K}\big(f(u)-f(u_{h}),\frac{\partial}{\partial x}\pi\big)_{D^{k}}+\sum_{k=1}^{K}((f(u)-f(u_{h})[\pi])_{k+\frac{1}{2}}+\sum_{k=1}^{K}((f(u_{h})
-\hat{f})[\pi])_{k+\frac{1}{2}}.\\
\end{split}
\end{equation}
We can rewrite \eqref{91bm44j} as:

\begin{equation}\label{91bm44}
\begin{split}
\sum_{k=1}^{K}\mathcal{H}_{k}(f;u,u_{h};\pi)=\sum_{k=1}^{K}\big(f(u)-f(u_{h}),\frac{\partial}{\partial x}\pi\big)_{D^{k}}+\sum_{k=1}^{K}((f(u)-f(\{u_{h}\})[\pi])_{k+\frac{1}{2}}+\sum_{k=1}^{K}((f(\{u_{h}\})
-\hat{f})[\pi])_{k+\frac{1}{2}}.\\
\end{split}
\end{equation}
\begin{lem}\label{lh} For $\mathcal{H}_{k}(f;u,u_{h};\pi)$ defined above, we have the following
estimate:
\begin{equation}\label{91bm4452}
\begin{split}
\sum_{k=1}^{K}\mathcal{H}_{k}(f;u,u_{h};\pi)\leq \gamma_{1}\|\pi_{x}\|^{2}_{L^{2}(\Omega)}
+\mu\sum_{k=1}^{K}\frac{[\pi]^{2}_{k+\frac{1}{2}}}{h}+
C\|\pi\|^{2}_{L^{2}(\Omega)}(1+h^{-1}\|\pi\|^{2}_{L^{2}(\Omega)})+h^{2N+2}.\\
\end{split}
\end{equation}
\end{lem}
\textbf{Proof.}
First, we give the estimate of the last term in \eqref{91bm44},  since the exact solution $u$ is continuous on each boundary point, we have that
\begin{equation}\label{91bm2}
\begin{split}
u_{h}=[\pi^{e}]-[\pi].\\
\end{split}
\end{equation}
Employing  Young's and the interpolation property \eqref{sth11}, we can easily show that
\begin{equation}\label{91bm34}
\begin{split}
\sum_{k=1}^{K}((f({u_{h}})
-\hat{f})[\pi])_{k+\frac{1}{2}}=&\sum_{k=1}^{K}(\kappa(\hat{f};u_{h})[u_{h}][\pi])_{k+\frac{1}{2}}=
\sum_{k=1}^{K}(\kappa(\hat{f};u_{h})[\pi^{e}][\pi])_{k+\frac{1}{2}}
-\sum_{k=1}^{K}(\kappa(\hat{f};u_{h})[\pi]^{2})_{k+\frac{1}{2}}\\
&\leq \frac{\mu}{3h}\sum_{k=1}^{K}[\pi]^{2}_{k+\frac{1}{2}}+\frac{3h}{4\mu}\sum_{k=1}^{K}[\kappa\pi^{e}]^{2}_{k+\frac{1}{2}}\\
&\leq \frac{\mu}{3h}\sum_{k=1}^{K}[\pi]^{2}_{k+\frac{1}{2}}+Ch^{2N+2}.
\end{split}
\end{equation}
For the first two terms of the right-hand \eqref{91bm44}, we would like to use the following Taylor expansions
\begin{equation}\label{91bm4}
\begin{split}
f(u)-f(u_{h})=f^{'}(u)(\pi-\pi^{e})-\frac{1}{2}f^{''}_{u}(\pi-\pi^{e})^{2},\\
\end{split}
\end{equation}

\begin{equation}\label{91bm4}
\begin{split}
f(u)-f(\{u_{h}\})=f^{'}(u)(\{\pi\}-\{\pi\}^{e})-\frac{1}{2}\tilde{f}^{''}_{u}(\{\pi\}-\{\pi\}^{e})^{2},\\
\end{split}
\end{equation}
where $f^{''}_{u}$ and $\tilde{f}^{''}_{u}$ are the mean values. These imply the following representation
\begin{equation}\label{91bm}
\begin{split}
\sum_{k=1}^{K}\big(f(u)-f(u_{h}),\frac{\partial}{\partial x}\pi\big)_{D^{k}}+\sum_{k=1}^{K}((f(u)-f(\{u_{h}\})[\pi])_{k+\frac{1}{2}}=I+II+III,\\
\end{split}
\end{equation}
where
\begin{equation}\label{91bm}
\begin{split}
I=\sum_{k=1}^{K}\big(f^{'}(u)\pi,\frac{\partial}{\partial x}\pi\big)_{D^{k}}+\sum_{k=1}^{K}f^{'}(u)\{\pi\}[\pi])_{k+\frac{1}{2}},\\
\end{split}
\end{equation}
\begin{equation}\label{91bm}
\begin{split}
II=-\Biggl(\sum_{k=1}^{K}\big(f^{'}(u)\pi^{e},\frac{\partial}{\partial x}\pi\big)_{D^{k}}+\sum_{k=1}^{K}f^{'}(u)\{\pi\}^{e}[\pi])_{k+\frac{1}{2}}\Biggl),\\
\end{split}
\end{equation}
\begin{equation}\label{91bm}
\begin{split}
III=-\frac{1}{2}\Biggl(\sum_{k=1}^{K}\big(f^{''}(\pi-\pi^{e})^{2},\frac{\partial}{\partial x}\pi\big)_{D^{k}}+\sum_{k=1}^{K}\tilde{f}^{''}(\{\pi\}-\{\pi\}^{e})^{2}[\pi])_{k+\frac{1}{2}}\Biggl),\\
\end{split}
\end{equation}
will be estimated separately as below.\\
For the $I$ term, a simple integration by parts gives
\begin{equation}\label{91bm5}
\begin{split}
I=-\frac{1}{2}\sum_{k=1}^{K}\big( f^{''}(u)\frac{\partial u}{\partial x},\pi\big)_{D^{k}}\leq C\|\pi\|^{2}_{L^{2}(\Omega)}.\\
\end{split}
\end{equation}
For the $II$ term, using Young's inequality, we obtain
\begin{equation}\label{91bm6}
\begin{split}
II\leq&\frac{\gamma_{1}}{2}\|\pi_{x}\|^{2}_{L^{2}(\Omega)}
+\frac{1}{2\gamma_{1}}\|f^{'}(u)\pi^{e}\|^{2}_{L^{2}(\Omega)}
+\frac{\mu}{3h}\sum_{k=1}^{K}[\pi]^{2}_{k+\frac{1}{2}}
+\frac{3h}{4\mu}\sum_{k=1}^{K}(f^{'}(u)\{\pi\}^{e})^{2}_{k+\frac{1}{2}}\\
&\leq \frac{\gamma_{1}}{2}\|\pi_{x}\|^{2}_{L^{2}(\Omega)}
+\frac{\mu}{3h}\sum_{k=1}^{K}[\pi]^{2}_{k+\frac{1}{2}}+Ch^{2N+2}.
\end{split}
\end{equation}
For the $III$ term, we use both projection and inverse inequalities, \eqref{sth11} and
\eqref{sth22}, to get
\begin{equation}\label{91bm}
\begin{split}
III\leq&\frac{\gamma_{1}}{2}\|\pi_{x}\|^{2}_{L^{2}(\Omega)}
+\frac{1}{2\gamma_{1}}\|f^{''}(\pi-\pi^{e})^{2}\|^{2}_{L^{2}(\Omega)}
+\frac{\mu}{3h}\sum_{k=1}^{K}[\pi]^{2}_{k+\frac{1}{2}}
+\frac{3h}{4\mu}\sum_{k=1}^{K}(\tilde{f}^{''}\{\pi-\pi^{e}\}^{2})^{2}_{k+\frac{1}{2}}\\
&\leq \frac{\gamma_{1}}{2}\|\pi_{x}\|^{2}_{L^{2}(\Omega)}
+\frac{\mu}{3h}\sum_{k=1}^{K}[\pi]^{2}_{k+\frac{1}{2}}+C\|\pi-\pi^{e}\|_{\infty}^{2}(\|\pi\|^{2}_{L^{2}(\Omega)}
+\|\pi^{e}\|^{2}_{L^{2}(\Omega)}+h\|\pi\|^{2}_{\Gamma_{h}}+h\|\pi^{e}\|^{2}_{\Gamma_{h}})\\
&\leq \frac{\gamma_{1}}{2}\|\pi_{x}\|^{2}_{L^{2}(\Omega)}
+\frac{\mu}{3h}\sum_{k=1}^{K}[\pi]^{2}_{k+\frac{1}{2}}+C(\|\pi\|^{2}_{L^{2}(\Omega)}
+h^{2N+2})(1+\|\pi-\pi^{e}\|_{\infty}^{2}).\\
\end{split}
\end{equation}
 Using the approximation results in \eqref{sth11} and
\eqref{sth22}, we have
\begin{equation}\label{91bm7}
\begin{split}
III\leq& \frac{\gamma_{1}}{2}\|\pi_{x}\|^{2}_{L^{2}(\Omega)}
+\frac{\mu}{3h}\sum_{k=1}^{K}[\pi]^{2}_{k+\frac{1}{2}}+C(\|\pi\|^{2}_{L^{2}(\Omega)}
+h^{2N+2})(1+h^{-1}\|\pi\|^{2}_{L^{2}(\Omega)}+h^{2N+1})\\
&\leq \frac{\gamma_{1}}{2}\|\pi_{x}\|^{2}_{L^{2}(\Omega)}
+\frac{\mu}{3h}\sum_{k=1}^{K}[\pi]^{2}_{k+\frac{1}{2}}+
C\|\pi\|^{2}_{L^{2}(\Omega_{h})}(1+h^{-1}\|\pi\|^{2}_{L^{2}(\Omega)})+h^{2N+2}.
\end{split}
\end{equation}
Combining \eqref{91bm34}, \eqref{91bm5},  \eqref{91bm6}, \eqref{91bm7}, and \eqref{91bm44}, we obtain
 \begin{equation}\label{91bm445}
\begin{split}
\sum_{k=1}^{K}\mathcal{H}_{k}(f;u,u_{h};\pi)\leq \gamma_{1}\|\pi_{x}\|^{2}_{L^{2}(\Omega)}
+\frac{\mu}{h}\sum_{k=1}^{K}[\pi]^{2}_{k+\frac{1}{2}}+
C\|\pi\|^{2}_{L^{2}(\Omega)}(1+h^{-1}\|\pi\|^{2}_{L^{2}(\Omega)})+h^{2N+2}.\quad\Box\\
\end{split}
\end{equation}

\begin{thm}\label{th2}
Let $u$ be the exact solution of the problem \eqref{25n}, and let $u_{h}$ be the numerical solution of the semi-discrete DDG scheme \eqref{1bch2znm}.  Then for small enough $h$, we have the following error estimates:
\begin{equation}\label{bnklz}
\begin{split}
&\|u-u_{h}\|_{L^{2}(\Omega)}\leq Ch^{N+1},\\
\end{split}
\end{equation}
\end{thm}
where the constant $C$ is dependent upon $T$ and some norms of the solutions. \\
\textbf{Proof}.
 The exact solution of the above  \eqref{25nl} satisfies
\begin{equation}\label{1bch2kjz}
\begin{split}
&\big( u_{t},v\big)_{D^{k}}-\varepsilon\big(p,v\big)_{D^{k}}+
\big(p,\psi\big)_{D^{k}}-\big(\Delta_{(\alpha-2)/2}q,\psi\big)_{D^{k}}+\big(q,\phi\big)_{D^{k}}-\big( f(u),\partial_{x} v\big)_{D^{k}}+\big(\hat{n}.f(u)^{*},v\big)_{\partial D^{k}}\\
&\quad+\big(\partial_{x}u,\partial_{x}\phi\big)_{D^{k}}+\bigg((\partial_{x}u)^{*}[\phi]+\{\partial_{x}\phi\}[u]\bigg)_{k+\frac{1}{2}}=0.\\
\end{split}
\end{equation}
Subtracting  \eqref{1bch2kjz} from   \eqref{1bch2znm}, we have the following error equation
\begin{equation}\label{1bch2kj2}
\begin{split}
\big( &(u-u_{h})_{t},v\big)_{D^{k}}-\varepsilon\big(p-p_{h},v\big)_{D^{k}}+
\big(p-p_{h},\psi\big)_{D^{k}}-\big(\Delta_{(\alpha-2)/2}(q-q_{h}),\psi\big)_{D^{k}}-\big( f(u)-f(u_{h}),\partial_{x} v\big)_{D^{k}}\\
&+\big(\hat{n}.(f(u)-f(u_{h}))^{*},v\big)_{\partial D^{k}}+\big(q-q_{h},\phi\big)_{D^{k}}+\big(\partial_{x}(u-u_{h}),\partial_{x}\phi\big)_{D^{k}}
+\bigg((\partial_{x}(u-u_{h}))^{*}[\phi]+\{\partial_{x}\phi\}[u-u_{h}]\bigg)_{k+\frac{1}{2}}=0.\\
\end{split}
\end{equation}
From \eqref{91bm44}, we can rewrite \eqref{1bch2kj2} as:
\begin{equation}\label{1bch2kj2f}
\begin{split}
\big( &(u-u_{h})_{t},v\big)_{D^{k}}-\varepsilon\big(p-p_{h},v\big)_{D^{k}}+
\big(p-p_{h},\psi\big)_{D^{k}}-\big(\Delta_{(\alpha-2)/2}(q-q_{h}),\psi\big)_{D^{k}}\\
&-\mathcal{H}_{k}(f;u,u_{h};v)+\big(q-q_{h},\phi\big)_{D^{k}}+\big(\partial_{x}(u-u_{h}),\partial_{x}\phi\big)_{D^{k}}
+\bigg((\partial_{x}(u-u_{h}))^{*}[\phi]+\{\partial_{x}\phi\}[u-u_{h}]\bigg)_{k+\frac{1}{2}}=0.\\
\end{split}
\end{equation}
From the Galerkin orthogonality \eqref{1bch2kj2f}, we get
\begin{equation}\label{1bch2kj2h}
\begin{split}
\big( &(\pi-\pi^{e})_{t},v\big)_{D^{k}}-\varepsilon\big(\epsilon-\epsilon^{e},v\big)_{D^{k}}+
\big(\epsilon-\epsilon^{e},\psi\big)_{D^{k}}-\big(\Delta_{(\alpha-2)/2}(\varphi-\varphi^{e}),\psi\big)_{D^{k}}-
\mathcal{H}_{k}(f;u,u_{h};v)\\
&+\big(\varphi-\varphi^{e},\phi\big)_{D^{k}}+\big(\partial_{x}(\pi-\pi^{e}),\partial_{x}\phi\big)_{D^{k}}
+\bigg((\partial_{x}(\pi-\pi^{e}))^{*}[\phi]+\{\partial_{x}\phi\}[\pi-\pi^{e}]\bigg)_{k+\frac{1}{2}}=0.\\
\end{split}
\end{equation}
Following the proof of Theorem \ref{th2}, we set
\begin{equation}\label{91h}
\begin{split}
&v=\pi,\quad\psi=\epsilon-\varphi,\quad \phi=\pi.\\
\end{split}
\end{equation}
Summing over $k$ and from  the admissible  condition  \eqref{1bch2hj}, definitions of the projections
$\mathcal{P}^{+}$, the approximation results \eqref{sth} and Young's inequality, we derive the following inequality:

\begin{equation}\label{1bch2kj2kp}
\begin{split}
&\big( \pi_{t},\pi\big)+
\big(\epsilon,\epsilon\big)+\big(\Delta_{(\alpha-2)/2}\varphi,\varphi\big)
+ (1-\gamma)\|\partial_{x}\pi\|^{2}_{L^{2}(\Omega)}+\mu\sum_{k=1}^{K}\frac{[\pi]^{2}_{k+\frac{1}{2}}}{h}\\
&\leq c_{3}\|\pi\|^{2}_{L^{2}(\Omega)}+c_{2}\|\varphi\|^{2}_{L^{2}(\Omega)}+c_{1}\|\epsilon\|^{2}_{L^{2}(\Omega)}
+Ch^{2N+2}+\sum_{k=1}^{K}\mathcal{H}_{k}(f;u,u_{h};\pi).\\
\end{split}
\end{equation}
Recalling Lemmas \ref{lg} and \ref{lh}, we get
\begin{equation}\label{1bch2kj2kp}
\begin{split}
&\big( \pi_{t},\pi\big)+
\big(\epsilon,\epsilon\big)
\leq c_{3}\|\pi\|^{2}_{L^{2}(\Omega)}+c_{1}\|\epsilon\|^{2}_{L^{2}(\Omega)}
+Ch^{2N+2}+C\|\pi\|^{2}_{L^{2}(\Omega)}(1+h^{-1}\|\pi\|^{2}_{L^{2}(\Omega)}+h^{2N+2}),\\
\end{split}
\end{equation}
provided $c_{1}$ is sufficiently small such that $c_{1}\leq1$, we obtain
\begin{equation}\label{1bch2kj2kp}
\begin{split}
&\big( \pi_{t},\pi\big)\leq c_{3}\|\pi\|^{2}_{L^{2}(\Omega)}
+C\|\pi\|^{2}_{L^{2}(\Omega)}(1+h^{-1}\|\pi\|^{2}_{L^{2}(\Omega)})+h^{2N+2}.\\
\end{split}
\end{equation}
From the Gronwall's lemma and standard approximation theory, the desired result follows. $\quad\Box$

\section{ DDG method for  nonlinear fractional
Schr\"{o}dinger equation}\label{sc4}
In this section, we present and analyze a direct discontinuous Galerkin method for the following:
\begin{equation}\label{sch1vn}
\begin{split}
&i\frac{\partial u}{\partial t}- \varepsilon_{1}(-\Delta)^{\frac{\alpha}{2}}u+ \varepsilon_{2}f(|u|^{2})u=0,\quad x\in \mathbb{R},\,\,t\in(0,T],\\
&u(x,0) = u_{0}(x),\quad x\in \mathbb{R},
\end{split}
\end{equation}
We introduce  the auxiliary variables $p,q$ and set
\begin{equation}\label{1a}
\begin{split}
&p=\Delta_{(\alpha-2)/2}q, \quad q=\frac{\partial^{2}}{\partial x^{2}}u,
\end{split}
\end{equation}
then, the nonlinear fractional
Schr\"{o}dinger  problem can be rewritten as
\begin{equation}\label{1bchcfaz}
\begin{split}
&i\frac{\partial u}{\partial t}+\varepsilon_{1}e+ \varepsilon_{2}f(|u|^{2})u=0,\\
&r=\Delta_{(\alpha-2)/2}q,\quad s=\frac{\partial^{2}}{\partial x^{2}}u.\\
\end{split}
\end{equation}
For actual numerical implementation, it might be more efficient if we decompose the complex function
$u(x,t)$ into its real and imaginary parts by writing
\begin{equation}\label{1b}
\begin{split}
u(x,t)=p(x,t)+iq(x,t),
\end{split}
\end{equation}
where $p$, $q$ are real functions. Under the new notation, the problem \eqref{1bchcfaz} can be written as
\begin{equation}\label{1bxc}
\begin{split}
&\frac{\partial p}{\partial t}+\varepsilon_{1}e+ \varepsilon_{2}f(p^{2}+q^{2})q=0,\\
&e=\Delta_{(\alpha-2)/2}r,\quad r=\frac{\partial^{2}}{\partial x^{2}}q,\\
&\frac{\partial q}{\partial t}-\varepsilon_{1}l- \varepsilon_{2}f(p^{2}+q^{2})p=0,\\
&l=\Delta_{(\alpha-2)/2}w,\quad w=\frac{\partial^{2}}{\partial x^{2}}p.\\
\end{split}
\end{equation}
We now apply the direct discontinuous Galerkin discretization to \eqref{1bxc}: find $p_{h}, q_{h}, e_{h},l_{h},r_{h},w_{h}\in V_{k}^{N}$, such that for all test functions $\vartheta,\varphi,\phi,\chi,\zeta,\psi\in V_{k}^{N}$,
\begin{equation}\label{1bch2}
\begin{split}
&\big( (p_{h})_{t},\vartheta\big)_{D^{k}}+\varepsilon_{1}\big(e_{h},\vartheta\big)_{D^{k}}+ \varepsilon_{2}\big(f(p_{h}^{2}+q_{h}^{2})q_{h},\vartheta\big)_{D^{k}}=0,\\
&\big(e_{h},\varphi\big)_{D^{k}}=\big(\Delta_{(\alpha-2)/2}r_{h},\varphi\big)_{D^{k}},\\
&\big(r_{h},\phi\big)_{D^{k}}=-\big(\partial_{x}q_{h},\partial_{x}\phi\big)_{D^{k}}-\bigg((\partial_{x}q_{h})^{*}[\phi]
+\{\partial_{x}\phi\}[q_{h}]\bigg)_{k+\frac{1}{2}},\\
&\big((q_{h})_{t},\chi\big)_{D^{k}}-\varepsilon_{1}\big(l_{h},\chi\big)_{D^{k}}- \varepsilon_{2}\big(f(p_{h}^{2}+q_{h}^{2})p_{h},\chi\big)_{D^{k}}=0,\\
&\big(l_{h},\zeta\big)_{D^{k}}=\big(\Delta_{(\alpha-2)/2}w_{h},\zeta\big)_{D^{k}},\\
&\big(w_{h},\psi\big)_{D^{k}}=-\big(\partial_{x}p_{h},\partial_{x}\psi\big)_{D^{k}}-\bigg((\partial_{x}p_{h})^{*}[\psi]
+\{\partial_{x}\psi\}[p_{h}]\bigg)_{k+\frac{1}{2}},\\
\end{split}
\end{equation}

\section{ Stability and error estimates}\label{sc5}
 In the following we discuss stability and accuracy of the proposed scheme, for the nonlinear fractional Schr\"{o}dinger   problem.
\subsection{ Stability analysis }  In order to carry out the analysis of the DDG scheme, we have the following results.
\begin{thm}\label{tt4h} ($L^{2}$ stability).
 The solution to the scheme \eqref{1bch2} satisfies the $
\|u_{h}(x,T)\|_{L^{2}(\Omega)}\leq c\|u_{0}(x)\|_{L^{2}(\Omega)}$ for any $T>0$.
\end{thm}
\textbf{Proof.} Set $(\vartheta,\varphi,\phi,\chi,\zeta,\psi)=(p_{h},e_{h}-r_{h},q_{h},q_{h},l_{h}-w_{h},p_{h})$ in \eqref{1bch2}, we get
\begin{equation}\label{1bch1f}
\begin{split}
&\big((p_{h})_{t},p_{h}\big)_{D^{k}}+\big((q_{h})_{t},q_{h}\big)_{D^{k}}
+\big(e_{h},e_{h}\big)_{D^{k}}+\big(l_{h},l_{h}\big)_{D^{k}}
+\big(\Delta_{(\alpha-2)/2}w_{h},w_{h}\big)_{D^{k}}+\big(\Delta_{(\alpha-2)/2}r_{h},r_{h}\big)_{D^{k}}\\
&=\big(\Delta_{(\alpha-2)/2}w_{h},l_{h}\big)_{D^{k}}
+\big(\Delta_{(\alpha-2)/2}r_{h},e_{h}\big)_{D^{k}}-\big( r_{h},q_{h}\big)_{D^{k}}-\big( w_{h},p_{h}\big)_{D^{k}}+\big( e_{h},r_{h}\big)_{D^{k}}+\big( l_{h},w_{h}\big)_{D^{k}}\\
&\quad-\varepsilon_{1}\big( e_{h},p_{h}\big)_{D^{k}}+\varepsilon_{1}\big( l_{h},q_{h}\big)_{D^{k}}
-\big(\partial_{x}q_{h},\partial_{x}\phi\big)_{D^{k}}-\bigg((\partial_{x}q_{h})^{*}[\phi]
+\{\partial_{x}\phi\}[q_{h}]\bigg)_{k+\frac{1}{2}}-\big(\partial_{x}p_{h},\partial_{x}\psi\big)_{D^{k}}\\
&\quad-\bigg((\partial_{x}p_{h})^{*}[\psi]
+\{\partial_{x}\psi\}[p_{h}]\bigg)_{k+\frac{1}{2}}.\\
\end{split}
\end{equation}
Summing over $k$ and from  the admissible  condition  \ref{tt4hzx}  of the numerical  flux  defined  in  \eqref{flf}, we obtain that
\begin{equation}\label{1bch1}
\begin{split}
&\big((p_{h})_{t},p_{h}\big)+\big((q_{h})_{t},q_{h}\big)
+\big(e_{h},e_{h}\big)+\big(l_{h},l_{h}\big)
+\big(\Delta_{(\alpha-2)/2}w_{h},w_{h}\big)+\big(\Delta_{(\alpha-2)/2}r_{h},r_{h}\big)\\
&+(1-\gamma_{2})\|\partial_{x}q_{h}\|^{2}_{L^{2}(\Omega)}
+\mu_{2}\sum_{k=1}^{K}\frac{[q_{h}]^{2}_{k+\frac{1}{2}}}{h}+(1-\gamma_{1})\|\partial_{x}p_{h}\|^{2}_{L^{2}(\Omega)}
+\mu_{1}\sum_{k=1}^{K}\frac{[p_{h}]^{2}_{k+\frac{1}{2}}}{h}\\
&\leq\big(\Delta_{(\alpha-2)/2}w_{h},l_{h}\big)
+\big(\Delta_{(\alpha-2)/2}r_{h},e_{h}\big)-\big( r_{h},q_{h}\big)_{D^{k}}-\big( w_{h},p_{h}\big)+\big( e_{h},r_{h}\big)+\big( l_{h},w_{h}\big)\\
&\quad-\varepsilon_{1}\big( e_{h},p_{h}\big)+\varepsilon_{1}\big( l_{h},q_{h}\big)
,\\
\end{split}
\end{equation}
Employing Young's inequality and  Lemma \ref{lga2}, we obtain
\begin{equation}\label{1bch1}
\begin{split}
&\big((p_{h})_{t},p_{h}\big)+\big((q_{h})_{t},q_{h}\big)
+\big(e_{h},e_{h}\big)+\big(l_{h},l_{h}\big)
+\big(\Delta_{(\alpha-2)/2}w_{h},w_{h}\big)+\big(\Delta_{(\alpha-2)/2}r_{h},r_{h}\big)\\
&+(1-\gamma_{2})\|\partial_{x}q_{h}\|^{2}_{L^{2}(\Omega)}
+\mu_{2}\sum_{k=1}^{K}\frac{[q_{h}]^{2}_{k+\frac{1}{2}}}{h}+(1-\gamma_{1})\|\partial_{x}p_{h}\|^{2}_{L^{2}(\Omega)}
+\mu_{1}\sum_{k=1}^{K}\frac{[p_{h}]^{2}_{k+\frac{1}{2}}}{h}\\
&\leq c_{4}\|w_{h}\|^{2}_{L^{2}(\Omega)}+c_{3}\|r_{h}\|^{2}_{L^{2}(\Omega)}
+c_{2}\|e_{h}\|^{2}_{L^{2}(\Omega)}+c_{1}\|l_{h}\|^{2}_{L^{2}(\Omega)}
+c_{5}\|p_{h}\|^{2}_{L^{2}(\Omega)}+c_{6}\|q_{h}\|^{2}_{L^{2}(\Omega)}
.\\
\end{split}
\end{equation}
Recalling Lemma \ref{lg}, provided $c_{i},\,\,i=1,2$ are sufficiently small such that $c_{i}\leq1$, we obtain
\begin{equation}\label{1bch1}
\begin{split}
&\big((p_{h})_{t},p_{h}\big)+\big((q_{h})_{t},q_{h}\big)
+(1-\gamma_{2})\|\partial_{x}q_{h}\|^{2}_{L^{2}(\Omega)}
+\mu_{2}\sum_{k=1}^{K}\frac{[q_{h}]^{2}_{k+\frac{1}{2}}}{h}+(1-\gamma_{1})\|\partial_{x}p_{h}\|^{2}_{L^{2}(\Omega)}
+\mu_{1}\sum_{k=1}^{K}\frac{[p_{h}]^{2}_{k+\frac{1}{2}}}{h}\\
&\leq c_{5}\|p_{h}\|^{2}_{L^{2}(\Omega)}+c_{6}\|q_{h}\|^{2}_{L^{2}(\Omega)}
.\\
\end{split}
\end{equation}
Employing Gronwall's lemma, we obtain $\|u_{h}(x,T)\|_{L^{2}(\Omega)}\leq c\|u_{0}(x)\|_{L^{2}(\Omega)}$. $\quad\Box$
\subsection{Error estimates}
We consider the linear fractional Schr\"{o}dinger equation
\begin{equation}\label{1bchj}
\begin{split}
&i\frac{\partial u}{\partial t}-\varepsilon_{1}(-\Delta) ^{\frac{\alpha}{2}}u+\varepsilon_{2}u=0.\\
\end{split}
\end{equation}
It is easy to verify that the exact solution of the above  \eqref{1bchj} satisfies
\begin{equation}\label{1bch2j}
\begin{split}
&\big( p_{t},\vartheta\big)_{D^{k}}+\varepsilon_{1}\big(e,\vartheta\big)_{D^{k}}+ \varepsilon_{2}\big(q,\vartheta\big)_{D^{k}}=0,\\
&\big(e,\varphi\big)_{D^{k}}=\big(\Delta_{(\alpha-2)/2}r,\varphi\big)_{D^{k}},\\
&\big(r,\phi\big)_{D^{k}}=-\big(\partial_{x}q,\partial_{x}\phi\big)_{D^{k}}-\bigg((\partial_{x}q)^{*}[\phi]
+\{\partial_{x}\phi\}[q]\bigg)_{k+\frac{1}{2}},\\
&\big(q_{t},\chi\big)_{D^{k}}-\varepsilon_{1}\big(l,\chi\big)_{D^{k}}- \varepsilon_{2}\big(p,\chi\big)_{D^{k}}=0,\\
&\big(l,\zeta\big)_{D^{k}}=\big(\Delta_{(\alpha-2)/2}w,\zeta\big)_{D^{k}},\\
&\big(w,\psi\big)_{D^{k}}=-\big(\partial_{x}p,\partial_{x}\psi\big)_{D^{k}}-\bigg((\partial_{x}p)^{*}[\psi]
+\{\partial_{x}\psi\}[p]\bigg)_{k+\frac{1}{2}}.\\
\end{split}
\end{equation}
Subtracting  \eqref{1bch2j}, from  the linear fractional Schr\"{o}dinger equation \eqref{1bch2}, we have the following error equation
\begin{equation}\label{1bch2k2}
\begin{split}
\big( &(p-p_{h})_{t},\vartheta\big)_{D^{k}}+\big((q-q_{h})_{t},\chi\big)_{D^{k}}
-\big(\Delta_{(\alpha-2)/2}(r-r_{h}),\varphi\big)_{D^{k}}
-\big(\Delta_{(\alpha-2)/2}(w-w_{h}),\zeta\big)_{D^{k}}+\big((q-q_{h})_{x},\phi_{x}\big)_{D^{k}}\\
&
+\big((p-p_{h})_{x},\psi_{x}\big)_{D^{k}}+ \varepsilon_{2}\big(q-q_{h},\vartheta\big)_{D^{k}}- \varepsilon_{2}\big(p-p_{h},\chi\big)_{D^{k}}
+\big(r-r_{h},\phi\big)_{D^{k}}
+\big(l-l_{h},\zeta\big)_{D^{k}}+\big(e-e_{h},\varphi\big)_{D^{k}}\\
&+\big(w-w_{h},\psi\big)_{D^{k}} -\varepsilon_{1}\big(l-l_{h},\chi\big)_{D^{k}}+\varepsilon_{1}\big(e-e_{h},\vartheta\big)_{D^{k}}
+\bigg((\partial_{x}(p-p_{h}))^{*}[\psi]
+\{\partial_{x}\psi\}[p-p_{h}]\bigg)_{k+\frac{1}{2}}\\
&+\bigg((\partial_{x}(q-q_{h}))^{*}[\phi]
+\{\partial_{x}\phi\}[q-q_{h}]\bigg)_{k+\frac{1}{2}}=0.\\
\end{split}
\end{equation}
\begin{thm}\label{fgrh}
Let $u$ be the exact solution of the problem \eqref{1bchj}, and let $u_{h}$ be the numerical solution of the semi-discrete DDG scheme \eqref{1bch2}.  Then for small enough $h$, we have the following error estimates:
\begin{equation}\label{tt7hb}
\begin{split}
&\|u(.,t)-u_{h}(.,t)\|_{L^{2}(\Omega_{h})}\leq Ch^{N+1},\\
\end{split}
\end{equation}
\end{thm}
where the constant $C$ is dependent upon $T$ and some norms of the solutions. \\
\textbf{Proof}. Let
\begin{equation}\label{91h}
\begin{split}
&\pi=\mathcal{P}^{+}p-p_{h},\quad \pi^{e}=\mathcal{P}^{+}p-p,\quad \epsilon=\mathcal{P}^{+}r-r_{h},\quad \epsilon^{e}=\mathcal{P}^{+}r-r,\quad\phi_{1}=\mathcal{P}^{+}e-e_{h},\quad\phi_{1}^{e}=\mathcal{P}^{+}e-e,\\
&  \sigma=\mathcal{P}^{+}q-q_{h}, \quad \sigma^{e}=\mathcal{P}^{+}q-q,\quad\phi_{2}=\mathcal{P}^{+}l-l_{h},\quad\phi_{2}^{e}=\mathcal{P}^{+}l-l,  \quad\varphi_{1}=\mathcal{P}^{+}w-w_{h},\quad \varphi^{e}_{1}=\mathcal{P}^{+}w-w.
\end{split}
\end{equation}
From the Galerkin orthogonality \eqref{1bch2k2}, we get
\begin{equation}\label{1bch2k2}
\begin{split}
\big( &(\pi-\pi^{e})_{t},\vartheta\big)_{D^{k}}+\big((\sigma-\sigma^{e})_{t},\chi\big)_{D^{k}}
-\big(\Delta_{(\alpha-2)/2}(\epsilon-\epsilon^{e}),\varphi\big)_{D^{k}}
-\big(\Delta_{(\alpha-2)/2}(\varphi_{1}-\varphi^{e}_{1}),\zeta\big)_{D^{k}}
+\big((\sigma-\sigma^{e})_{x},\phi_{x}\big)_{D^{k}}\\
&
+\big((\pi-\pi^{e})_{x},\psi_{x}\big)_{D^{k}}+ \varepsilon_{2}\big(\sigma-\sigma^{e},\vartheta\big)_{D^{k}}- \varepsilon_{2}\big(\pi-\pi^{e},\chi\big)_{D^{k}}
+\big(\epsilon-\epsilon^{e},\phi\big)_{D^{k}}
+\big(\phi_{2}-\phi_{2}^{e},\zeta\big)_{D^{k}}+\big(\phi_{1}-\phi_{1}^{e},\varphi\big)_{D^{k}}\\
&+\big(\varphi_{1}-\varphi^{e}_{1},\psi\big)_{D^{k}} -\varepsilon_{1}\big(\phi_{2}-\phi_{2}^{e},\chi\big)_{D^{k}}+\varepsilon_{1}\big(\phi_{1}-\phi_{1}^{e},\vartheta\big)_{D^{k}}
+\bigg((\partial_{x}(\pi-\pi^{e}))^{*}[\psi]
+\{\partial_{x}\psi\}[\pi-\pi^{e}]\bigg)_{k+\frac{1}{2}}\\
&\bigg((\partial_{x}(\sigma-\sigma^{e}))^{*}[\phi]
+\{\partial_{x}\phi\}[\sigma-\sigma^{e}]\bigg)_{k+\frac{1}{2}}=0.\\
\end{split}
\end{equation}
 We take the test functions
\begin{equation}\label{91h}
\begin{split}
\vartheta=\pi,\quad\varphi=\phi_{1}-\epsilon,\quad \phi=\sigma,\quad  \chi=\sigma,\quad\zeta=\phi_{2}-\varphi,\quad \psi=\pi,
\end{split}
\end{equation}
 we obtain
\begin{equation}\label{1bch2k2}
\begin{split}
\big( &(\pi-\pi^{e})_{t},\pi\big)_{D^{k}}+\big((\sigma-\sigma^{e})_{t},\sigma\big)_{D^{k}}
-\big(\Delta_{(\alpha-2)/2}(\epsilon-\epsilon^{e}),\phi_{1}-\epsilon\big)_{D^{k}}
-\big(\Delta_{(\alpha-2)/2}(\varphi_{1}-\varphi^{e}_{1}),\phi_{2}-\varphi_{1}\big)_{D^{k}}
\\
&+\big((\sigma-\sigma^{e})_{x},\sigma_{x}\big)_{D^{k}}
+\big((\pi-\pi^{e})_{x},\pi_{x}\big)_{D^{k}}+ \varepsilon_{2}\big(\sigma-\sigma^{e},\pi\big)_{D^{k}}- \varepsilon_{2}\big(\pi-\pi^{e},\sigma\big)_{D^{k}}
+\big(\epsilon-\epsilon^{e},\sigma\big)_{D^{k}}
\\
&+\big(\phi_{2}-\phi_{2}^{e},\phi_{2}-\varphi_{1}\big)_{D^{k}}+\big(\phi_{1}-\phi_{1}^{e},\phi_{1}-\epsilon\big)_{D^{k}}
+\big(\varphi_{1}-\varphi^{e}_{1},\pi\big)_{D^{k}} -\varepsilon_{1}\big(\phi_{2}-\phi_{2}^{e},\sigma\big)_{D^{k}}
+\varepsilon_{1}\big(\phi_{1}-\phi_{1}^{e},\pi\big)_{D^{k}}\\
&
+\bigg((\partial_{x}(\pi-\pi^{e}))^{*}[\pi]
+\{\partial_{x}\pi\}[\pi-\pi^{e}]\bigg)_{k+\frac{1}{2}}+\bigg((\partial_{x}(\sigma-\sigma^{e}))^{*}[\sigma]
+\{\partial_{x}\sigma\}[\sigma-\sigma^{e}]\bigg)_{k+\frac{1}{2}}=0.\\
\end{split}
\end{equation}
Summing over $k$ and from  the admissible  condition  \ref{tt4hzx}  of the numerical  flux  defined  in  \eqref{flf}, we obtain that
\begin{equation}\label{1bch2k2}
\begin{split}
\big( &\pi_{t},\pi\big)+\big(\sigma_{t},\sigma\big)
+\big(\Delta_{(\alpha-2)/2}\epsilon,\epsilon\big)
+\big(\Delta_{(\alpha-2)/2}\varphi_{1},\varphi_{1}\big)
+(1-\gamma_{2})\|\partial_{x}\sigma_{h}\|^{2}_{L^{2}(\Omega)}
\\
&+\mu_{2}\sum_{k=1}^{K}\frac{[\sigma_{h}]^{2}_{k+\frac{1}{2}}}{h}+(1-\gamma_{1})\|\partial_{x}\pi_{h}\|^{2}_{L^{2}(\Omega)}
+\mu_{1}\sum_{k=1}^{K}\frac{[\pi_{h}]^{2}_{k+\frac{1}{2}}}{h}
+\big(\phi_{2},\phi_{2}\big)+\big(\phi_{1},\phi_{1}\big)\\
&\leq \big( \pi^{e}_{t},\pi\big)+\big(\sigma^{e}_{t},\sigma\big)
-\big(\Delta_{(\alpha-2)/2}\epsilon^{e},\phi_{1}-\epsilon\big)
+\big(\Delta_{(\alpha-2)/2}\epsilon,\phi_{1}\big)
-\big(\Delta_{(\alpha-2)/2}\varphi^{e}_{1},\phi_{2}-\varphi_{1}\big)\\
&\quad
+\big(\Delta_{(\alpha-2)/2}\varphi_{1},\phi_{2}\big)
+\big(\sigma^{e}_{x},\sigma_{x}\big)+\big(\pi^{e}_{x},\pi_{x}\big)
- \varepsilon_{2}\big(\sigma^{e},\pi\big)- \varepsilon_{2}\big(\pi^{e},\sigma\big)-\big(\epsilon-\epsilon^{e},\sigma\big)\\
&\quad
+\big(\phi_{2},\varphi_{1}\big)+\big(\phi_{1},\epsilon\big)
+\big(\phi_{2}^{e},\phi_{2}-\varphi_{1}\big)+\big(\phi_{1}^{e},\phi_{1}-\epsilon\big)
+\varepsilon_{1}\big(\phi_{2}-\phi_{2}^{e},\sigma\big)
-\varepsilon_{1}\big(\phi_{1}-\phi_{1}^{e},\pi\big)\\
&\quad
-\big(\varphi_{1}-\varphi^{e}_{1},\pi\big)
+\sum_{k=1}^{K}\biggl(\bigg((\partial_{x}(\pi^{e}))^{*}[\pi]
+\{\partial_{x}\pi\}[\pi^{e}]\bigg)_{k+\frac{1}{2}}+\bigg((\partial_{x}(\sigma^{e}))^{*}[\sigma]
+\{\partial_{x}\sigma\}[\sigma^{e}]\bigg)_{k+\frac{1}{2}}\biggl)
\end{split}
\end{equation}
Using the definitions of the projections $\mathcal{P}^{+}$ in \eqref{prh}, we get
\begin{equation}\label{1bch2k2}
\begin{split}
\big( &\pi_{t},\pi\big)+\big(\sigma_{t},\sigma\big)
+\big(\Delta_{(\alpha-2)/2}\epsilon,\epsilon\big)
+\big(\Delta_{(\alpha-2)/2}\varphi_{1},\varphi_{1}\big)
+(1-\gamma_{2})\|\partial_{x}\sigma_{h}\|^{2}_{L^{2}(\Omega)}
\\
&+\mu_{2}\sum_{k=1}^{K}\frac{[\sigma_{h}]^{2}_{k+\frac{1}{2}}}{h}+(1-\gamma_{1})\|\partial_{x}\pi_{h}\|^{2}_{L^{2}(\Omega)}
+\mu_{1}\sum_{k=1}^{K}\frac{[\pi_{h}]^{2}_{k+\frac{1}{2}}}{h}
+\big(\phi_{2},\phi_{2}\big)+\big(\phi_{1},\phi_{1}\big)\\
&\leq \big( \pi^{e}_{t},\pi\big)+\big(\sigma^{e}_{t},\sigma\big)
-\big(\Delta_{(\alpha-2)/2}\epsilon^{e},\phi_{1}-\epsilon\big)+\big(\Delta_{(\alpha-2)/2}\epsilon,\phi_{1}\big)
-\big(\Delta_{(\alpha-2)/2}\varphi^{e}_{1},\phi_{2}-\varphi_{1}\big)\\
&\quad
+\big(\Delta_{(\alpha-2)/2}\varphi_{1},\phi_{2}\big)
- \varepsilon_{2}\big(\sigma^{e},\pi\big)- \varepsilon_{2}\big(\pi^{e},\sigma\big)-\big(\epsilon-\epsilon^{e},\sigma\big)
+\big(\phi_{2},\varphi_{1}\big)\\
&\quad
+\big(\phi_{1},\epsilon\big)
+\big(\phi_{2}^{e},\phi_{2}-\varphi_{1}\big)+\big(\phi_{1}^{e},\phi_{1}-\epsilon\big)
+\varepsilon_{1}\big(\phi_{2}-\phi_{2}^{e},\sigma\big)
-\varepsilon_{1}\big(\phi_{1}-\phi_{1}^{e},\pi\big)-\big(\varphi_{1}-\varphi^{e}_{1},\pi\big)\\
\end{split}
\end{equation}
From the approximation results \eqref{sth}, Young's inequality and recalling Lemma \ref{lg}, we obtain
\begin{equation}\label{1bch2k2}
\begin{split}
\big( &\pi_{t},\pi\big)_{D^{k}}+\big(\sigma_{t},\sigma\big)
+(1-\gamma_{2})\|\partial_{x}\sigma_{h}\|^{2}_{L^{2}(\Omega)}
+\mu_{2}\sum_{k=1}^{K}\frac{[\sigma_{h}]^{2}_{k+\frac{1}{2}}}{h}+(1-\gamma_{1})\|\partial_{x}\pi_{h}\|^{2}_{L^{2}(\Omega)}
+\mu_{1}\sum_{k=1}^{K}\frac{[\pi_{h}]^{2}_{k+\frac{1}{2}}}{h}\\
&
+\big(\phi_{2},\phi_{2}\big)+\big(\phi_{1},\phi_{1}\big)_{D^{k}}\leq c_{1}\|\phi_{1}\|^{2}_{L^{2}(\Omega)}+c_{1}\|\phi_{2}\|^{2}_{L^{2}(\Omega)}
+c_{3}\|\pi\|^{2}_{L^{2}(\Omega)}+c_{4}\|\sigma\|^{2}_{L^{2}(\Omega)}\\
\end{split}
\end{equation}
provided $c_{i},\,\,i=1,2$ are sufficiently small such that $c_{i}\leq1$, employing Gronwall's lemma  and standard approximation theory, we can get \eqref{tt7hb}. $\quad\Box$\\

\section{ DDG method for  strongly nonlinear  coupled
fractional Schr\"{o}dinger equations}\label{sc6}

In this section, we present DDG method for the strongly  coupled  nonlinear fractional Schr\"{o}dinger equations

\begin{equation}\label{1bchcc}
\begin{split}
&i\frac{\partial u_{1}}{\partial t}- \varepsilon_{1}(-\Delta)^{\frac{\alpha}{2}}u_{1}+ \varpi_{1}u_{1}+ \varpi_{2}u_{2}+ \varepsilon_{2}f(|u_{1}|^{2},|u_{2}|^{2})u_{1}=0,\\
&i\frac{\partial u_{2}}{\partial t}- \varepsilon_{3}(-\Delta)^{\frac{\alpha}{2}}u_{2}+\varpi_{2} u_{1}+ \varpi_{1}u_{2}+\varepsilon_{4} g(|u_{1}|^{2},|u_{2}|^{2})u_{2}=0.\\
\end{split}
\end{equation}
To define the DDG method, we rewrite  \eqref{1bchcc} as:
\begin{equation}\label{1bch}
\begin{split}
&i\frac{\partial u_{1}}{\partial t}+\varepsilon_{1}e+ \varpi_{1}u_{1}+\varpi_{2} u_{2}+ \varepsilon_{2}f(|u_{1}|^{2},|u_{2}|^{2})u_{1}=0,\\
&e=\Delta_{(\alpha-2)/2}r,\quad r=\frac{\partial^{2}}{\partial x^{2}}u_{1},\\
&i\frac{\partial u_{2}}{\partial t}+\varepsilon_{3}l+ \varpi_{2}u_{1}+\varpi_{1} u_{2}+ \varepsilon_{4}g(|u_{1}|^{2},|u_{2}|^{2})u_{2}=0,\\
&l=\Delta_{(\alpha-2)/2}w,\quad w=\frac{\partial^{2}}{\partial x^{2}}u_{2}.\\
\end{split}
\end{equation}
We decompose the complex functions $u(x, t)$ and $v(x, t)$ into their real and imaginary parts. Setting
$u_{1}(x, t)=p(x, t) + iq(x, t)$ and $u_{2}(x, t) = \upsilon(x, t) + i\theta(x, t)$ in system \eqref{1bchcc}, we can obtain the following coupled system
\begin{equation}\label{1bx}
\begin{split}
&\frac{\partial p}{\partial t}+\varepsilon_{1}Q+ \varpi_{1}q+ \varpi_{2}\theta+ \varepsilon_{2}f(|u_{1}|^{2},|u_{2}|^{2})q=0,\\
&Q=\Delta_{(\alpha-2)/2}r,\quad r=\frac{\partial^{2}}{\partial x^{2}}q,\\
&\frac{\partial q}{\partial t}-\varepsilon_{1}H- \varpi_{1}p- \varpi_{2}\upsilon- \varepsilon_{2}f(|u_{1}|^{2},|u_{2}|^{2})p=0,\\
&H=\Delta_{(\alpha-2)/2}w,\quad w=\frac{\partial^{2}}{\partial x^{2}}p,\\
&\frac{\partial \upsilon}{\partial t}+\varepsilon_{3}L+\varpi_{3}q+ \varpi_{4}\theta+ \varepsilon_{4}g(|u_{1}|^{2},|u_{2}|^{2})\theta=0,\\
&L=\Delta_{(\alpha-2)/2}\rho,\quad\rho=\frac{\partial^{2}}{\partial x^{2}}\theta,\\
&\frac{\partial \theta}{\partial t}-\varepsilon_{3}E-\varpi_{2} p- \varpi_{1}\upsilon- \varepsilon_{4}g(|u_{1}|^{2},|u_{2}|^{2})\upsilon=0,\\
&E=\Delta_{(\alpha-2)/2}\xi,\quad\xi=\frac{\partial^{2}}{\partial x^{2}}\upsilon.\\
\end{split}
\end{equation}
We now define  DDG scheme as follows: find $p_{h}, q_{h}, Q_{h},r_{h},H_{h},w_{h}$,
$\upsilon_{h},\theta_{h},L_{h},\rho_{h},E_{h},\xi_{h} \in V_{k}^{N}$, such that for all test functions $\vartheta_{1},\beta_{1},\phi,\chi,\beta_{2},\psi$,
$\gamma,\beta_{3},\delta,o,\beta_{4},\omega\in V_{k}^{N}$,

\begin{equation}\label{1bch2ccf}
\begin{split}
&\big(\frac{\partial p_{h}}{\partial t},\vartheta_{1}\big)_{D^{k}}+\varepsilon_{1}\big(Q_{h},\vartheta_{1}\big)_{D^{k}}+ \varpi_{1}\big(q_{h},\vartheta_{1}\big)_{D^{k}}+\varpi_{2}\big(\theta_{h},\vartheta_{1}\big)_{D^{k}}+ \varepsilon_{2}\big(f(|u_{1}|^{2},|u_{2}|^{2})q_{h},\vartheta_{1}\big)_{D^{k}}=0,\\
&\big(Q_{h},\beta_{1}\big)_{D^{k}}=\big(\Delta_{(\alpha-2)/2}r_{h},\beta_{1}\big)_{D^{k}},\\
&\big(r_{h},\phi\big)_{D^{k}}=-\big(\partial_{x}q_{h},\partial_{x}\phi\big)_{D^{k}}-\bigg((\partial_{x}q_{h})^{*}[\phi]
+\{\partial_{x}\phi\}[q_{h}]\bigg)_{k+\frac{1}{2}},\\
&\big(\frac{\partial q_{h}}{\partial t},\chi\big)_{D^{k}}-\varepsilon_{1}\big(H_{h},\chi\big)_{D^{k}}
-\varpi_{1} \big(p_{h},\chi\big)_{D^{k}}- \varpi_{2}\big(\upsilon_{h},\chi\big)_{D^{k}}- \varepsilon_{2}\big(f(|u_{1}|^{2},|u_{2}|^{2})p_{h},\chi\big)_{D^{k}}=0,\\
&\big(H_{h},\beta_{2}\big)_{D^{k}}=\big(\Delta_{(\alpha-2)/2}w_{h},\beta_{2}\big)_{D^{k}},\\
&\big(w_{h},\psi\big)_{D^{k}}=-\big(\partial_{x}p_{h},\partial_{x}\psi\big)_{D^{k}}-\bigg((\partial_{x}p_{h})^{*}[\psi]
+\{\partial_{x}\psi\}[p_{h}]\bigg)_{k+\frac{1}{2}},\\
&\big(\frac{\partial \upsilon_{h}}{\partial t},\gamma\big)_{D^{k}}+\varepsilon_{3}\big(L_{h},\gamma\big)_{D^{k}}
+\varpi_{2}\big(q_{h},\gamma\big)_{D^{k}}+\varpi_{1}\big(\theta_{h},\gamma\big)_{D^{k}}+ \varepsilon_{4}\big(g(|u_{1}|^{2},|u_{2}|^{2})\theta_{h},\gamma\big)_{D^{k}}=0,\\
&\big(L_{h},\beta_{3}\big)_{D^{k}}=\big(\Delta_{(\alpha-2)/2}\rho_{h},\beta_{3}\big)_{D^{k}},\\
&\big(\rho_{h},\delta\big)_{D^{k}}=-\big(\partial_{x}\theta_{h},\partial_{x}\delta\big)_{D^{k}}-\bigg((\partial_{x}\theta_{h})^{*}[\delta]
+\{\partial_{x}\delta\}[\theta_{h}]\bigg)_{k+\frac{1}{2}},\\
&\big(\frac{\partial \theta_{h}}{\partial t},o\big)_{D^{k}}-\lambda_{3}\big(E_{h},o\big)_{D^{k}}-\varpi_{2}\big(p_{h},o\big)_{D^{k}}- \varpi_{1} \big(\upsilon_{h},o\big)_{D^{k}}- \lambda_{4}\big(g(|u_{1}|^{2},|u_{2}|^{2})\upsilon_{h},o\big)_{D^{k}}=0,\\
&\big(E_{h},\beta_{4}\big)_{D^{k}}=\big(\Delta_{(\alpha-2)/2}\xi_{h},\beta_{4}\big)_{D^{k}},\\
&\big(\xi_{h},\omega\big)_{D^{k}}=-\big(\partial_{x}\upsilon_{h},\partial_{x}\omega\big)_{D^{k}}-\bigg((\partial_{x}\upsilon_{h})^{*}[\omega]
+\{\partial_{x}\omega\}[\upsilon_{h}]\bigg)_{k+\frac{1}{2}},\\
\end{split}
\end{equation}
\begin{thm}\label{tt4ha} ($L^{2}$ stability). Suppose $u_{1}(x, t)=p(x, t) + iq(x, t)$ and $u_{2}(x, t) = \upsilon(x, t) + i\theta(x, t)$ and let  $u_{1h}, u_{2h}\in V_{k}^{N}$ be the approximation of $u_{1},u_{2}$ then the solution to the scheme \eqref{1bch2ccf}  satisfies the $L^{2}$ stability
$$
\|u_{1h}\|_{L^{2}(\Omega)}^{2}+\|u_{2h}\|^{2}_{L^{2}(\Omega)}
\leq C(\|u_{1h}(x,0)\|_{L^{2}(\Omega)}^{2}+\|u_{2h}(x,0)\|_{L^{2}(\Omega)}^{2}).$$

\end{thm}
\begin{thm}\label{fgrha}
Let $u_{1}$ and $u_{2}$  be the exact solutions of the linear coupled  fractional
Schr\"{o}dinger equations \eqref{1bchcc}, and let $u_{1h}$ and $u_{2h}$ be the numerical solutions of the DDG scheme \eqref{1bch2ccf}.  Then for small enough $h$, we have the following error estimates:
\begin{equation}\label{tt7h}
\begin{split}
&\|u_{1}(.,T)-u_{1h}(.,T)\|_{L^{2}(\Omega)}+\|u_{2}(.,T)-u_{2h}(.,T)\|_{L^{2}(\Omega)}\leq Ch^{N+1}.\\
\end{split}
\end{equation}
\end{thm}
Theorem  \ref{fgrha} and \ref{tt4ha} can be proven by similar techniques as that in the proof of Theorem  \ref{tt4h} and \ref{fgrh}. We will thus not
give the details here.

\section{Numerical examples}\label{s5n}
\begin{figure}
\begin{center}
  \vspace{4mm}
\begin{tabular}{cc}
\hspace{-2.0cm}
\begin{overpic}[scale=0.8]{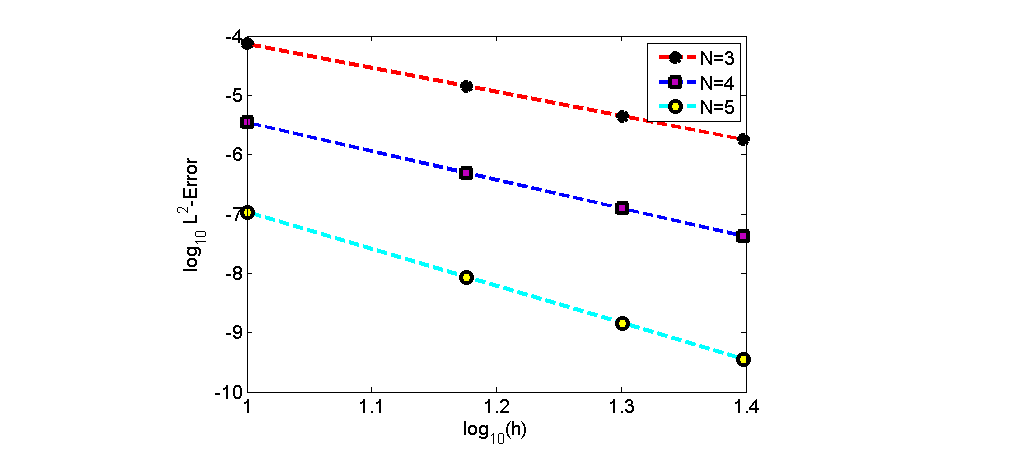}
\end{overpic}

\end{tabular}

\end{center}
\caption[]{\small{ Convergence tests of  \eqref{ex12} with different values of $N$ and $K$.
}} \label{fig:gh2}
\end{figure}

In the following, we present  some numerical experiments to show the accuracy and the performance of the present DDG method for the fractional convection-diffusion and Schr\"{o}dinger type equations. To
deal with the method-of-line fractional PDE, i.e., the classical ODE system, we use the high-order Runge-Kutta time discretizations \cite{Cockburn1999}, when the polynomials are of
degree $N$, a higher-order accurate Runge-Kutta (RK) method must be used in order to guarantee
that the scheme is stable. In this paper we use a fourth-order non-Total variation diminishing (TVD) Runge-Kutta scheme \cite{Gottlieb:1998:TVD:279724.279737}. Numerical experiments demonstrate its numerical stability
\begin{equation}\label{a1}
\begin{split}
\frac{\partial \mathbf{u}_{h}}{\partial t}=\mathcal{F}(\mathbf{u}_{h},t),
\end{split}
\end{equation}
where $\mathbf{u}_{h}$ is the vector of unknowns, we can use the standard fourth-order
four stage explicit RK method (ERK)
\begin{equation}\label{a1}
\begin{split}
&\mathbf{k}^{1}=\mathcal{F}(\mathbf{u}_{h}^{n},t^{n}),\\
&\mathbf{k}^{2}=\mathcal{F}(\mathbf{u}_{h}^{n}+\frac{1}{2}\Delta t\mathbf{k}^{1},t^{n}+\frac{1}{2}\Delta t),\\
&\mathbf{k}^{3}=\mathcal{F}(\mathbf{u}_{h}^{n}+\frac{1}{2}\Delta t\mathbf{k}^{2},t^{n}+\frac{1}{2}\Delta t),\\
&\mathbf{k}^{4}=\mathcal{F}(\mathbf{u}_{h}^{n}+\Delta t\mathbf{k}^{3},t^{n}+\Delta t),\\
&\mathbf{u}_{h}^{n+1}=\mathbf{u}_{h}^{n}+\frac{1}{6}(\mathbf{k}^{1}+2\mathbf{k}^{2}+2\mathbf{k}^{3}+\mathbf{k}^{4}),
\end{split}
\end{equation}
to advance from $\mathbf{u}_{h}^{n}$ to $\mathbf{u}_{h}^{n+1}$, separated by the time step, $\Delta t$. In our examples, the condition $\Delta t\leq C \Delta x^{\alpha}_{min}\,\,\, (0<C<1)$ is used to ensure stability.
\begin{exmp}\label{ex1} Consider the fractional diffusion equation
\begin{equation}\label{91gynb}
\begin{split}
&\frac{\partial u(x,t)}{\partial t}+\varepsilon(-\Delta)^{\frac{\alpha}{2}}u(x,t)=g(x,t),\quad x\in[-1,1],\quad t\in(0,0.5],\\
&u(x,0) = u_{0}(x),
\end{split}
\end{equation}
with the initial condition $u_{0}(x)=(x^{2}-1)^{4}$ and the corresponding forcing term $g(x,t)$  is of the form
\begin{equation}\label{91}
\begin{split}
g(x,t)=e^{-t}\bigg(-u_{0}(x)+\varepsilon(-\Delta)^{\frac{\alpha}{2}} u_{0}(x)\bigg),
\end{split}
\end{equation}
\end{exmp}
then the exact solution is $u(x,t)=e^{-t}(x^{2}-1)^{4}$
with $\varepsilon=\frac{\Gamma(9-\alpha)}{\Gamma(9)}$.\\
We solve the equation  for several different  $\alpha$ and polynomial orders. Moreover, we use   numerical  flux  \eqref{flf} with $\beta_{1}=0$ and  $\beta_{0}$ is a parameter depending on the degree of the approximation polynomial in \cite{liu2010direct}. The errors
and order of convergence are listed in Table \ref{Tab:ab1} and show that the DDG method can achieve the accuracy of order $N +1$ .
\begin{table}[!htb]
    \centering
\begin{center}
 \begin{tabular}{|c|| c c c c c c c |}
  \hline
 \hline
 &\multicolumn{7}{|c|}{$N=1$} \\

 \hline
  K&\multicolumn{7}{|c|}{
 \quad$2^{4}$\qquad\quad\qquad\qquad$2^{5}$\qquad\qquad\qquad\qquad$2^{6}$\qquad\qquad\qquad\qquad$2^{7}$ \qquad\quad\qquad\qquad}\\

  $\alpha$ & $L^{2}$-Error  & $L^{2}$-Error &order & $L^{2}$-Error & order & $L^{2}$-Error &order\\ [0.5ex]

 \hline

1.1& 4.5e-03&7.65e-04 &2.55 &  2.16e-04& 1.82 & 5.07e-05 & 2.09 \\

1.3 & 5.0e-03&1.3e-03&1.98  & 3.42e-04& 1.88 & 8.36e-05 & 2.03
 \\
1.6 &8.4 e-03 & 2.7e-03&1.66& 7.64e-04& 1.78 & 2.11e-04

 & 1.86\\
%

 \hline
 \hline
 &\multicolumn{7}{|c|}{$N=2$} \\

 \hline
  K&\multicolumn{7}{|c|}{
 \quad$2^{4}$\qquad\quad \qquad\qquad $2^{5}$\qquad \qquad \qquad\qquad $2^{6}$\qquad\qquad\qquad\qquad $2^{7}$ \qquad\quad\qquad\qquad}\\

  $\alpha$ & $L^{2}$-Error  & $L^{2}$-Error &order & $L^{2}$-Error & order & $L^{2}$-Error &order\\ [0.5ex]

 \hline

1.1 & 3.3e-03  &  2.89e-04& 3.51 & 2.73e-05 & 3.41
&2.60e-06&3.39 \\

1.3 & 2.5e-03 & 2.03e-04& 3.65  & 1.86e-05 & 3.44
&1.79e-06&3.38 \\
1.6   & 1.2e-03& 1.04e-04& 3.58 &9.21e-06

 & 3.49&8.61e-07&3.42\\
 \hline
 \hline
&\multicolumn{7}{|c|}{$N=3$} \\

 \hline

  K&\multicolumn{7}{|c|}{
 \quad$2^{4}$\qquad\quad \qquad\qquad $2^{5}$\qquad \qquad \qquad\qquad $2^{6}$\qquad\qquad\qquad\qquad $2^{7}$ \qquad\quad\qquad\qquad}\\

  $\alpha$ & $L^{2}$-Error  & $L^{2}$-Error &order & $L^{2}$-Error & order & $L^{2}$-Error &order\\ [0.5ex]

 \hline
1.1 &1.84e-05& 1.12e-06& 4.05&5.55e-08
 & 4.32&2.59e-09&4.42 \\

1.3 & 2.03e-05  &  1.22e-06& 4.06 & 6.46e-08 & 4.24
&3.96e-09&4.03 \\

1.6 & 2.35e-05 & 1.45e-06& 4.02 & 9.38e-08 & 3.95
&6.22e-09&3.92 \\
\hline
 \hline
&\multicolumn{7}{|c|}{$N=4$} \\

 \hline
  K&\multicolumn{7}{|c|}{
 \quad20\qquad\quad\qquad\qquad 25\qquad \qquad \qquad\qquad 30\qquad\qquad\qquad \qquad35\qquad\quad\qquad \qquad}\\

  $\alpha$ & $L^{2}$-Error  & $L^{2}$-Error &order & $L^{2}$-Error & order & $L^{2}$-Error &order\\ [0.5ex]

 \hline
1.1 &3.72e-07& 1.27e-07& 4.82&5.23e-08
 & 4.86&2.47e-08&4.86 \\

1.3 & 3.98e-07  &  1.35e-07& 4.84 & 5.54e-08 & 4.88
&2.60e-08&4.91 \\

1.6 & 4.34e-07 & 1.47e-07& 4.85 & 6.04e-08 & 4.88
&2.84e-08&4.9 \\
%
 \hline
 \hline

\end{tabular}
\end{center}
\caption{ Numerical results for the fractional diffusion equation   in Example \ref{ex1}.}\label{Tab:ab1}
\end{table}
\begin{exmp}\label{ex12} we consider the fractional diffusion equation
\begin{equation}\label{91gyn}
\begin{split}
&\frac{\partial u(x,t)}{\partial t}+\varepsilon(-\Delta)^{\frac{\alpha}{2}}u(x,t)=g(x,t),\quad x\in[0,1],\quad t\in(0,0.5],\\
&u(x,0) = u_{0}(x),
\end{split}
\end{equation}
with the initial condition $u_{0}(x)=x^{11}$ and the Dirichlet boundary conditions $u(0,t)=0,\,\,u(1,t)=e^{-t}$. The corresponding forcing term $g(x,t)$  is of the form
\begin{equation}\label{91}
\begin{split}
g(x,t)=e^{-t}\bigg(-u_{0}(x)+(-\Delta)^{\frac{\alpha}{2}} u_{0}(x)\bigg),
\end{split}
\end{equation}
\end{exmp}
then the exact solution is $u(x,t)=e^{-t}x^{11}$ with $\alpha=1.1,\varepsilon=\frac{\Gamma(12-\alpha)}{\Gamma(12)}$. We consider cases with $N = 3,4, 5$ , $K = 10,15, 20,25$. The numerical orders of convergence are shown in Figure~\ref{fig:gh2}, showing an $O(h^{N+1})$ convergence rate for all orders.
\begin{exmp}\label{ex2}  We consider the fractional Burgers' equation
\begin{equation}\label{91gyn2}
\begin{split}
&\frac{\partial u(x,t)}{\partial t}+\varepsilon(-\Delta)^{\frac{\alpha}{2}}u(x,t)
+\frac{\partial}{\partial x}\bigg(\frac{u^{2}(x,t)}{2}\bigg)=g(x,t),\quad x\in[-1,1],\quad t\in(0,1],\\
&u(x,0) = u_{0}(x),
\end{split}
\end{equation}
with the initial condition $u_{0}(x)=\frac{(x^{2}-1)^4}{100}$ and the corresponding forcing term $g(x,t)$  is of the form
\begin{equation}\label{91}
\begin{split}
g(x,t)=e^{-t}\bigg(-u_{0}(x)+e^{-t}u_{0}(x)u_{0}^{'}(x)u_{0}(x)+\varepsilon(-\Delta)^{\frac{\alpha}{2}}u_{0}(x)\bigg).
\end{split}
\end{equation}
\end{exmp}
In this case, the exact solution will be  $u(x,t)=\frac{e^{-t}(x^{2}-1)^4}{100}$.\\
To complete the scheme, we choose a Lax-Friedrichs flux for the nonlinear term. The problem is solved for several
different values of $\alpha$, polynomial orders $(N)$, and numbers of elements $(K)$.
Table \ref{Tab:a2} shows the numerical $L^{2}$ -Error and the convergence rates of the DDG method with  numerical  flux  \eqref{flf} when $\beta_{0}=1$ and $\beta_{1}=\frac{1}{12}$. From there we see that the DDG method can achieve the accuracy of order $N+1$.
\begin{table}[!htb]
    \centering
\begin{center}
 \begin{tabular}{|c|| c c c c c c c |}
  \hline
 \hline
 &\multicolumn{7}{|c|}{$N=1$} \\

 \hline
  K&\multicolumn{7}{|c|}{
 \quad$10$\qquad\quad \qquad\qquad $20$\qquad \qquad \qquad\qquad $30$\qquad\qquad\qquad\qquad $40$ \qquad\quad\qquad\qquad}\\

  $\alpha$ & $L^{2}$-Error  & $L^{2}$-Error &order & $L^{2}$-Error & order & $L^{2}$-Error &order\\ [0.5ex]

 \hline

1.1& 2.76e-04&3.15e-05 &3.13 &  8.97e-06& 3.1 & 4.68e-06 & 2.26 \\

1.3 & 1.43e-04&3.05e-05&2.23 & 1.34e-05& 2.03& 7.88e-06 & 1.85
 \\
1.6 &1.51e-04  & 5.75e-05&1.5& 2.92e-05& 1.7 & 1.7e-05

 & 1.9\\
%

 \hline
 \hline
 &\multicolumn{7}{|c|}{$N=2$} \\

 \hline
  K&\multicolumn{7}{|c|}{
 \quad$10$\qquad\quad \qquad\qquad $20$\qquad \qquad \qquad\qquad $30$\qquad\qquad\qquad\qquad $40$ \qquad\quad\qquad\qquad}\\

  $\alpha$ & $L^{2}$-Error  & $L^{2}$-Error &order & $L^{2}$-Error & order & $L^{2}$-Error &order\\ [0.5ex]

 \hline

1.1 & 1.35e-04  &  1.82e-05& 2.89 & 5.45e-06 & 2.9717
&2.32e-06&2.98 \\

1.3 & 8.1e-05 & 1.16e-05& 2.78  & 3.57e-06 & 2.91
&1.54e-06&2.92 \\
1.6   & 2.88e-05& 4.43e-06& 2.70 &1.39e-06

 & 2.89&6.07e-07&2.88\\
 \hline
 \hline

\end{tabular}
\end{center}
\caption{ $L^{2}$-Error and order of convergence for Example \ref{ex2} with $K$ elements and polynomial order $N$.}\label{Tab:a2}
\end{table}

\begin{exmp}\label{ex3}  We consider the fractional Burgers' equation
\begin{equation}\label{91gyn2}
\begin{split}
&\frac{\partial u(x,t)}{\partial t}+\varepsilon(-\Delta)^{\frac{\alpha}{2}}u(x,t)
+\frac{\partial}{\partial x}\bigg(\frac{u^{2}(x,t)}{2}\bigg)=g(x,t),\quad x\in[0,1],\quad t\in(0,1],\\
&u(x,0) = u_{0}(x),
\end{split}
\end{equation}
with the initial condition $u_{0}(x)=\frac{x^4}{100}$ and the corresponding forcing term $g(x,t)$  is of the form
\begin{equation}\label{91}
\begin{split}
g(x,t)=e^{-t}\bigg(-u_{0}(x)+e^{-t}u_{0}(x)u_{0}^{'}(x)+\varepsilon(-\Delta)^{\frac{\alpha}{2}}u_{0}(x)\bigg).
\end{split}
\end{equation}
\end{exmp}

In this case, the exact solution will be  $u(x,t)=\frac{e^{-t}x^4}{100}$ and $\varepsilon=\frac{\Gamma(5-\alpha)}{\Gamma(5)}$. Table \ref{Tab:a3} shows the numerical $L^{2}$ -Error and the convergence rates of the DDG method. From there we see that the DDG method can achieve the accuracy of order $N+1$.

\begin{table}[!htb]
    \centering
\begin{center}
 \begin{tabular}{|c|| c c c c c c c |}
  \hline
 \hline

 &\multicolumn{7}{|c|}{$N=2$} \\

 \hline
  K&\multicolumn{7}{|c|}{
 \quad$10$\qquad\quad \qquad\qquad $20$\qquad \qquad \qquad\qquad $30$\qquad\qquad\qquad\qquad $40$ \qquad\quad\qquad\qquad}\\

  $\alpha$ & $L^{2}$-Error  & $L^{2}$-Error &order & $L^{2}$-Error & order & $L^{2}$-Error &order\\ [0.5ex]

 \hline

1.2 & 1.11e-05  &  1.54e-06& 2.85 & 4.54e-07 & 3.01
&1.87e-07&3.09 \\

1.4   & 7.03e-05& 1.02e-05& 2.79 &3.24e-06
 &2.83&1.4e-06&2.9\\
1.6 & 6.55e-05 & 9.77e-06& 2.74  & 3.13e-06 & 2.81
&1.39e-06&2.82 \\
1.8 & 6.21e-05 & 9.39e-06& 2.73  & 3.03e-06& 2.79
&1.37e-06&2.75 \\

 \hline
 \hline
&\multicolumn{7}{|c|}{$N=3$} \\

 \hline
  K&\multicolumn{7}{|c|}{
 \quad20\qquad\quad\qquad\qquad 25\qquad \qquad \qquad\qquad 30\qquad\qquad\qquad \qquad35\qquad\quad\qquad \qquad}\\

  $\alpha$ & $L^{2}$-Error  & $L^{2}$-Error &order & $L^{2}$-Error & order & $L^{2}$-Error &order\\ [0.5ex]

 \hline
1.2 & 2.46e-06  &  2.07e-07& 3.57 & 3.92e-08 &4.1
&1.168e-08&4.21 \\

1.4   & 2.98e-06& 2.33e-07& 3.67 &4.64e-08
 &3.98&1.43e-08&4.09\\
1.6 & 2.9e-06 & 2.28e-07& 3.67  & 4.61e-08 & 2.81
&1.48e-08&3.91 \\
1.8 & 2.71e-06 & 2.06e-07&3.72 & 4.32e-08& 3.85
&1.38e-08&3.98 \\
%
 \hline
 \hline

\end{tabular}
\end{center}
\caption{ $L^{2}$-Error and order of convergence for Example \ref{ex3} with $K$ elements and polynomial order $N$.}\label{Tab:a3}
\end{table}

\begin{figure}
\begin{center}
  \vspace{1mm}
\begin{tabular}{cc}
\hspace{-1.0cm}
\begin{overpic}[scale=0.8]{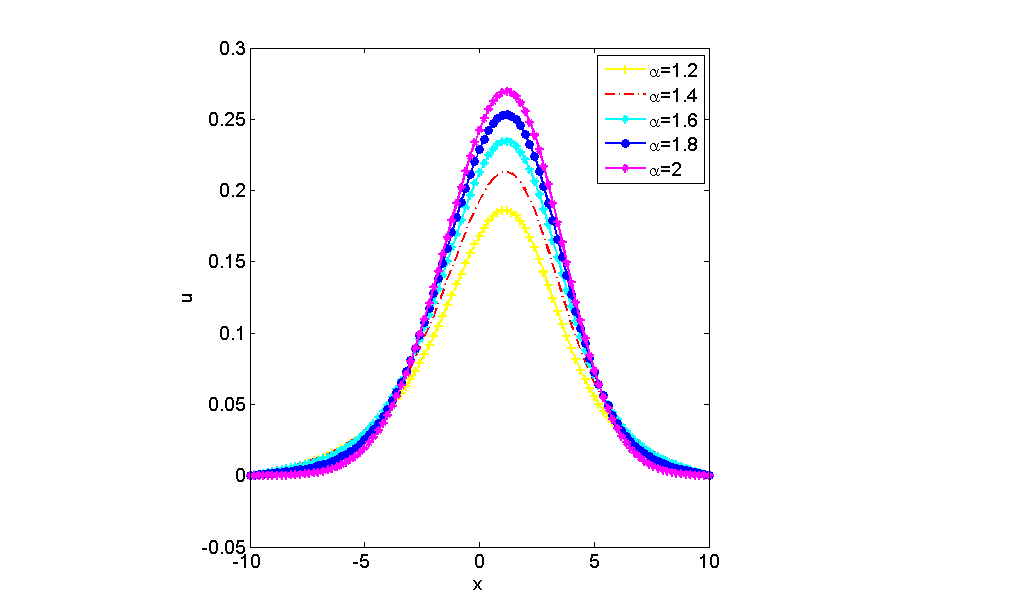}
\end{overpic}

\end{tabular}

\end{center}
\caption[]{\small{  The profile of $u$ with different fractional order $\alpha$ for smooth initial function at $T=3$.
}} \label{fig:bn}
\end{figure}

\begin{figure}
\begin{center}
  \vspace{1mm}
\begin{tabular}{cc}
\hspace{-1.0cm}
\begin{overpic}[scale=0.8]{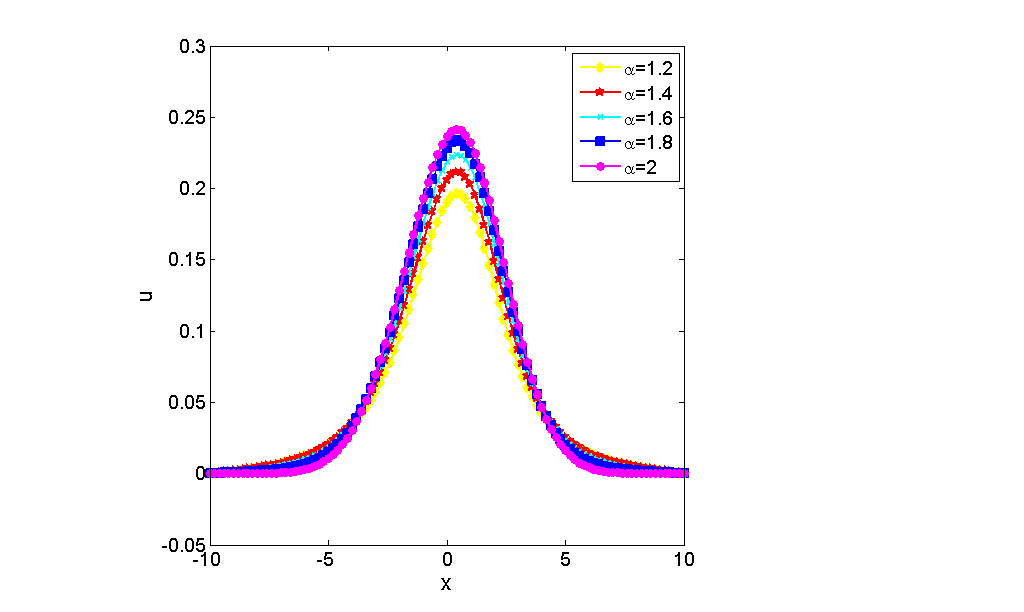}
\end{overpic}

\end{tabular}

\end{center}
\caption[]{\small{  The profile of $u$ with different fractional order $\alpha$ for smooth initial function at $T=3$.
}} \label{fig:fig6}
\end{figure}

\begin{figure}
\begin{center}
  \vspace{5mm}
\begin{tabular}{cc}
\hspace{-4cm}
\begin{overpic}[width=6.0in]{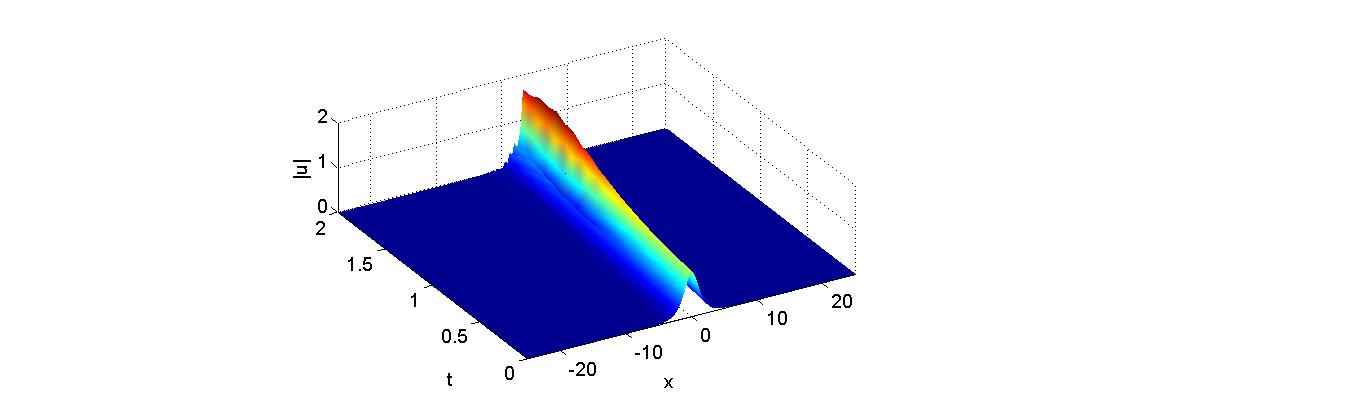}
\put(35,-1.75) {\scriptsize \large{$\alpha=1.4$}}

\end{overpic}
 &\hspace{-5.0cm}
\begin{overpic}[width=6.0in]{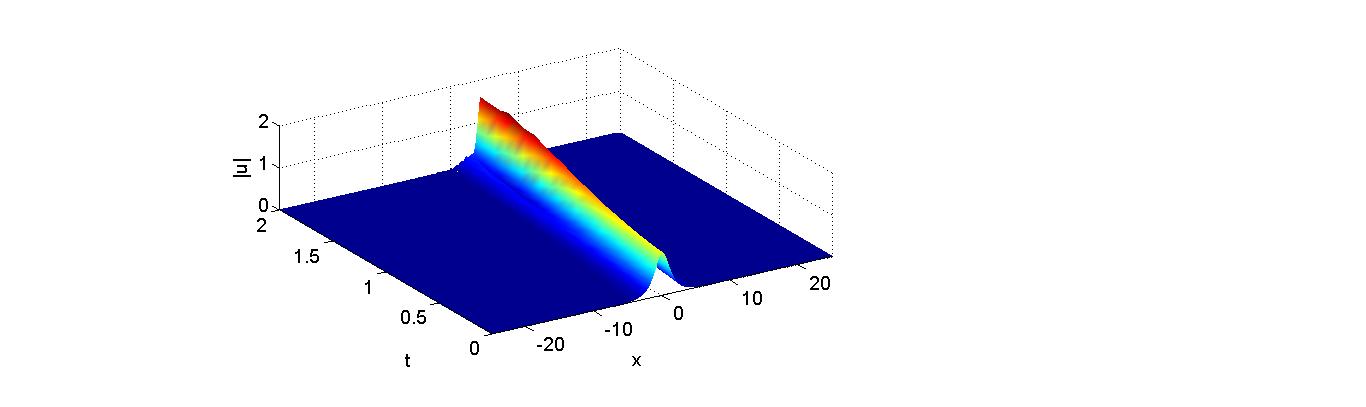}
\put(35,-1.75) {\scriptsize \large{$\alpha=1.6$}}

\end{overpic}

\\
\\
\\
\hspace{-4cm}
\begin{overpic}[width=6in]{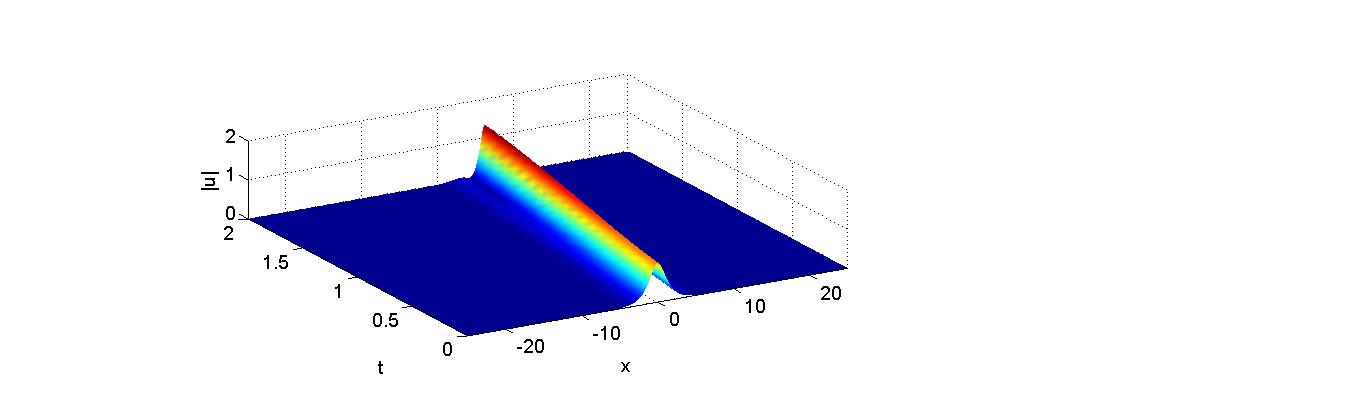}
\put(35,-1.75) {\scriptsize \large{$\alpha=1.8$}}
\end{overpic}
 &\hspace{-5.0cm}
\begin{overpic}[width=6in]{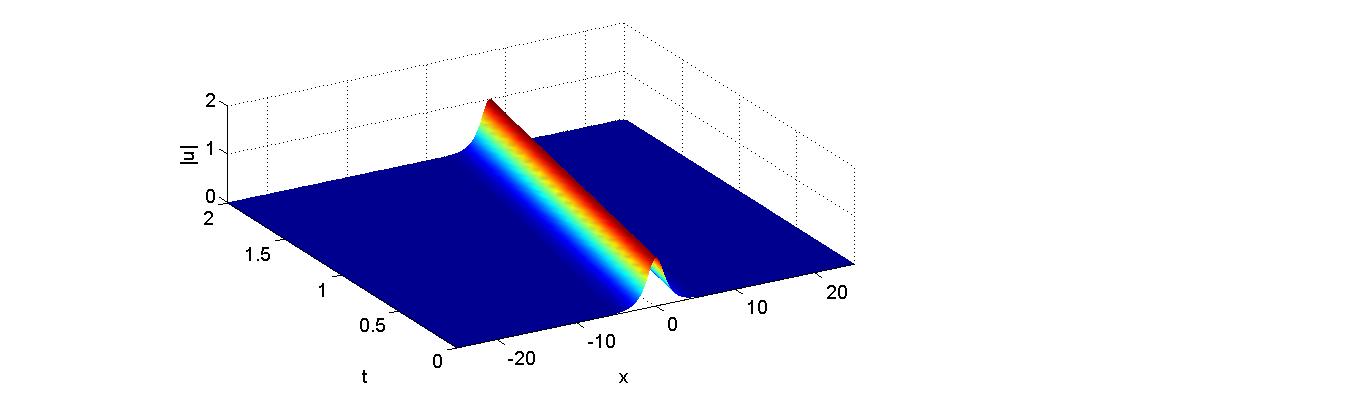}
\put(35,-1.75) {\scriptsize \large{$\alpha=2$}}

\end{overpic}

\end{tabular}
\caption[]{{ Numerical results for the nonlinear fractional Schr\"{o}dinger equation   in Example \ref{ex4n}.}}\label{fig:fig2}
\end{center}
\end{figure}

\begin{exmp}\label{ex233}  We consider the fractional convection-diffusion equation with a discontinuous initial condition,
$$
u(x,0)=\left\{
           \begin{array}{lll}
             x+1, & \hbox{$-1\leq x<0$;} \\
             2x, & \hbox{$0\leq x\leq1$;}\\
                          0, & \hbox{otherwise.}
           \end{array}
         \right.
$$
\end{exmp}
We consider \eqref{25n} with parameters $ \varepsilon=1$, $x\in[-10,10]$  and  solve
the equation for several different values of $\alpha$. The numerical solution $u_{h}(x, t)$ for $\alpha=1.2,\,1.4,\,1.6,\,1.8,\,2.0$ is shown in figure \ref{fig:bn}. From this figure  clear that the dissipative effect increases with $\alpha$ and the classical case with $\alpha = 2$ is a limit of the fractional case.
\begin{exmp}  \label{ex4nm} We consider the fractional convection-diffusion  equation \eqref{25n} with
 initial condition,
\begin{equation}\label{91}
\begin{split}
u(x,0)=e^{-2x^{2}},
\end{split}
\end{equation}
\end{exmp}
 with parameters   $\varepsilon=1$, $\,x\in[-10,10]$. We consider cases with $N =  2$ and $K = 50$ and  solve
the equation for several different values of $\alpha$. The numerical solution $u_{h}(x, t)$ for $\alpha=1.2,\,1.4,\,1.6,\,1.8,\,2.0$ is shown in Figure \ref{fig:fig6}.  We observe that
the order $\alpha$ will affect the shape of the soliton case. This property of the fractional convection-diffusion equation can be used in physics to modify the shape of wave without change of the nonlinearity and dispersion effects. The numerical solutions of the fractional equation are convergent to the solutions of the classical non-fractional equation when $\alpha$ tends to $2$.
\begin{exmp}  \label{ex4cz} We consider the nonlinear fractional Schr\"{o}dinger equation
\begin{equation}\label{sch1}
\begin{split}
&i\frac{\partial u}{\partial t}- \varepsilon(-\Delta)^{\frac{\alpha}{2}}u+|u|^{2}u=g(x,t),\quad x\in[-1,1],\quad t\in(0,0.5],\\
&u(x,0) = u_{0}(x),
\end{split}
\end{equation}

with the initial condition $u_{0}(x)=(x^{2}-1)^{5}$ and the corresponding forcing term $g(x,t)$  is of the form
\begin{equation}\label{91}
\begin{split}
g(x,t)=e^{-it}\bigg(iu_{0}(x)-\varepsilon(-\Delta)^{\frac{\alpha}{2}}u_{0}(x)+(u_{0}(x))^{3}\bigg),
\end{split}
\end{equation}
\end{exmp}
to obtain an exact solution $u(x,t)=e^{-it}(x^{2}-1)^{5}$ with  $\varepsilon=\frac{\Gamma(11-\alpha)}{\Gamma(11)}$. We solve the equation for several diﬀerent α and polynomial orders. The errors and order of convergence are listed in Table \ref{Tab:bzz2}, confirming optimal $O(h^{N+1})$ order of
convergence across $1<\alpha<2$.
\begin{table}[!htb]
    \centering
\begin{center}
 \begin{tabular}{|c|| c c c c c c c |}
  \hline
 \hline
 &\multicolumn{7}{|c|}{$N=1$} \\

 \hline
  K&\multicolumn{7}{|c|}{
 \quad$2^{4}$\qquad\quad\qquad\qquad$2^{5}$\qquad\qquad\qquad\qquad$2^{6}$\qquad\qquad\qquad\qquad$2^{7}$ \qquad\quad\qquad\qquad}\\

  $\alpha$ & $L^{2}$-Error  & $L^{2}$-Error &order & $L^{2}$-Error & order & $L^{2}$-Error &order\\ [0.5ex]

 \hline

1.1& 1.6e-03&3.49e-04&2.19 &  6.92e-05& 2.33 & 1.34e-05 & 2.37 \\

1.3 & 6.2e-03&1.6e-03&1.94  & 3.73e-04& 2.12 & 8.29e-05
 & 2.17
 \\
1.6 &9.0e-03 & 2.9e-03&1.65& 8.12e-04& 1.83 & 1.94e-04

 & 2.07\\

 \hline
 \hline
 &\multicolumn{7}{|c|}{$N=2$} \\

 \hline
  K&\multicolumn{7}{|c|}{
 \quad$2^{4}$\qquad\quad \qquad\qquad $2^{5}$\qquad \qquad \qquad\qquad $2^{6}$\qquad\qquad\qquad\qquad $2^{7}$ \qquad\quad\qquad\qquad}\\

  $\alpha$ & $L^{2}$-Error  & $L^{2}$-Error &order & $L^{2}$-Error & order & $L^{2}$-Error &order\\ [0.5ex]

 \hline

1.1 & 5.52e-04 & 5.75e-05& 3.27 & 1.03e-05 & 2.48
&5.44e-07&4.25 \\

1.3 & 6.23e-04 & 7.68e-05& 3.02  & 1.04e-05 & 2.88
&1.31e-06&3.0 \\
1.6   & 6.0e-04& 7.82e-05& 2.94 &9.67e-06

 & 3.02&1.18e-06&3.04\\
 \hline
 \hline
&\multicolumn{7}{|c|}{$N=3$} \\
\hline
  K&\multicolumn{7}{|c|}{
 \quad$2^{4}$\qquad\quad \qquad\qquad $2^{5}$\qquad \qquad \qquad\qquad $2^{6}$\qquad\qquad\qquad\qquad $2^{7}$ \qquad\quad\qquad\qquad}\\

  $\alpha$ & $L^{2}$-Error  & $L^{2}$-Error &order & $L^{2}$-Error & order & $L^{2}$-Error &order\\ [0.5ex]

 \hline
1.1 &2.47e-05& 1.62e-06& 3.93&8.86e-08
 & 4.19&5.52e-09&4.0 \\

1.3 & 2.50e-05  &  1.62e-06&3.95 &1.016e-07 & 4.0
&6.142e-09&4.05 \\

1.6 & 2.5e-05 & 1.51e-06& 4.06 & 9.52e-08 & 3.98
&5.43e-09&4.13 \\
\hline
 \hline

\end{tabular}
\end{center}
\caption{ Numerical results for the nonlinear fractional Schr\"{o}dinger equation   in Example \ref{ex4cz}.}\label{Tab:bzz2}
\end{table}
\begin{exmp}  \label{ex4n} We consider the nonlinear fractional Schr\"{o}dinger equation \eqref{sch1vn} with
 initial conditions:\\
 $(a)$ Single soliton: We consider the initial condition
\begin{equation}\label{91}
\begin{split}
u(x,0)=e^{2i(x-x_{0})}sech(x-x_{0}),
\end{split}
\end{equation}
\end{exmp}
 with parameters   $\varepsilon_{1}=\varepsilon_{2}=2$, $x_{0} = 0$ and $\,x\in[-25,25]$. We consider cases with $N =  2$ and $K = 200$ and  solve
the equation for several different values of $\alpha$. The numerical solution $u_{h}(x, t)$ for $\alpha=1.4,\,1.6,\,1.8,\,2.0$ is shown in Figure \ref{fig:fig2}.  We observe that
the order $\alpha$ will affect the shape of the soliton case.  When $\alpha$ becomes smaller, the shape of the soliton will change more quickly. This property of the fractional Schr\"{o}dinger equation can be used in physics to modify the shape of wave without change of the nonlinearity and dispersion effects. The numerical solutions of the fractional equation are convergent to the solutions of the classical non-fractional equation when $\alpha$ tends to $2$.

\begin{figure}
\begin{center}
  \vspace{5mm}
\begin{tabular}{cc}
\hspace{-4cm}
\begin{overpic}[width=6.0in]{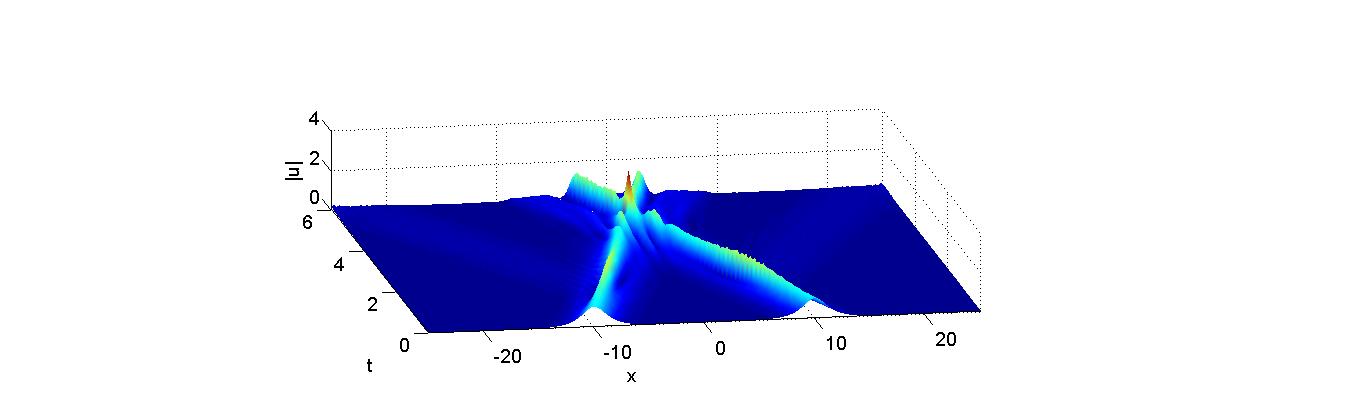}
\put(39,-1.75) {\scriptsize \large{$\alpha=1.6$}}

\end{overpic}
 &\hspace{-5.0cm}
\begin{overpic}[width=6.0in]{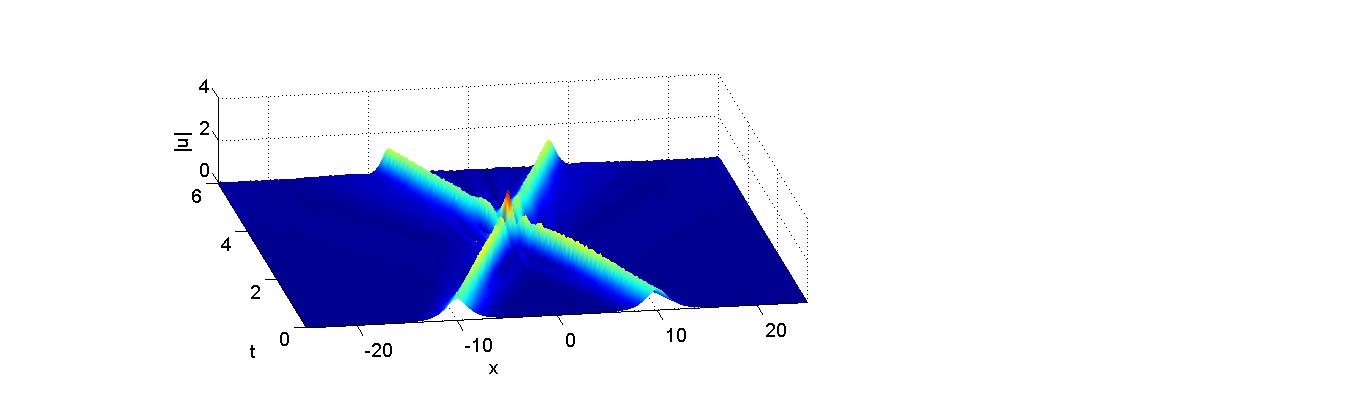}
\put(35,-1.75) {\scriptsize \large{$\alpha=1.8$}}

\end{overpic}

\\
\\
\\
\hspace{-4cm}
\begin{overpic}[width=6in]{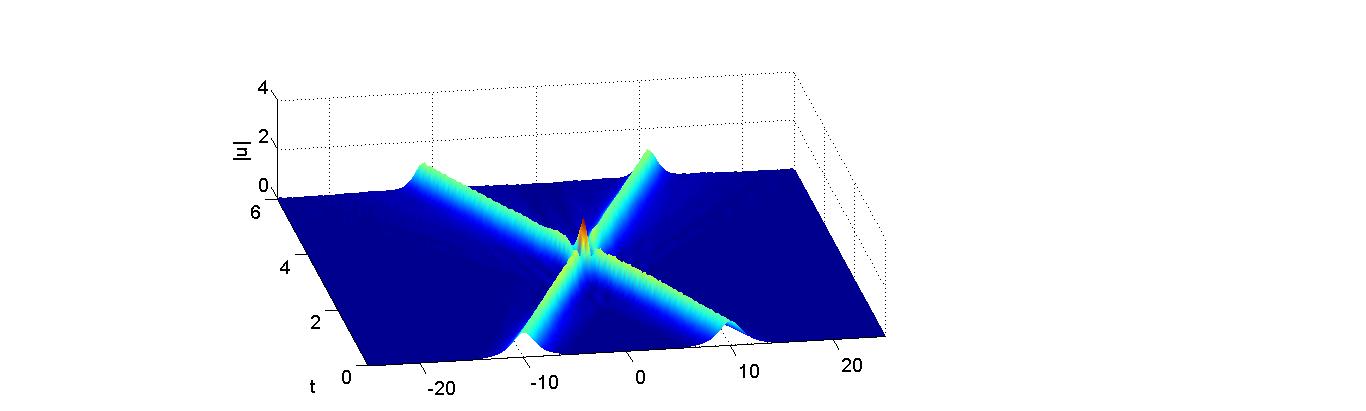}
\put(35,-1.75) {\scriptsize \large{$\alpha=1.9$}}
\end{overpic}
 &\hspace{-5.0cm}
\begin{overpic}[width=6in]{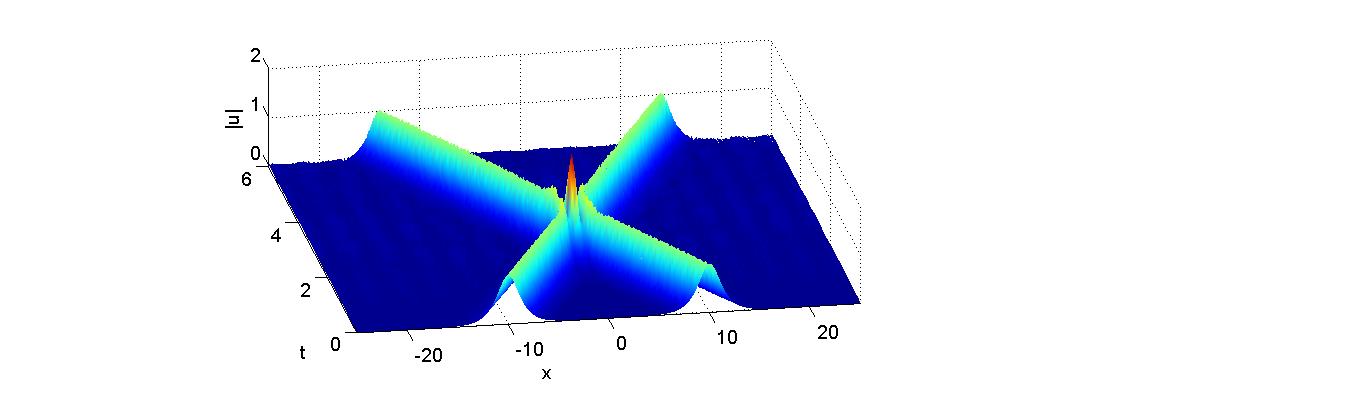}
\put(35,-1.75) {\scriptsize \large{$\alpha=2$}}

\end{overpic}

\end{tabular}
\caption[]{{  The double soliton collision of \eqref{sch1vn} with initial condition \eqref{91za}, $c_{1} =4,x_{1} = -10, c_{2} = -4, x_{2} = 10$.}}\label{fig:cbzz}
\end{center}
\end{figure}

\newpage
$(b)$ Interaction of two solitons: To study the interaction of two solitons, we will take equation \eqref{sch1vn} with the initial condition
\begin{equation}\label{91za}
\begin{split}
u(x,0)=\sum _{j=1}^{2}e^{\frac{1}{2}ic_{j}(x-x_{j})}sech(x-x_{j}),
\end{split}
\end{equation}
with parameters   $\varepsilon_{1}=1,\,\varepsilon_{2}=2$ and $\,x\in[-25,25]$. We consider cases with $N =  2$ and $K = 200$ and  solve
the equation for several different values of $\alpha$. The numerical solution $u_{h}(x, t)$ for $\alpha=1.6,\,1.8,\,1.9,\,2.0$ is shown in Figure \ref{fig:cbzz}. We observe that
the order $\alpha$ will affect the shape of the two solitons case and the two waves approach each other interact and leave the interaction unchanged
in the shape and velocity. In addition, the interaction is strictly elastic because each of them recovers its exact initial shape after they pass
through each other.
\begin{exmp}\label{ex6} We consider the nonlinear  coupled fractional Schr\"{o}dinger equations
\begin{equation}\label{sch1}
\begin{split}
&i\frac{\partial u_{1}(x,t)}{\partial t}- \varepsilon_{1}(-\Delta)^{\frac{\alpha}{2}}u_{1}(x,t)+u_{2}(x,t)+u_{1}(x,t)+(|u_{1}(x,t)|^{2}+|u_{2}(x,t)|^{2})
u_{1}(x,t)=g_{1}(x,t),\,\, x\in[0,1],\,\,t\in(0,0.5],\\
&i\frac{\partial u_{2}(x,t)}{\partial t}- \varepsilon_{2}(-\Delta)^{\frac{\alpha}{2}}u_{2}(x,t)+u_{2}(x,t)+u_{1}(x,t)+(|u_{1}(x,t)|^{2}
+|u_{2}(x,t)|^{2})u_{2}(x,t)=g_{2}(x,t),\, x\in[0,1],\, t\in(0,0.5],\\
\end{split}
\end{equation}
 and the corresponding forcing terms $g_{1}(x,t)$ and $g_{2}(x,t)$  are of the form
\begin{equation}\label{91}
\begin{split}
&g_{1}(x,t)=e^{-it}\bigg(iu_{1}(x,0)-\varepsilon_{1}(-\Delta)^{\frac{\alpha}{2}}u_{1}(x,0)+u_{2}(x,0)+u_{1}(x,0)
+(|u_{1}(x,0)|^{2}+|u_{2}(x,0)|^{2})u_{1}(x,0)\bigg),\\
&g_{2}(x,t)=e^{-it}\bigg(iu_{2}(x,0)-\varepsilon_{2}(-\Delta)^{\frac{\alpha}{2}}u_{2}(x,0)+u_{2}(x,0)+u_{1}(x,0)
+(|u_{1}(x,0)|^{2}+|u_{2}(x,0)|^{2})u_{2}(x,0)\bigg),
\end{split}
\end{equation}
\end{exmp}
to obtain an exact solutions $u_{1}(x,t)=e^{-it}x^{5}$ and $u_{2}(x,t)=e^{-it}x^{5}$
with $\alpha=1.1,\,\varepsilon_{1}=\frac{\Gamma(6-\alpha)}{2\Gamma(6)}$, $\varepsilon_{2}=\frac{\Gamma(6-\alpha)}{2\Gamma(6)}$. The errors
and order of convergence are listed in Tables \ref{Tab:e1} and \ref{Tab:e2}, confirming optimal $O(h^{N+1})$ order of
convergence across.
\begin{table}[!htb]
    \centering
\begin{center}
 \begin{tabular}{|c|| c c||c|| c c ||c||c c |}
  \hline
 \hline

 \hline
  N&\multicolumn{8}{|c|}{   N=1\quad \quad \quad\quad\quad \quad\qquad \qquad N=2\quad\quad\quad \qquad\quad \quad  \qquad \quad N=3}\\
\hline\hline
  K & $L^{2}$-Error & order & K &$L^{2}$-Error & order&K&$L^{2}$-Error& order \\ [0.5ex]

 \hline
10 & 3.95e-03& -&       10&4.7e-04&-&         10&1.02e-04&-  \\

20 & 1.02e-03 &1.95&   20&8.95e-05&3.19 &    20&5.91e-06&4.11 \\
 40&2.21e-04&2.20&     40&1.05e-05&3.09&    40&3.82e-07&3.95\\

 \hline
 \hline

\end{tabular}
\end{center}
\caption{ $L^{2}$-Error and order of convergence for $u_{1}$ with $K$ elements and polynomial order $N$.}\label{Tab:e1}
\end{table}

\begin{table}[!htb]
    \centering
\begin{center}
 \begin{tabular}{|c|| c c||c|| c c ||c||c c |}
  \hline
 \hline

 \hline
  N&\multicolumn{8}{|c|}{   N=1\quad \quad \quad\quad\quad \quad\qquad \qquad N=2\quad\quad\quad \qquad\quad \quad  \qquad \quad N=3}\\
\hline\hline
  K & $L^{2}$-Error & order & K &$L^{2}$-Error & order&K&$L^{2}$-Error& order \\ [0.5ex]

 \hline
10 & 4.32e-03 & -&       10&4.18e-04&-&         10&2.18e-04&-  \\

20 & 1.17e-03 &1.89&   20&9.04e-05&3.24 &    20&1.32e-05&4.05\\
 40&2.69e-04&2.23&     40&1.1e-05&3.04&    40&8.57e-07&3.95\\

 \hline
 \hline

\end{tabular}
\end{center}
\caption{ $L^{2}$-Error and order of convergence for $u_{2}$ with $K$ elements and polynomial order $N$.}\label{Tab:e2}
\end{table}

\begin{figure}
\begin{center}
  \vspace{6mm}
\begin{tabular}{cc}
\hspace{-13cm}
\begin{overpic}[scale=0.25]{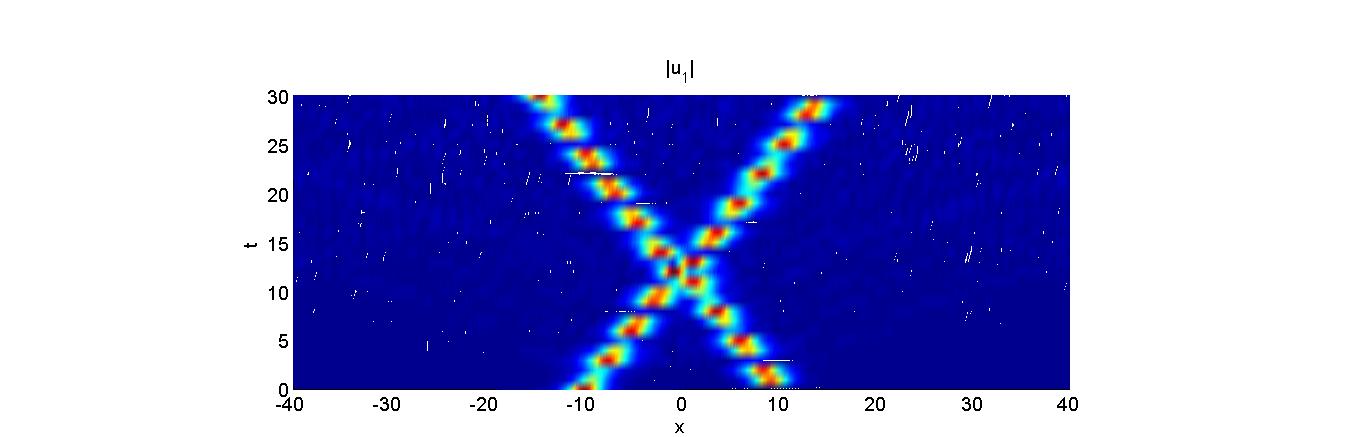}
\end{overpic}
 &\hspace{-12.0cm}
\begin{overpic}[scale=0.25]{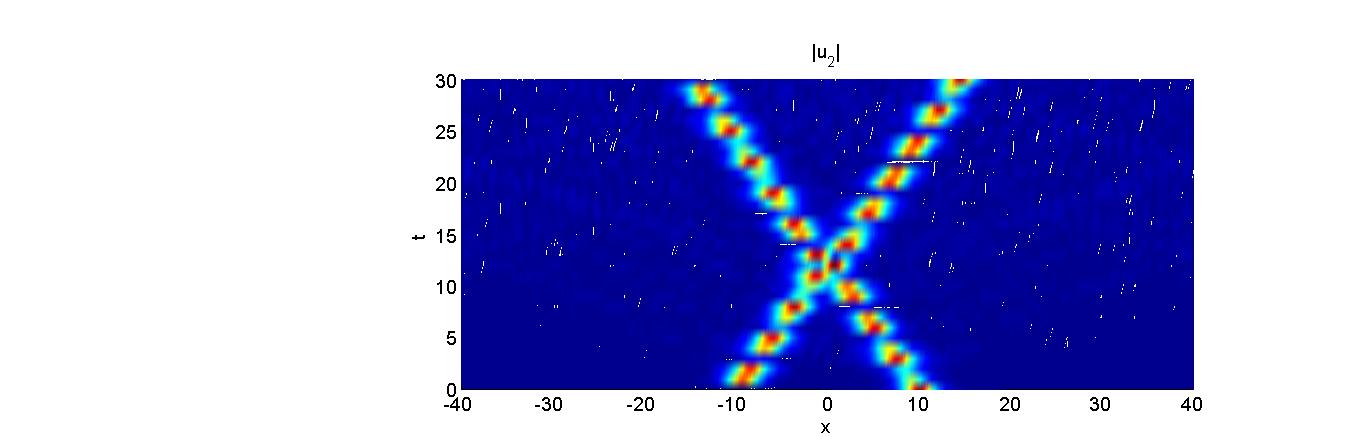}
\end{overpic}
\\
\\
\\
\hspace{-1cm}
\begin{overpic}[scale=0.4]{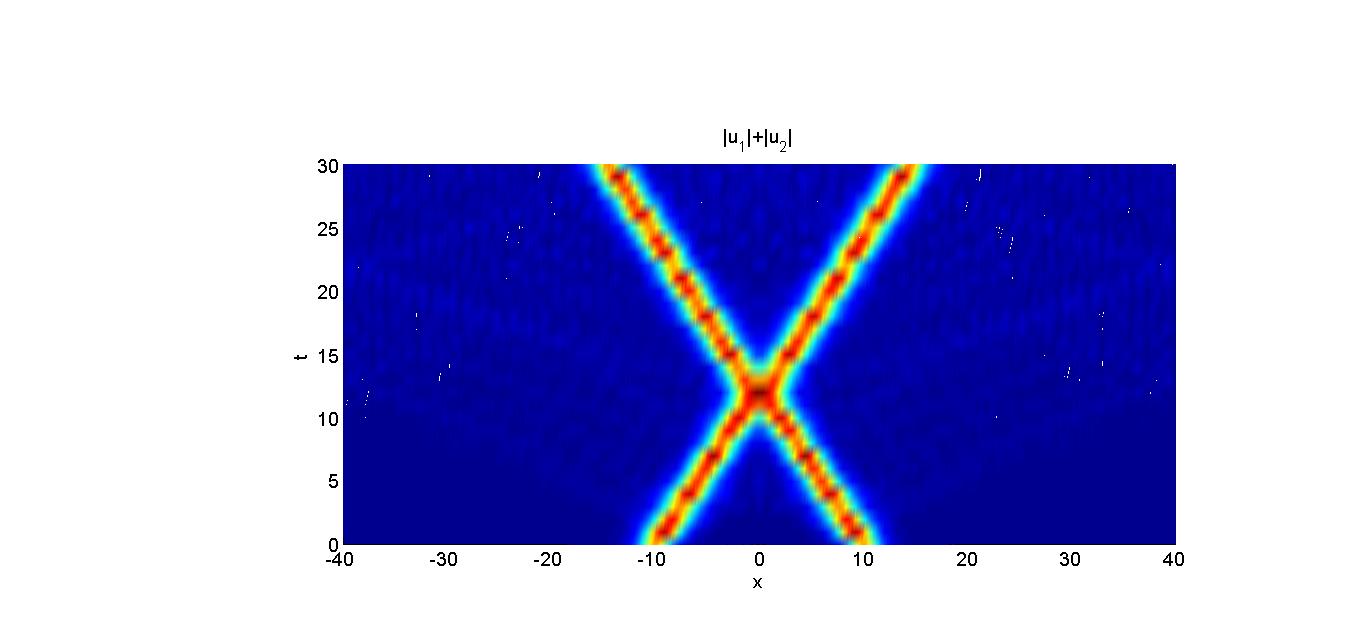}
\end{overpic}

\end{tabular}
\caption[]{{   The numerical simulation of the two soliton waves for Example \ref{ex12vng} with $\varpi_{1}= 1$ , $\alpha = 2$.}}\label{fig:elastic3}
\end{center}
\end{figure}

\begin{figure}
\begin{center}
  \vspace{6mm}
\begin{tabular}{cc}
\hspace{-12.5cm}
\begin{overpic}[scale=0.25]{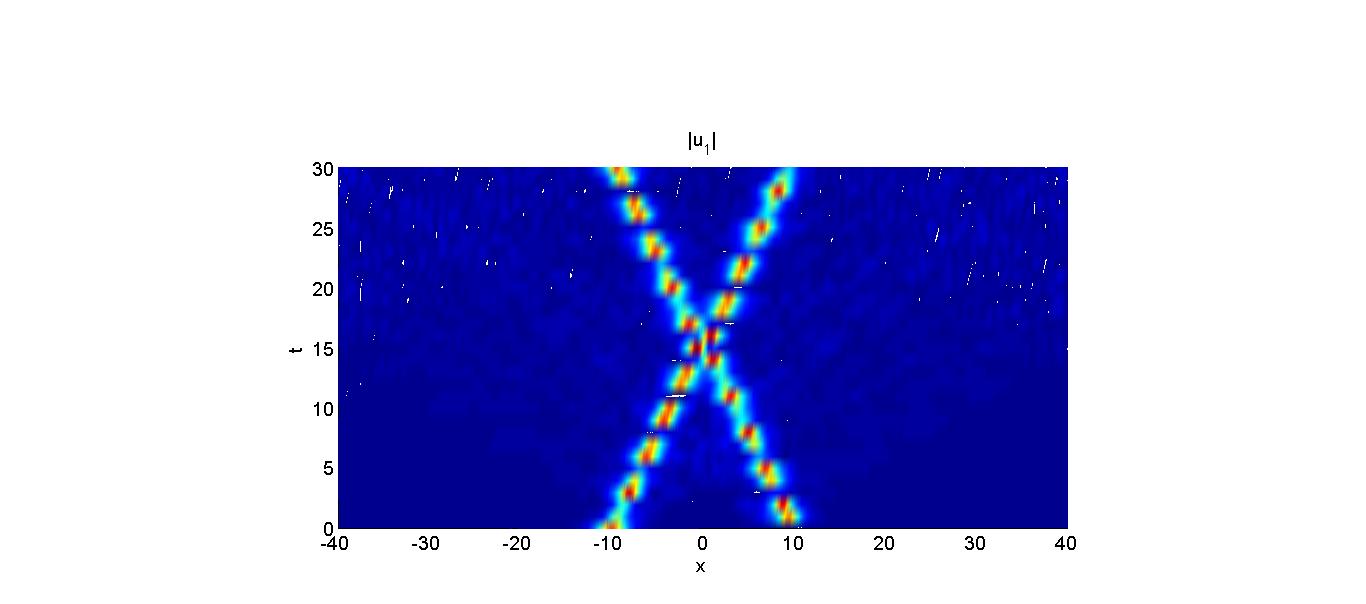}
\end{overpic}
 &\hspace{-12.0cm}
\begin{overpic}[scale=0.25]{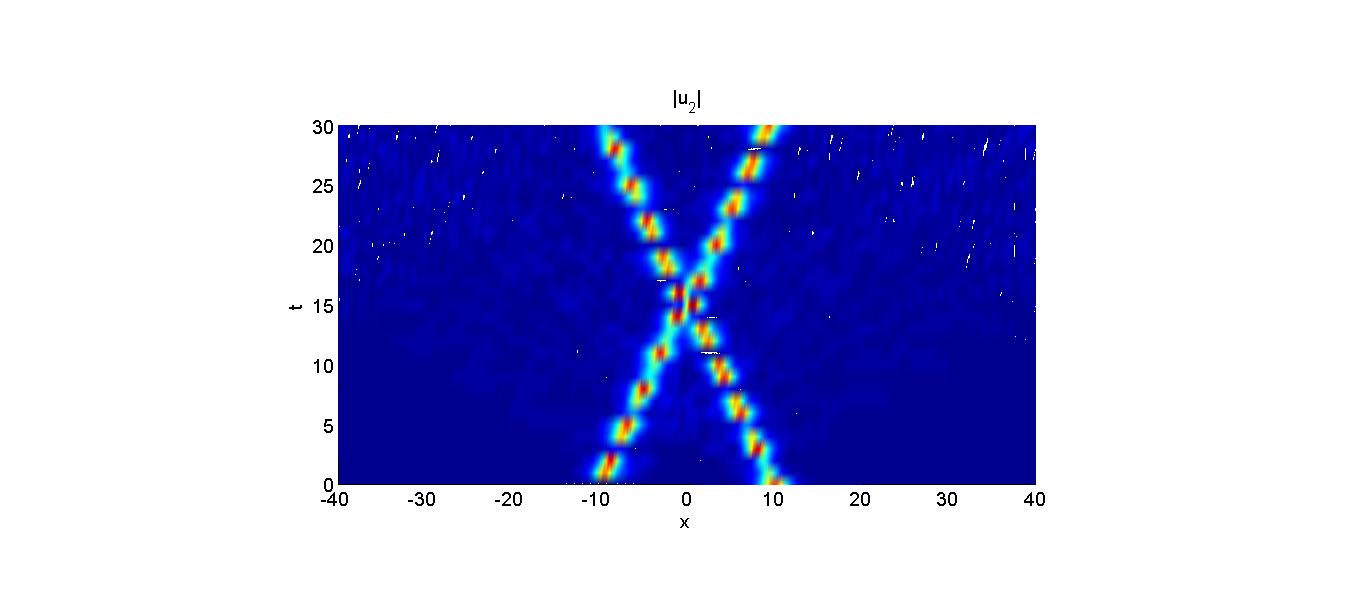}
\end{overpic}
\\
\hspace{-1cm}
\begin{overpic}[scale=0.4]{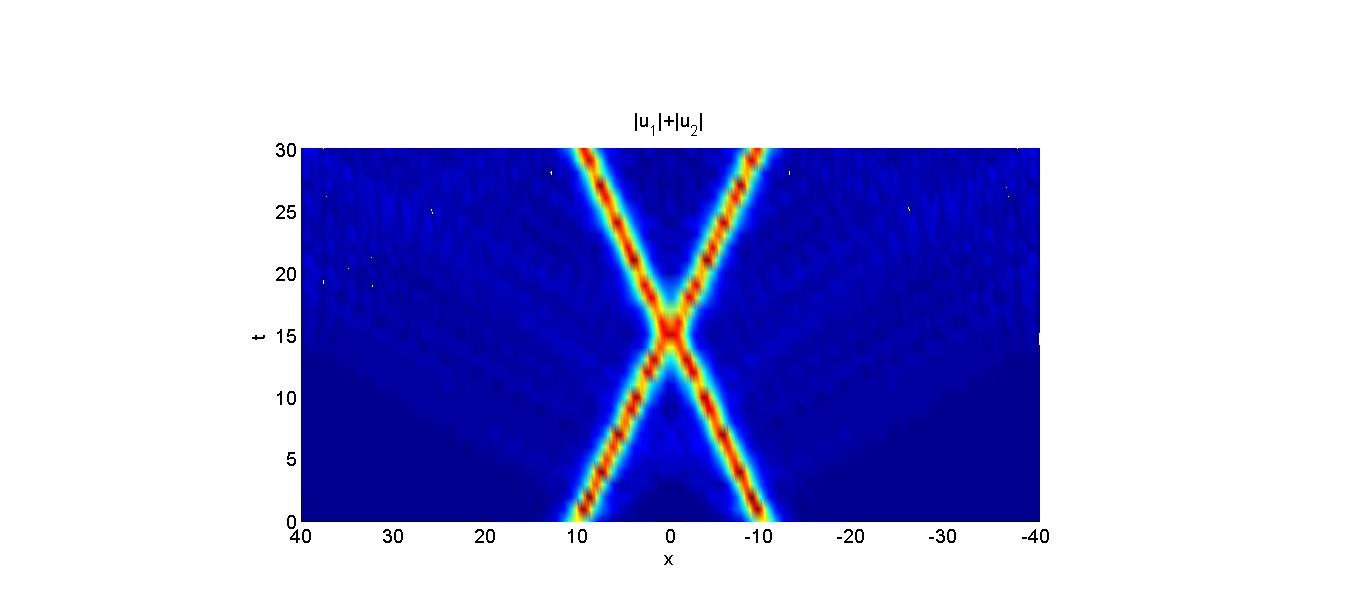}
\end{overpic}

\end{tabular}
\caption[]{{   The numerical simulation of the two soliton waves for Example \ref{ex12vng} with $\varpi_{1}= 1$ , $\alpha = 1.6$. }}\label{fig:elastic4}
\end{center}
\end{figure}

\begin{figure}
\begin{center}
  \vspace{6mm}
\begin{tabular}{cc}
\hspace{-13cm}
\begin{overpic}[scale=0.25]{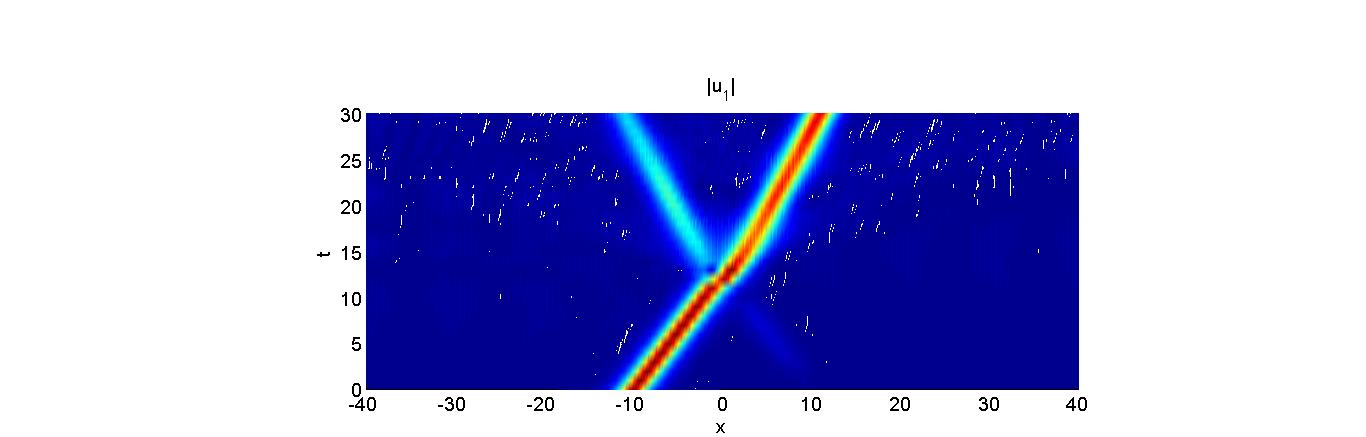}
\end{overpic}
 &\hspace{-12.0cm}
\begin{overpic}[scale=0.25]{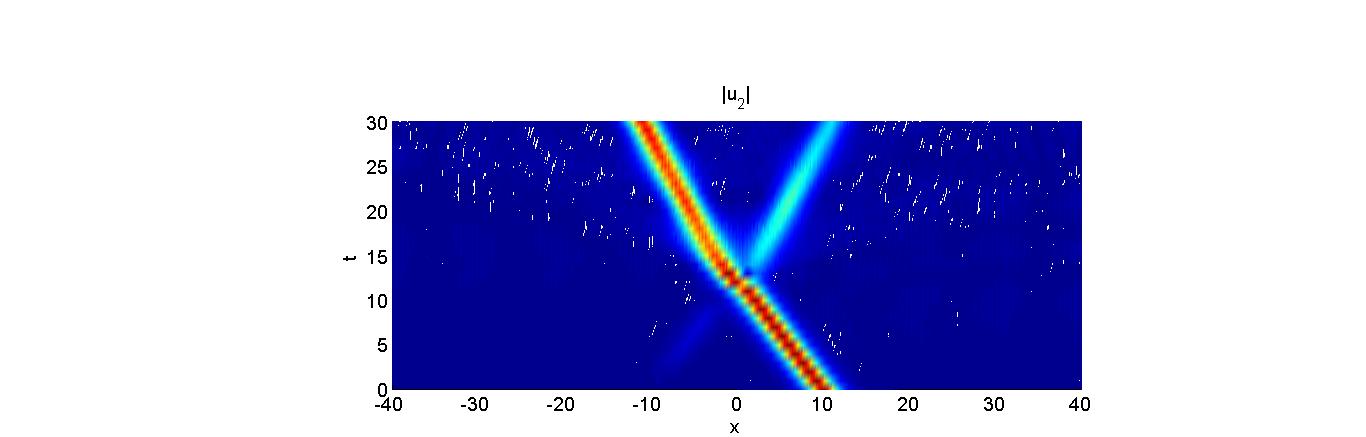}
\end{overpic}
\\
\hspace{-1cm}
\begin{overpic}[scale=0.4]{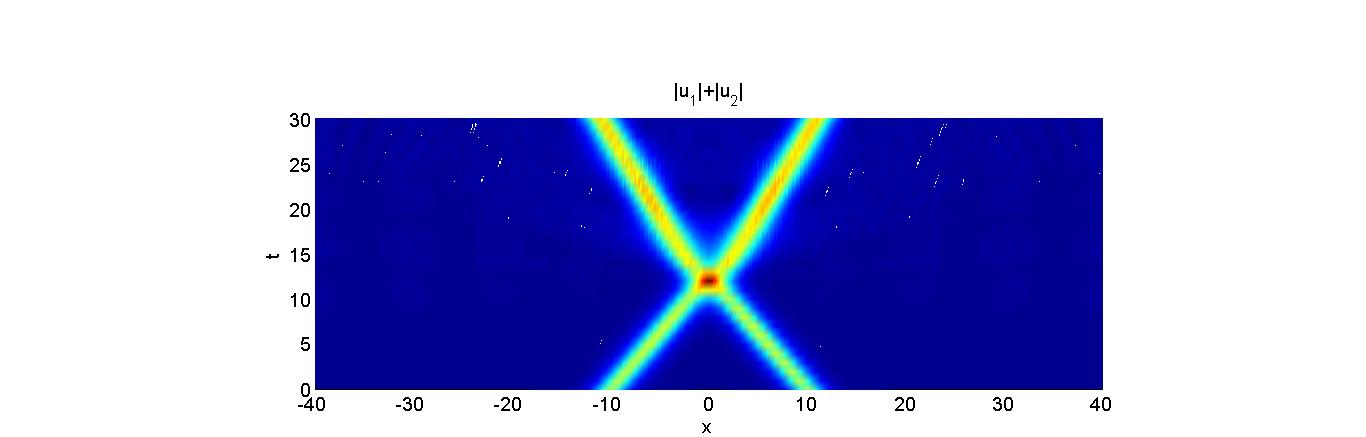}
\end{overpic}

\end{tabular}
\caption[]{{   The numerical simulation of the two soliton waves for Example \ref{ex12vng}  with $\varpi_{1}= 0.0175$ and $\alpha = 2$.}}\label{fig:inelastic6}
\end{center}
\end{figure}

\begin{figure}
\begin{center}
  \vspace{6mm}
\begin{tabular}{cc}
\hspace{-13cm}
\begin{overpic}[scale=0.25]{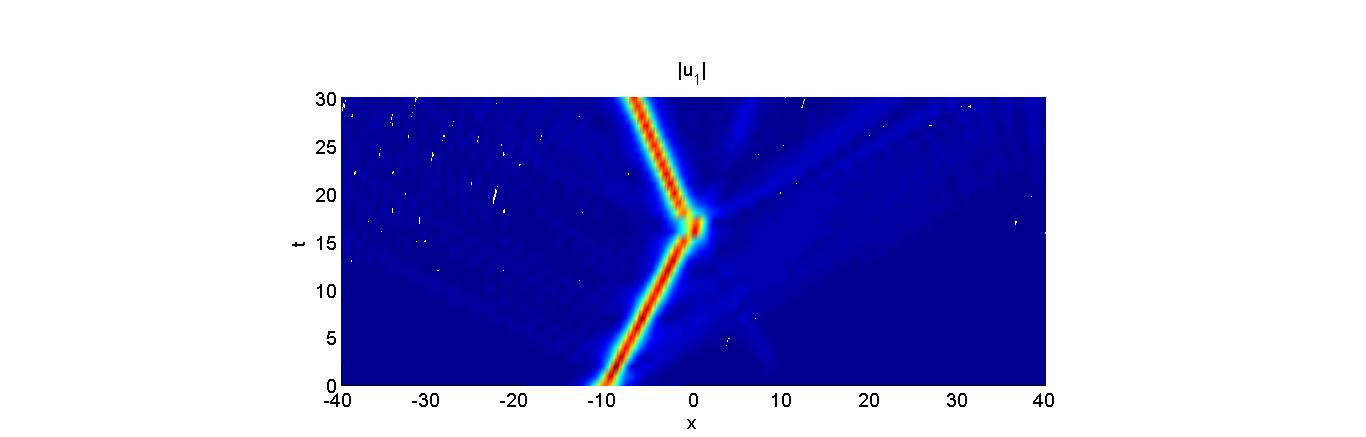}
\end{overpic}
 &\hspace{-12.5cm}
\begin{overpic}[scale=0.25]{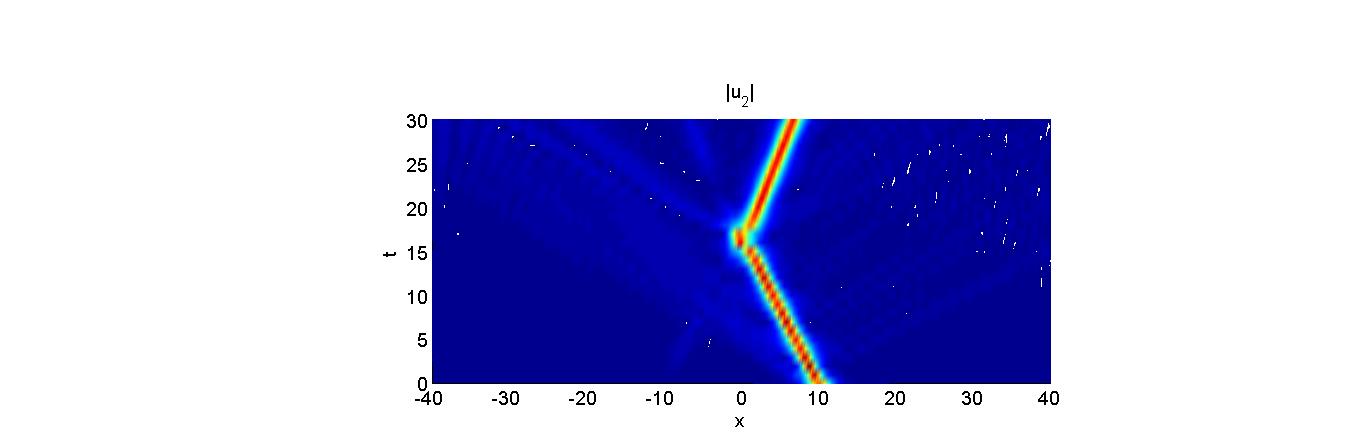}
\end{overpic}
\\
\hspace{-1cm}
\begin{overpic}[scale=0.4]{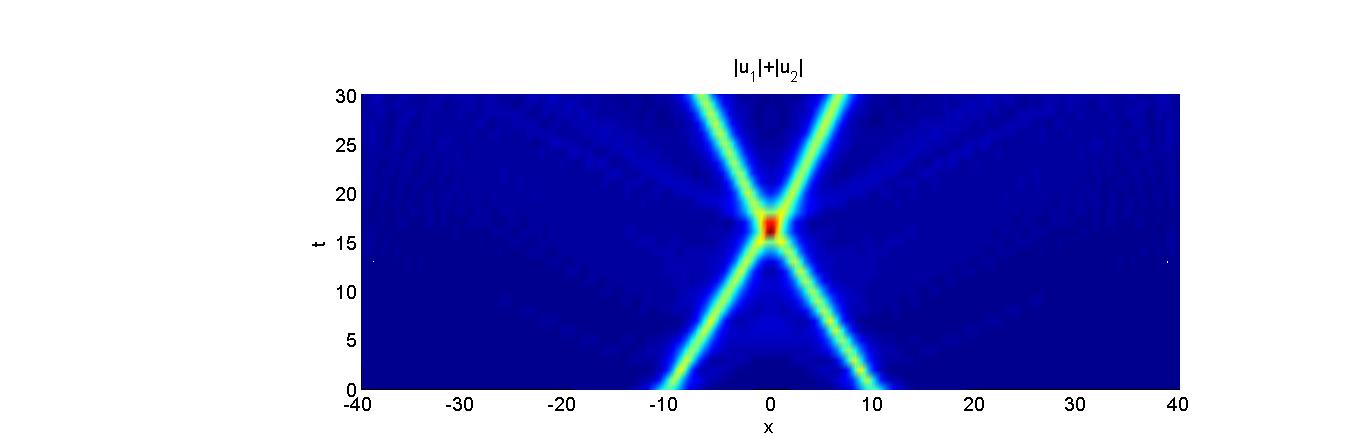}
\end{overpic}

\end{tabular}
\caption[]{{   The numerical simulation of the two soliton waves for Example \ref{ex12vng}  with $\varpi_{1}= 0.0175$ and $\alpha = 1.6$.}}\label{fig:inelastic6h}
\end{center}
\end{figure}

\begin{figure}
\begin{center}
  \vspace{6mm}
\begin{tabular}{cc}
\hspace{-13cm}
\begin{overpic}[scale=0.25]{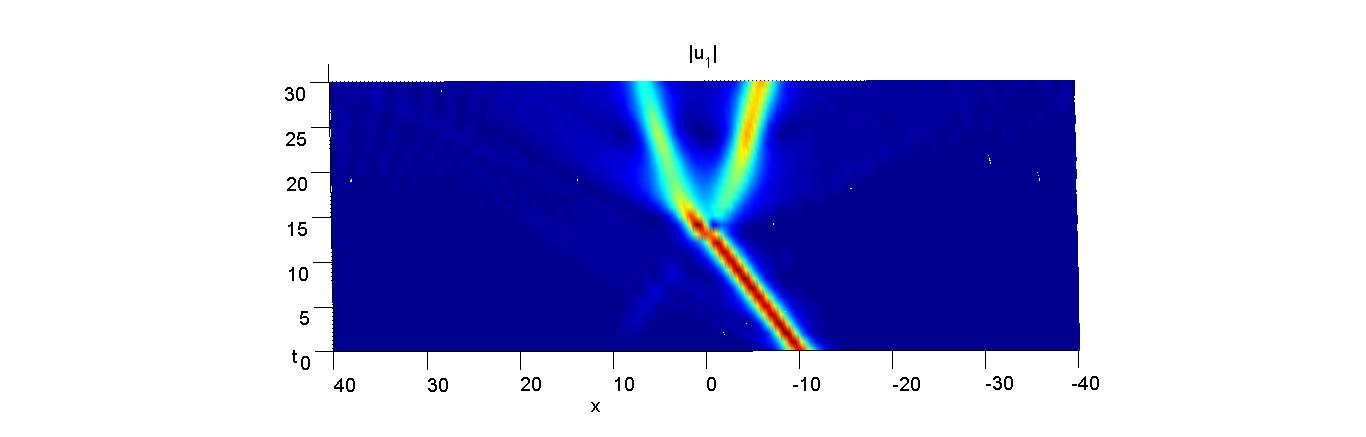}
\end{overpic}
 &\hspace{-12.0cm}
\begin{overpic}[scale=0.25]{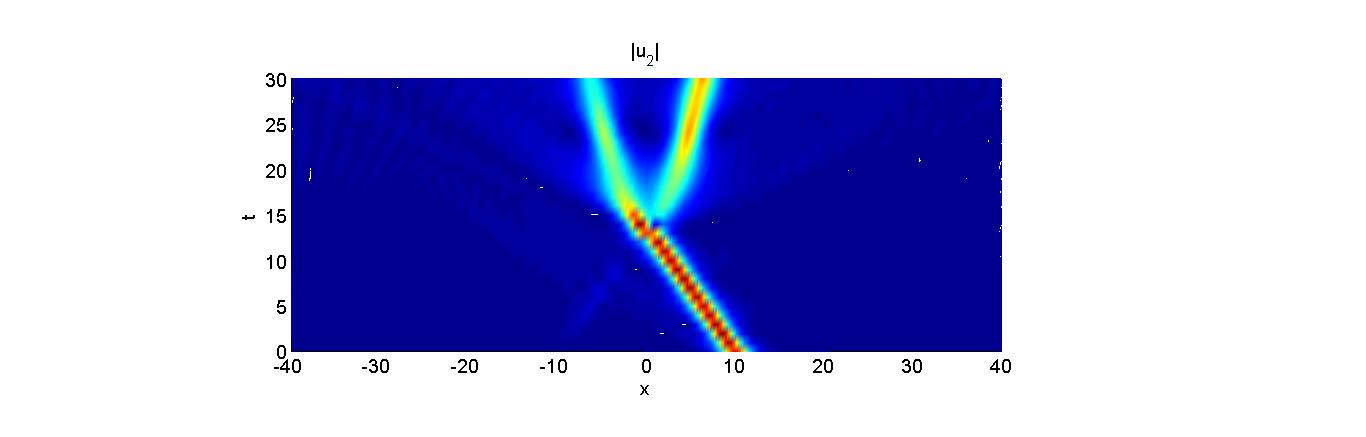}
\end{overpic}
\\
\hspace{-1cm}
\begin{overpic}[scale=0.4]{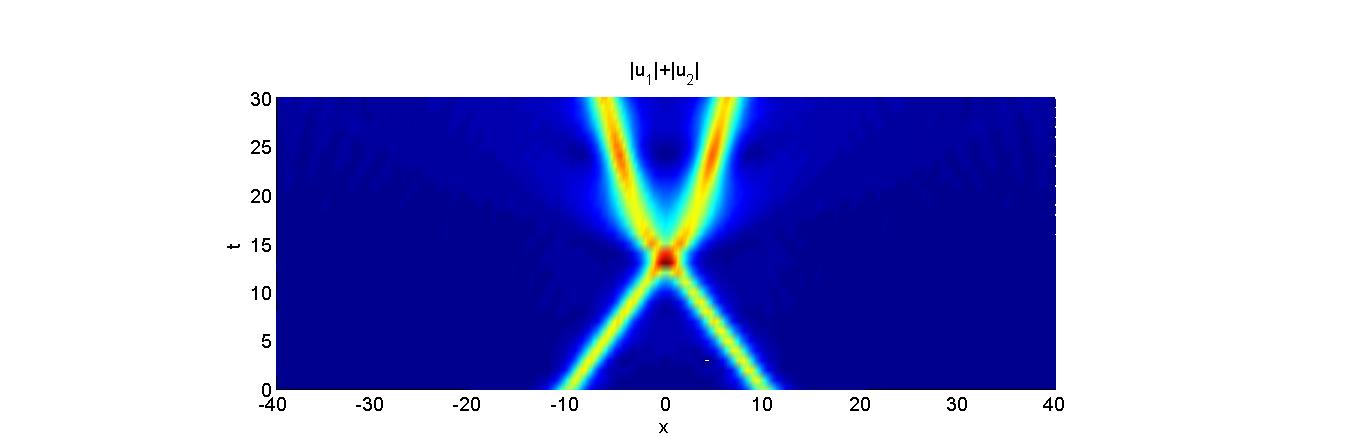}
\end{overpic}

\end{tabular}
\caption[]{{ The numerical simulation of the two soliton waves for Example \ref{ex12vng} with $\varpi_{1}= 0.0175$ and $\alpha = 1.8$.}}\label{fig:inelastic8h}
\end{center}
\end{figure}

\begin{exmp}\label{ex12vng} We consider the strongly coupled system as follows

\begin{equation}\label{sch1mn2}
\begin{split}
&i\frac{\partial u_{1}}{\partial t}- (-\Delta)^{\frac{\alpha}{2}}u_{1}+(|u_{1}|^{2}+|u_{2}|^{2})u_{1}+u_{1}+\varpi_{1}u_{2}=0,\\
&i\frac{\partial u_{2}}{\partial t}- (-\Delta)^{\frac{\alpha}{2}}u_{2}+(|u_{1}|^{2}+|u_{2}|^{2})u_{2}+\varpi_{1}u_{1}+u_{2}=0,\\
\end{split}
\end{equation}
 subject to the initial conditions

\begin{equation}\label{91}
\begin{split}
u_{1}(x,0)=\sqrt{2}r_{1}sech(r_{1}x+D)e^{iV_{0}x},\\
u_{2}(x,0)=\sqrt{2}r_{2}sech(r_{2}x+D)e^{iV_{0}x},\\
\end{split}
\end{equation}

\end{exmp}
where $r_{1} =r_{2} = 1$, $V_{0}=0.4$, $D=10$ and $x\in[-40,40]$. \\
Figures \ref{fig:elastic3} and \ref{fig:elastic4} shows that the proposed scheme simulates the solitary waves well. The two waves emerge without any changes
in their shapes for any $1 < \alpha \leq 2$. This phenomenon shows that the interaction is elastic.
The Figures \ref{fig:inelastic6}- \ref{fig:inelastic8h} present the numerical solutions for different values of order $\alpha$ for fixed $\varpi_{1}= 0.0175$.  From these figures it is obvious that the collision is always inelastic.  That is, the shapes and directions of two waves have changed after interaction.


\begin{figure}
\begin{center}
  \vspace{1mm}
\begin{tabular}{cc}
\hspace{-13cm}
\begin{overpic}[scale=0.25]{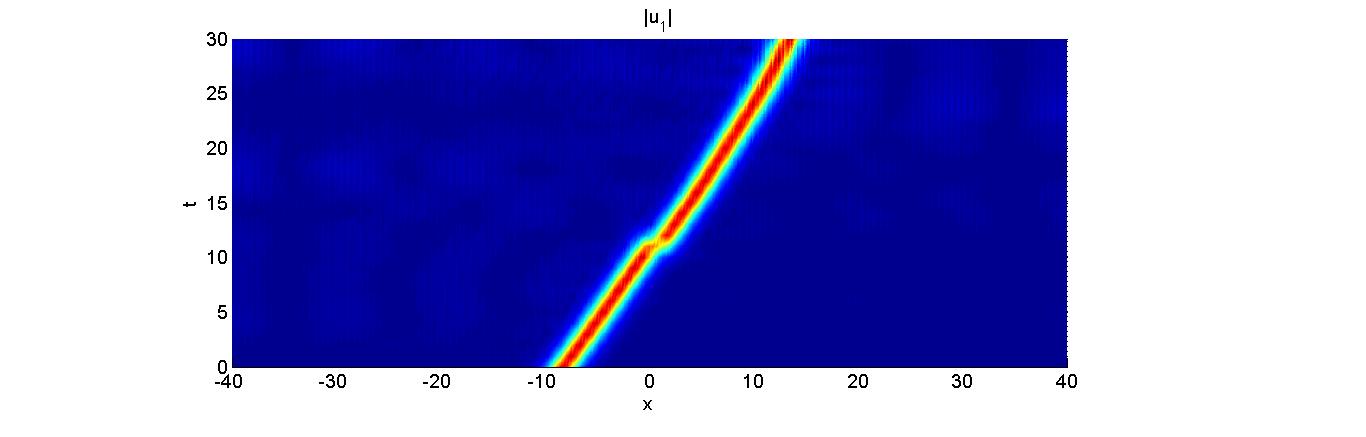}
\end{overpic}
 &\hspace{-12.2cm}
\begin{overpic}[scale=0.25]{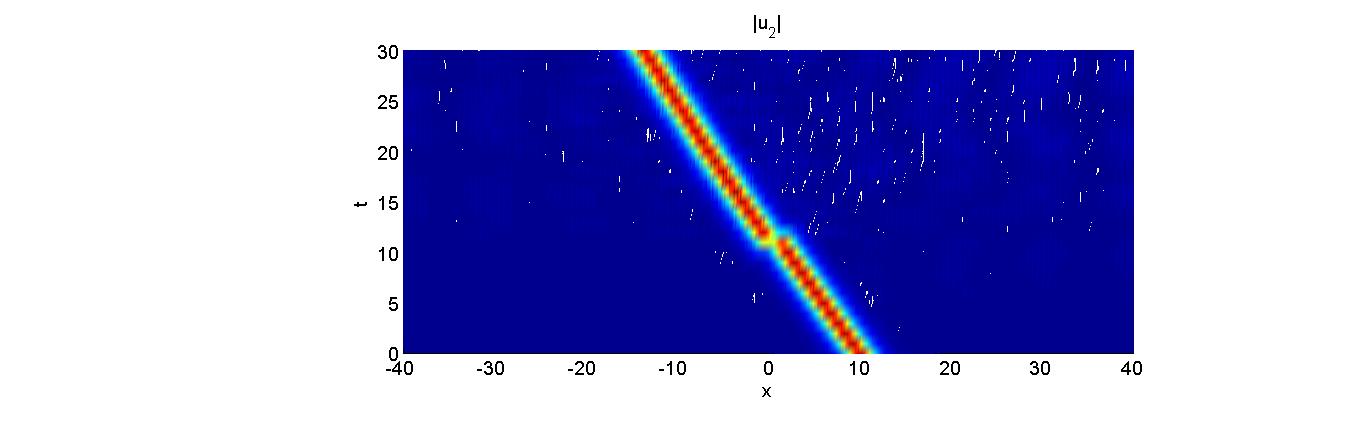}
\end{overpic}
\\
\\
\\
\hspace{-1cm}
\begin{overpic}[scale=0.4]{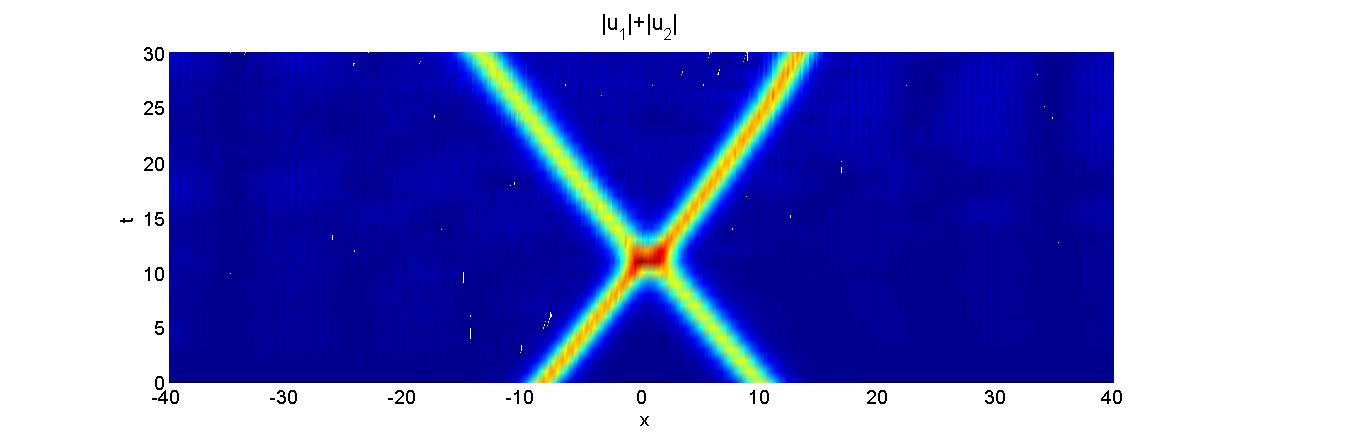}
\end{overpic}

\end{tabular}
\caption[]{{   The numerical simulation of the two soliton waves for Example \ref{ex12vng123}  with $\beta= 1$ and $\alpha = 2$.}}\label{fig:elastic1}
\end{center}
\end{figure}

\begin{figure}
\begin{center}
  \vspace{6mm}
\begin{tabular}{cc}
\hspace{-13cm}
\begin{overpic}[scale=0.25]{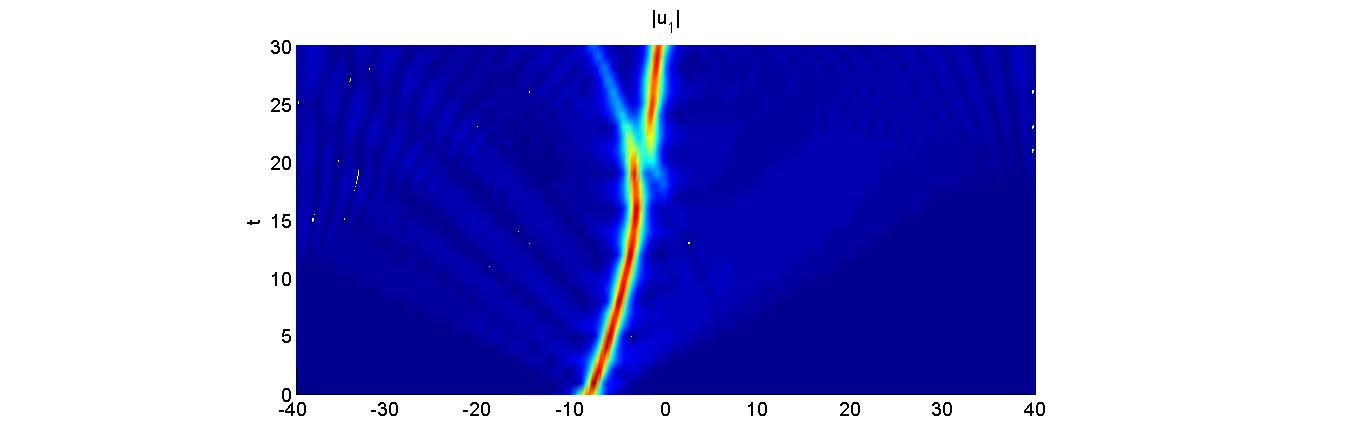}
\end{overpic}
 &\hspace{-12.5cm}
\begin{overpic}[scale=0.25]{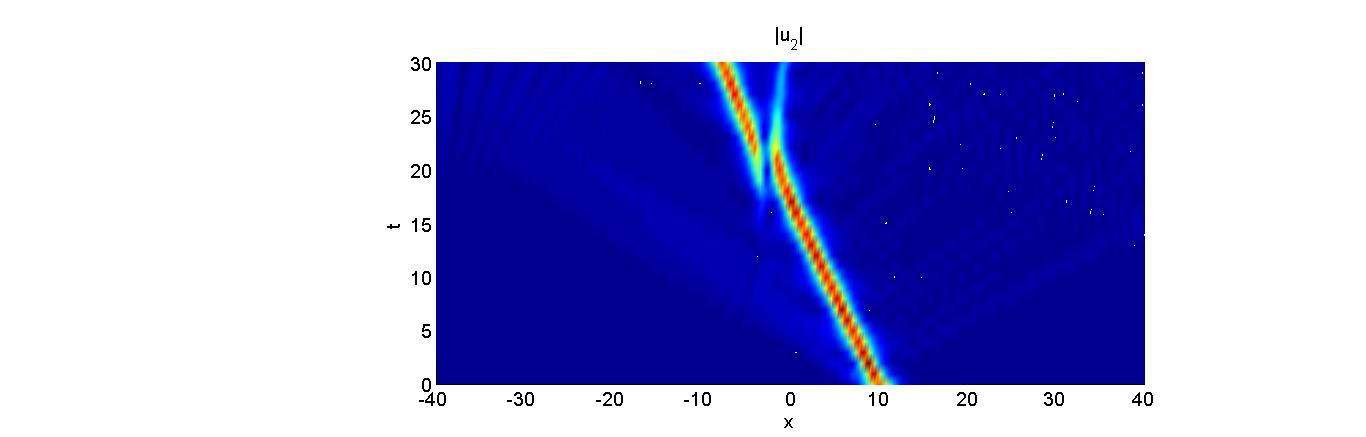}
\end{overpic}
\\
\\
\\
\hspace{-1cm}
\begin{overpic}[scale=0.4]{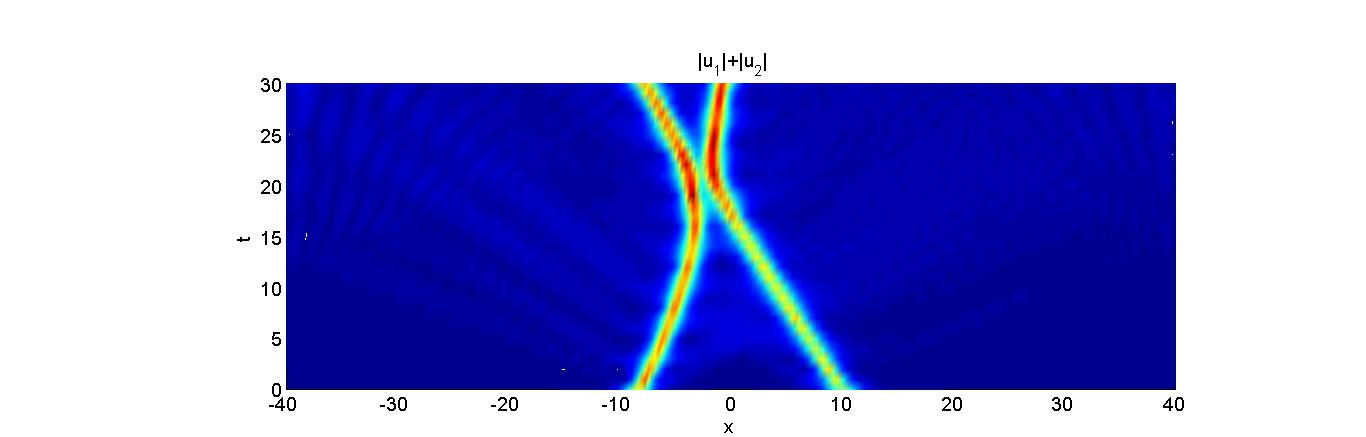}
\end{overpic}

\end{tabular}
\caption[]{{  The numerical simulation of the two soliton waves for Example \ref{ex12vng123}  with $\beta= 1$ and $\alpha = 1.6$.}}\label{fig:elastic2}
\end{center}
\end{figure}

\begin{figure}
\begin{center}
  \vspace{6mm}
\begin{tabular}{cc}
\hspace{-13cm}
\begin{overpic}[scale=0.25]{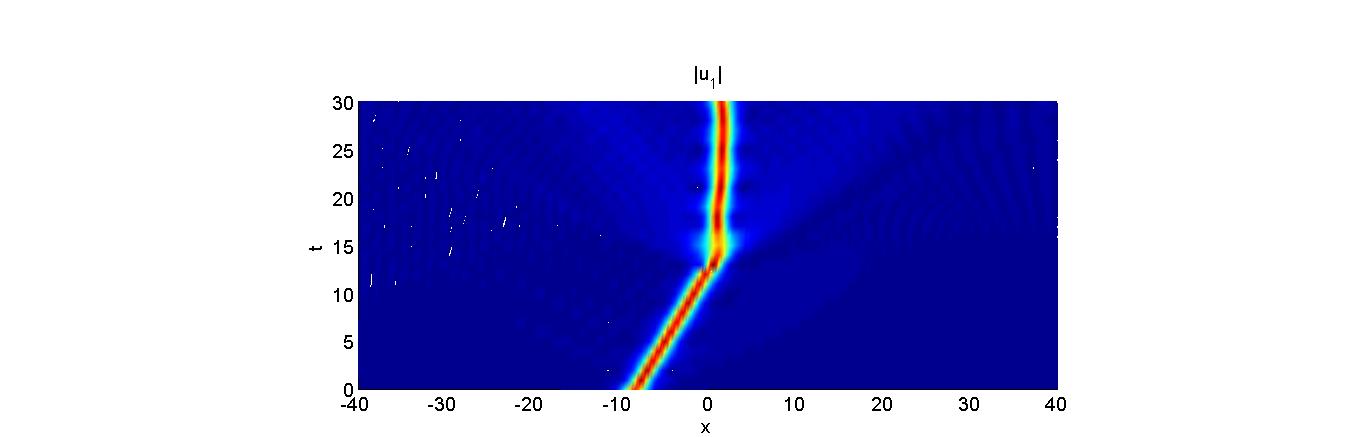}
\end{overpic}
 &\hspace{-12.0cm}
\begin{overpic}[scale=0.25]{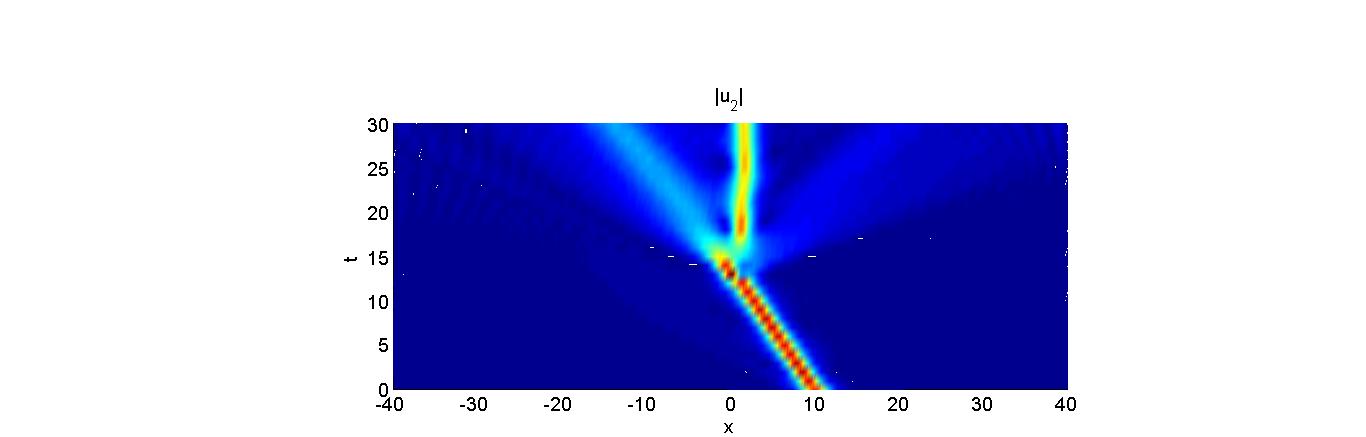}
\end{overpic}
\\
\\
\\
\hspace{-1cm}
\begin{overpic}[scale=0.4]{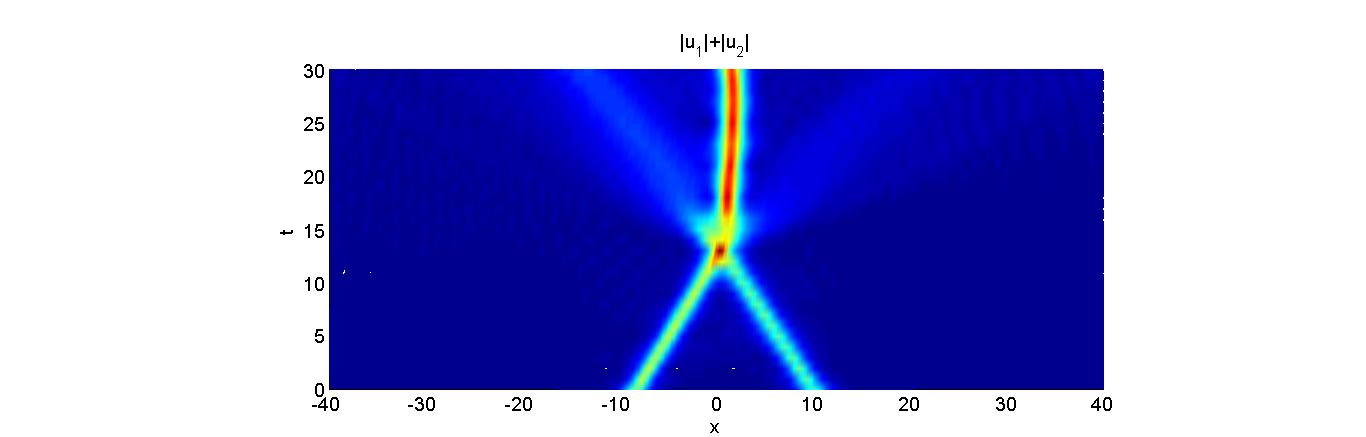}
\end{overpic}

\end{tabular}
\caption[]{{   The numerical simulation of the two soliton waves  for Example \ref{ex12vng123}  with $\beta= 0.3$ and $\alpha = 1.8$.}}\label{fig:elastic2nb}
\end{center}
\end{figure}

\begin{exmp}\label{ex12vng123} Finally, we consider the following weakly coupled problem
\begin{equation}\label{sch1mn2}
\begin{split}
&i\frac{\partial u_{1}}{\partial t}- (-\Delta)^{\frac{\alpha}{2}}u_{1}+(|u_{1}|^{2}+\beta|u_{2}|^{2})u_{1}=0,\\
&i\frac{\partial u_{2}}{\partial t}- (-\Delta)^{\frac{\alpha}{2}}u_{2}+(\beta|u_{1}|^{2}+|u_{2}|^{2})u_{2}=0,\\
\end{split}
\end{equation}
 subject to the initial conditions

\begin{equation}\label{91}
\begin{split}
u_{1}(x,0)=\sqrt{2}r_{1}sech(r_{1}x+D)e^{iV_{0}x},\\
u_{2}(x,0)=\sqrt{2}r_{2}sech(r_{2}x+D)e^{iV_{0}x},\\
\end{split}
\end{equation}
when $\beta= 1$ and $\alpha = 2$, the problem collapses to the Manakov equation, and
the solitary waves collide elastically see Figure \ref{fig:elastic1}. The exact solutions are given by
\begin{equation}\label{91}
\begin{split}
u_{1}(x,t)=\sqrt{2}r_{1}sech(r_{1}x-2r_{1}V_{0}t+D)e^{i(V_{0}x+(r_{1}^{2}-V_{0}^{2})t)},\\
u_{2}(x,t)=\sqrt{2}r_{2}sech(r_{2}x-2r_{2}V_{0}t-D)e^{i(-V_{0}x+(r_{2}^{2}-V_{0}^{2})t)},\\
\end{split}
\end{equation}
\end{exmp}
where $r_{1} = 1$, $r_{2} = 1$, $V_{0}=0.4$, $D=10$ and $x\in[-40,40]$. The Figures \ref{fig:elastic2} and \ref{fig:elastic2nb} present the numerical solutions for different values of order $\alpha$ and $\beta$.  From these figures  it is obvious that  the collision of solitons are inelastic. In particular, the colliding particles stick together after interaction when $\alpha =1.8$,
which means that there may occur a completely inelastic collision see Figure \ref{fig:elastic2nb}.

\newpage
\section{Conclusions}\label{s6n}
 We   propose  a  DDG  finite element method
for solving  fractional convection-diffusion and Schr\"{o}dinger type equations. The  scheme  is formulated  using  the direct weak  for these problems and the construct of the suitable numerical flux on the cell edges. Unlike  the  traditional  LDG  method,  the method  in this
paper  is applied  without  introducing  any  auxiliary  variables  or  rewriting the original
equation  into a  1st order  system.  An DDG method is proposed and stability and  a priori $L^{2}$ error estimates are presented. Numerical experiments for the fractional convection-diffusion and Schr\"{o}dinger type equations confirm the analysis. The numerical tests demonstrate both
accuracy and capacity of these methods, in particular, the numerical results are accurate for long time simulation.

$\mathbf{References}$
%

\end{document}